\documentclass[10pt]{amsart}

\usepackage{lipsum}
\usepackage{amsfonts,bm,amssymb,amsmath}
\usepackage{graphicx}
\usepackage{epstopdf}
\usepackage[caption=false]{subfig}
\usepackage{pgfplots}
\usepackage{algorithmic}
\usepackage{amsopn}
\usepackage{amsrefs}
\usepackage{float}
\usepackage{caption}

\usepackage{placeins}
\usepackage{multirow} 

\usepackage{hyperref}
\usepackage{cleveref}

\usepackage{color}
\usepackage{diagbox}

\usepackage{geometry}

\newtheorem{thm}{Theorem}[section]

\newtheorem{remark}{Remark}[section]
\newtheorem{corollary}{Corollary}[section]

\hypersetup{
	colorlinks,
	linkcolor={red!50!black},
	citecolor={blue!50!black},
	urlcolor={blue!80!black}
}

\numberwithin{equation}{section}

\definecolor{newcolor1}{rgb}{.8,.349,.1}
\colorlet{bblue}{blue!50!black}
\crefformat{equation}{(#2#1#3)}
\crefmultiformat{equation}{(#2#1#3)}{ and~(#2#1#3)}{, (#2#1#3)}{ and~(#2#1#3)}

\crefformat{figure}{Figure~#2#1#3}
\crefmultiformat{figure}{Figures~ #2#1#3}{ and~#2#1#3}{, (#2#1#3)}{ and~(#2#1#3)}
\crefformat{table}{Table~#2#1#3}

\def\a{\mbox{\boldmath $a$}}

\def\e{\mbox{\boldmath $e$}}
\def\f{\mbox{\boldmath $f$}}

\def\h{\mbox{\boldmath $h$}}
\def\m{\mbox{\boldmath $m$}}

\def\x{\mbox{\boldmath $x$}}
\def\y{\mbox{\boldmath $y$}}

\def\0{\mbox{\boldmath $0$}}

\begin{document}

\title[IHSS method for micromagnetics]{Efficient hermitian and skew-hermitian splitting methods for linear systems in micromagnetic simulations}

\author[Y. Miao]{Yingxi Miao}
\address{School of Mathematics and Physics\\ Xi'an-Jiaotong-Liverpool University\\Re'ai Rd. 111, Suzhou, 215123, Jiangsu\\ China.}
\email{yingximiao@gmail.com}


\author[C. Xie]{Changjian Xie$^*$}

\thanks{$*$ Corresponding author. Email: \texttt{Changjian.Xie@xjtlu.edu.cn}}
\address{School of Mathematics and Physics\\
Xi'an-Jiaotong-Liverpool University\\
Re'ai Rd. 111, Suzhou, 215123, Jiangsu\\
China.}
\email{Changjian.Xie@xjtlu.edu.cn}


\date{\today}


\begin{abstract}

For the Landau‑Lifshitz equation, the discrete linear systems obtained by our semi‑implicit method possess the following properties: they are large‑sparse systems with non‑Hermitian yet positive‑definite coefficient matrices. To solve these systems efficiently, we apply the Hermitian/skew-Hermitian splitting (HSS) method and its inexact variant (IHSS). Numerical experiments in one and three dimensions show that the spectral radius of the HSS iteration remains below its theoretical upper bound and strictly below one for the tested grid resolutions and damping parameters. Moreover, the theoretical bound closely follows the actual spectral radius, providing an accurate estimate of the convergence behavior. The IHSS results demonstrate effective convergence for the tested cases and show that its efficiency is sensitive to the splitting parameter. Overall, the two semi-implicit schemes exhibit comparable convergence behavior.

\end{abstract}

\maketitle

\section{Introduction}

Large sparse linear systems of the form $Ax=b$ constitute the core computational kernel in numerous scientific and engineering applications, including computational electromagnetism, fluid dynamics, finite element analysis, signal processing, and partial differential equation discretization \cite{gbikpi2022asynchronous,golub2013matrix}. 
Originated from micromagnetism, the Landau‑Lifshitz (LL) equation \cite{Landau1935On} describes the time evolution of magnetization fields for magnetic spin systems. This nonlinear PDE incorporates precessional rotation and damping effects, with a critical geometric constraint that preserves the magnitude of magnetization everywhere. Analytical solutions are only available for limited special cases, so high‑performance numerical schemes are essential for practical simulations. It serves as the core model for investigating domain‑wall motion, spin‑wave propagation and modern spintronic device behaviours.
This work considers semi‑implicit time discretization for the LL equation constructed by backward differentiation formulas (BDF) combined with one‑sided extrapolation, yielding a sequence of variable‑coefficient linear systems at each time step instead of solving fully nonlinear algebraic equations. The resulting linear operators possess the favourable property that all eigenvalues have strictly positive real parts, ensuring numerical invertibility. Nevertheless, large sparse matrix dimensions impose heavy computational overhead. Therefore, high‑performance iterative linear solvers are essential to achieve acceptable efficiency. Suitable solvers are required to cope with variable‑coefficient matrices across time steps, preserve computational stability, and accommodate the intrinsic unit‑length constraint of magnetization in micromagnetic simulations.

In practical numerical simulations, the coefficient matrix $A\in\mathbb{C}^{n\times n}$ is frequently non-Hermitian positive definite (NHPD), which means its Hermitian part is positive definite while the matrix itself lacks Hermitian symmetry \cite{bai2003hermitian}. Unlike Hermitian positive definite systems that support robust, well-established solvers with guaranteed convergence, NHPD linear systems pose significant computational challenges due to the coexistence of symmetric positive definite components and skew-symmetric asymmetric components. This structural complexity severely degrades the stability and efficiency of conventional iterative solvers, making the fast and accurate solution of NHPD linear systems a persistent hotspot in computational mathematics and numerical algebra \cite{yang2010generalized,bai2004preconditioned}.

Traditional iterative methods for linear systems exhibit distinct limitations when addressing NHPD problems. Classical Krylov subspace methods, such as the generalized minimal residual (GMRES) method and the biconjugate gradient stabilized (BiCGStab) method, are widely applicable to general non-Hermitian systems but suffer from obvious drawbacks \cite{saad2003iterative}. The GMRES method requires continuous storage of orthogonal basis vectors, leading to rapidly increasing memory consumption and computational overhead with the growth of matrix scale; restarted GMRES inevitably loses optimal convergence characteristics and may encounter slow convergence or stagnation for stiff NHPD systems \cite{saad1985generalized}. The BiCGStab method avoids excessive memory occupation but is prone to unstable iteration processes and residual oscillation when dealing with strongly non-Hermitian matrices. On the other hand, classic splitting iteration methods designed for Hermitian systems, including the Jacobi, Gauss–Seidel, and successive over-relaxation (SOR) methods, fail to maintain convergence stability for NHPD matrices, as their convergence criteria are strictly dependent on matrix Hermiticity and symmetric positive definiteness. These deficiencies restrict the application of traditional solvers in large-scale high-precision numerical simulations, creating an urgent demand for structured, efficient, and unconditionally convergent splitting iteration methods tailored for NHPD linear systems.

The Hermitian/skew-Hermitian splitting (HSS) iteration method, originally proposed by Bai et al. \cite{bai2003hermitian}, provides a targeted and robust solution for NHPD linear systems by fully exploiting the intrinsic structural properties of coefficient matrices. For any square matrix $A$, a unique structural decomposition $A=H+S$ can be performed, where $H=\frac{1}{2}(A+A^*)$ denotes the Hermitian part and $S=\frac{1}{2}(A-A^*)$ denotes the skew-Hermitian part \cite{bai2003hermitian}. For NHPD matrices, the Hermitian part $H$ is strictly positive definite, while the skew-Hermitian part $S$ captures all non-Hermitian characteristics of the matrix. Different from conventional single-structure splitting strategies, the HSS method constructs alternating iteration schemes based on the positive definite Hermitian component and the skew-Hermitian component, realizing separate and efficient processing of symmetric and asymmetric matrix properties \cite{bai2004preconditioned}. Theoretical studies have verified that the HSS iteration method possesses unconditional convergence for arbitrary NHPD linear systems with any initial iteration vector and positive relaxation parameter, which fundamentally overcomes the convergence instability of traditional splitting methods for non-Hermitian problems \cite{bai2007convergence}.

Benefiting from its superior structural adaptability and convergence robustness, the HSS method has been extensively developed and generalized in recent decades. To further reduce inner iteration computational cost, Bai et al. \cite{bai2003hermitian} proposed the inexact Hermitian/skew-Hermitian splitting (IHSS) method, which introduces Krylov subspace inner iterations to replace exact direct solutions in each outer iteration step, effectively improving computational efficiency for large-scale sparse systems \cite{bai2003hermitian}. Subsequently, a series of improved variants have been established, including modified HSS (MHSS) methods for complex symmetric systems, preconditioned HSS (PHSS) methods for accelerating convergence \cite{bai2004preconditioned,bai2011preconditioned,13}. Moreover, the HSS splitting strategy has been successfully extended to solve saddle-point problems, Sylvester matrix equations, and nonlinear equation systems with positive definite Jacobian matrices, demonstrating excellent universality and expandability. Recent studies have further optimized the HSS framework for parallel computing scenarios, developing asynchronous HSS iteration algorithms to adapt to large-scale distributed numerical simulations \cite{gbikpi2022asynchronous}.

Despite the comprehensive development of HSS-type methods, there still exist room for performance optimization and scenario expansion in practical applications. Most existing HSS variants rely on fixed parameter selection strategies, which cannot dynamically adapt to the spectral characteristics of different NHPD matrices, leading to suboptimal convergence rates for ill-conditioned systems \cite{bai2007convergence}. In addition, the computational efficiency of HSS methods for ultra-large sparse NHPD systems still needs to be further improved, and the theoretical analysis of convergence acceleration mechanisms remains incomplete. Against this background, this study focuses on the basic HSS iteration framework for NHPD linear systems, further analyzing its convergence characteristics, optimizing its iteration parameters and computational structure, and verifying its numerical superiority over traditional iterative methods through systematic theoretical derivation and numerical experiments. The research aims to provide a more efficient, stable, and adaptable numerical solution for large-scale NHPD linear systems in practical engineering and scientific computing problems.

The rest of the paper is organized in the following way. The physical model of micromagnetics and its numerical methods will be introduced in \Cref{sec: numerical scheme}. For the linear solver, we introduced the Hermitian/skew-Hermitian splitting method (HSS) in \Cref{sec:HSS}. For the numerical experiments in \Cref{sec:experiments}, we first test the linear system with non-Hermitian, sparsity and positivity. The spectral radius is tested in \Cref{sec:radii}. The IHSS results for scheme A and B are presented in \Cref{sec:ihss}. We drae our conclusion in \Cref{sec:conclusions}.

\section{The physical model and the numerical scheme}
\label{sec: numerical scheme}

\subsection{Governing equation}

The Landau-Lifshitz-Gilbert (LLG) equation forms the fundamental basis of micromagnetics, providing a rigorous description of the spatiotemporal evolution of magnetization in ferromagnetic materials by incorporating two key physical phenomena: gyromagnetic precession and dissipative relaxation \cite{Landau1935On,Brown1963micromagnetics}. In nondimensional form, this governing equation is expressed as
\begin{align}\label{c1-large}
{\m}_t =-{\m}\times{\bm h}_{\text{eff}}-\alpha{\m}\times({\m}\times{\bm h}_{\text{eff}})
\end{align}
subject to the homogeneous Neumann boundary condition
\begin{equation}\label{boundary-large}
\frac{\partial{\m}}{\partial {\bm \nu}}\Big|_{\partial \Omega}=0,
\end{equation}
where \(\Omega \subset \mathbb{R}^d\) (\(d=1,2,3\)) represents the bounded domain of the ferromagnetic material, and \(\bm \nu\) is the unit outward normal vector on the boundary \(\partial \Omega\). This boundary condition ensures no magnetic surface charge, a physically appropriate assumption for isolated ferromagnetic systems.

The magnetization field \(\m: \Omega \to \mathbb{R}^3\) is a three-dimensional vector field satisfying the pointwise constraint \(|\m| = 1\),  stemming from the quantum mechanical alignment of electron spins in ferromagnets. The first term on the right-hand side of \cref{c1-large} describes gyromagnetic precession, where magnetic moments precess around the effective field \(\bm h_{\text{eff}}\). The second term represents dissipative relaxation, with \(\alpha > 0\) being the dimensionless Gilbert damping coefficient that governs the rate of energy dissipation into the lattice.

From the perspective of the Gibbs free energy functional, the effective field \(\bm h_{\text{eff}}\) is obtained as the functional derivative of the Gibbs free energy \(F[\m]\) with respect to the magnetization, i.e., \(\bm h_{\text{eff}} = -\delta F/\delta \m\). his functional incorporates all relevant energy contributions in ferromagnetic systems—exchange, anisotropy, magnetostatic (stray), and Zeeman energies—and is given by
\begin{equation}\label{LL-Energy}
F[\m] = \frac {\mu_0 M_s^2}{2} \left\{\int_\Omega \left( \epsilon|\nabla\m|^2 +
q\left(m_2^2 + m_3^2\right)
-2\h_e\cdot\m - \h_s\cdot\m \right)\mathrm{d}\x \right\} . 
\end{equation}
Here, \(\mu_0 = 4\pi \times 10^{-7}\, \text{H/m}\) is the vacuum permeability, \(M_s\) is the saturation magnetization, and \(\epsilon\) and \(q\) are dimensionless parameters defined subsequently. The vectors \(\h_e\) and \(\h_s\) denote the externally applied magnetic field and the stray field, respectively. For uniaxial ferromagnetic materials with a single easy axis, the effective field \(\bm h_{\text{eff}}\) decomposes into distinct physical components, yielding the explicit expression
\begin{align}
{\bm h}_{\text{eff}} =\epsilon\Delta\m-q(m_2\e_2+m_3\e_3)+\h_s+\h_e,
\end{align}
where \(\epsilon = C_{\text{ex}}/(\mu_0 M_s^2 L^2)\) and \(q = K_u/(\mu_0 M_s^2)\). Here, \(L\) is the characteristic length scale, \(C_{\text{ex}}\) is the exchange constant (controlling short-range spin alignment), and \(K_u\) is the uniaxial anisotropy constant (representing the energy cost for deviation from the easy axis). The unit vectors \(\e_2 = (0,1,0)\) and \(\e_3 = (0,0,1)\) efine the hard axes perpendicular to the uniaxial easy axis, and \(\Delta\) is the Laplacian operator in \(d\)-dimensional space.
For Permalloy (NiFe), a common soft ferromagnetic material in spintronics, standard parameter values from the literature are: \(C_{\text{ex}} = 1.3 \times 10^{-11}\, \text{J/m}\), \(K_u = 100\, \text{J/m}^3\), and \(M_s = 8.0 \times 10^5\, \text{A/m}\). The stray field \(\h_s\) originates from magnetic charge distributions at domain boundaries and material surfaces and is mathematically represented by the integral expression
\begin{align}\label{eqn:div}
{\h}_{\text{s}}=\frac{1}{4\pi}\nabla \int_{\Omega} \nabla\left( \frac{1}{|\x-\y|}\right)\cdot {\bm m}({\bm y})\,d{\bm y},
\end{align}
which is a formulation that exhibits long-range spatial correlations. A critical computational advancement for practical micromagnetic simulations is that for rectangular domains \(\Omega\), the evaluation of \(\h_s\) can be efficiently computed via the Fast Fourier Transform (FFT) \cite{Wang2000}, which reduces the asymptotic computational complexity from \(O(N^d)\) to \(O(N^d \log N)\) for \(d\)-dimensional grids, enabling large-scale simulations.

To facilitate numerical discretization, we introduce the composite source term
\begin{align}\label{eq-4}
\f=-Q(m_2\e_2+m_3\e_3)+\h_s+\h_e.
\end{align}
which aggregates the anisotropy, stray field, and external field contributions. Substituting this source term into \cref{c1-large}, the LLG equation is re-expressed as
\begin{align}\label{eq-5}
\m_t=-\m\times(\epsilon\Delta\m+\f)-\alpha\m\times\m\times(\epsilon\Delta\m+\f).
\end{align}
Leveraging the vector triple product identity \(\a \times ({\bm b} \times {\bm c}) = (\a \cdot {\bm c}){\bm b} - (\a \cdot {\bm b}){\bm c}\) and the pointwise constraint \(|\m| = 1\) (which implies \(\m \cdot \partial_t \m = 0\) via time differentiation), we simplify \cref{eq-5} to an equivalent formulation that is more amenable to stable numerical discretization:
\begin{equation}\label{eq-model}
\m_t=\alpha  (\epsilon\Delta\m+\f)+\alpha \left(\epsilon |\nabla \m|^2 -\m \cdot\f \right)\m-\m\times(\epsilon\Delta\m+\f).
\end{equation}
We establish a standardized discretization framework to underpin subsequent numerical approximations, defining temporal and spatial discretization notations and boundary condition enforcement strategies.

Let \(k = t^{n+1} - t^n\) denote the uniform temporal step-size, with discrete time levels given by \(t^n = nk\) for \(n = 0,1,\dots,N_T\), where \(N_T = \lfloor T/k \rfloor\) and \(T\) is the final simulation time. For spatial discretization, a uniform Cartesian grid is employed with mesh-size \(h_x = h_y = h_z = h = L/N\) (for cubic domains of characteristic length \(L\)). The notation \(\m_{i,j,\ell}^n\) denotes the numerical approximation of \(\m(x_{i-1/2}, y_{j-1/2}, z_{\ell-1/2}, t^n)\), where \(x_{i-1/2} = (i - 1/2)h_x\), \(y_{j-1/2} = (j - 1/2)h_y\), and \(z_{\ell-1/2} = (\ell - 1/2)h_z\) define cell-centered grid positions. The index range \(-1 \leq i,j,\ell \leq N+2\) is adopted to accommodate boundary extrapolation.
To enforce the homogeneous Neumann boundary condition \cref{boundary-large} while preserving high-order accuracy, a third-order extrapolation scheme is implemented. For instance, along the \(z\)-direction (normal to the boundary at \(z=0\) and \(z=L\)), the extrapolation rules for the magnetization field are:
\begin{align*}
\m_{i,j,1}=\m_{i,j,0},\quad \m_{i,j,-1}=\m_{i,j,2},\quad \m_{i,j,N+1}=\m_{i,j,N} ,\quad \m_{i,j,N+2}=\m_{i,j,N-1}.
\end{align*}
Analogous extrapolation formulas are applied to enforce the boundary condition along the \(x\)- and \(y\)-directions, ensuring consistent high-order accuracy across the entire computational domain.

The fourth-order spatial difference operators should be introduced further. 
To achieve fourth-order spatial accuracy, which is essential for resolving fine magnetic structures (e.g., domain walls with nanoscale width) and minimizing discretization-induced errors, we employ long-stencil finite difference operators for first and second partial derivatives. For the \(x\)-direction, the fourth-order accurate operators for \(\partial_x\) (denoted \({\mathcal D}_{x,(4)}^1\)) and \(\partial_x^2\) (denoted \({\mathcal D}_{x,(4)}^2\)) are defined as:
\begin{eqnarray} 
\hspace{-0.35in}  
{\mathcal D}_{x,(4)}^1 f_{i,j,k} &=& \tilde{D}_x ( 1 - \frac{h_x^2}{6} D_x^2 ) f_{i,j,k} \nonumber 
\\
&=& 
\frac{  f_{i-2,j,k} - 8 f_{i-1,j,k}  + 8 f_{i+1,j,k} - f_{i+2,j,k} }{12 h_x} ,  
\label{FD-4th-1} 
\\
\hspace{-0.35in}  
{\mathcal D}_{x,(4)}^2 f_{i,j,k} &=& D_x^2 ( 1 - \frac{h_x^2}{12} D_x^2 ) f_{i,j,k}  \nonumber 
\\
&=& 
\frac{ - f_{i-2,j,k} + 16 f_{i-1,j,k,k} - 30 f_{i,j,k} + 16 f_{i+1,j,k} - f_{i+2,j,k} }{12 h_x^2 } . 
\label{FD-4th-2} 
\end{eqnarray} 
These operators are derived by canceling leading-order discretization errors via inclusion of adjacent grid points, resulting in fourth-order convergence (\(O(h^4)\)) for sufficiently smooth functions.

By symmetry, the fourth-order difference operators for the \(y\)- and \(z\)-directions—designated \({\mathcal D}_{y,(4)}^1\), \({\mathcal D}_{y,(4)}^2\), \({\mathcal D}_{z,(4)}^1\), and \({\mathcal D}_{z,(4)}^2\)—are defined by substituting the respective spatial coordinates and mesh-sizes. The discrete Laplacian operator \(\Delta_h\), which approximates the continuous Laplacian \(\Delta\) with fourth-order accuracy, is constructed as the sum of the second-order operators in all three directions: \(\Delta_h = {\mathcal D}_{x,(4)}^2 + {\mathcal D}_{y,(4)}^2 + {\mathcal D}_{z,(4)}^2\).

To provide a rigorous performance comparison for two numerical methods, we adapt the third-order (BDF3) backward differentiation formula with one-sided extrapolation, which are standard in higher-order numerical schemes, to include the fourth-order spatial operators and third-order temporal discretizations for two equivalent LLG formulation \cref{c1-large} and \cref{eq-model}. We develop the efficient linear solver for these fully discretized linear system of equations.

\subsection{Scheme A}
A third-order semi-implicit scheme for model \cref{eq-model} as presented in \cite{xie2025thirdorder}, is given by
\begin{equation}\label{scheme-third-order}
\left\{ 
\begin{aligned}
&\frac{\frac{11}{6}\tilde{\m}_h^{n+3}-3{\m}_h^{n+2}+\frac32{\m}_h^{n+1}-\frac13 {\m}_h^n}{k}\\
&\quad=-\hat{\m}_h^{n+3} \times (\epsilon \Delta_{h,(4)}\tilde{\m}_h^{n+3}+\hat{\f}_h^{n+3})+\alpha(\epsilon \Delta_{h,(4)}\tilde{\m}_h^{n+3} +\hat{\f}_h^{n+3})\\
&\qquad+\alpha(\epsilon |\tilde{\nabla}_{h,(4)}\hat{\m}_h^{n+3}|^2+\hat{\m}_h^{n+3}\cdot \hat{\f}_h^{n+3}) \hat{\m}_h^{n+3}, \\
&\m_h^{n+3}=\frac{\tilde{\m}_h^{n+3}}{|\tilde{\m}_h^{n+3}|},
\end{aligned}
\right.
\end{equation}
where
\begin{align*}
\hat{\m}_h^{n+3} &= 3 \m_h^{n+2}-3\m_h^{n+1} + \m_h^n,\\
\hat{\f}_h^{n+3} &= 3 \f_h^{n+2}-3\f_h^{n+1} + \f_h^n.
\end{align*}

The third order semi-implicit method \cref{scheme-third-order} gives the linear system of equations
\begin{align*}
	&\left(\frac{11}{6}I+k\epsilon \hat{\m}_h^{n+3} \times \Delta_h-k\alpha \epsilon \Delta_h\right)\tilde{{\m}}_h^{n+3}=3{\m}_h^{n+2}-\frac32{\m}_h^{n+1}+\frac13 {\m}_h^n\\
	&-k\hat{\m}_h^{n+3} \times \hat{\f}_h^{n+3} +\alpha \hat{\f}_h^{n+3}+\alpha k(\epsilon |\tilde{\nabla}_{h,(4)}\hat{\m}_h^{n+3}|^2+\hat{\m}_h^{n+3}\cdot \hat{\f}_h^{n+3}) \hat{\m}_h^{n+3}
\end{align*}
Thus, 
\begin{align*}
A&=\frac{11}{6}I+k\epsilon \hat{\m}_h^{n+3} \times \Delta_h-k\alpha \epsilon \Delta_h\\
&=\begin{pmatrix}
\frac{11}{6}-k\alpha \epsilon \Delta_h &-k\epsilon\hat{w}\Delta_h & k\epsilon\hat{v} \Delta_h\\
k\epsilon\hat{w}\Delta_h &\frac{11}{6}-k\alpha \epsilon \Delta_h& -k\epsilon\hat{u} \Delta_h\\
-k\epsilon\hat{v}\Delta_h & k\epsilon\hat{u}\Delta_h &\frac{11}{6}-k\alpha \epsilon \Delta_h
\end{pmatrix}
\end{align*}

Here, we can simply denote
\begin{align*}
	A=H_1+S_2,
\end{align*}
where $H_1$ is a symmetric matrix and $S_1$ is skew-symmetric matrix. Specifically,
\begin{align*}
	H_1&=\frac{11}{6}I-k\alpha \epsilon \Delta_h=\begin{pmatrix}
	\frac{11}{6}-k\alpha \epsilon \Delta_h &0& 0\\
	0&\frac{11}{6}-k\alpha \epsilon \Delta_h& 0 \\
	0 & 0 &\frac{11}{6}-k\alpha \epsilon \Delta_h
	\end{pmatrix}\\
	S_1&=k\epsilon \hat{\m}_h^{n+3}\times \Delta_h=\begin{pmatrix}
	0 &-k\epsilon\hat{w}\Delta_h & k\epsilon\hat{v} \Delta_h\\
	k\epsilon\hat{w}\Delta_h &0 & -k\epsilon\hat{u} \Delta_h\\
	-k\epsilon\hat{v}\Delta_h & k\epsilon\hat{u}\Delta_h &0
	\end{pmatrix}
\end{align*}

\subsection{Scheme B} \label{discretisations}

The SIPM in \cite{Xie2018} employs a semi-implicit treatment for both the gyromagnetic and damping terms, meaning $\Delta \m$ is handled implicitly while coefficient functions are updated via a second-order explicit extrapolation. A natural extension is to apply a third-order BDF method, where $\Delta \m$ remains implicit and coefficient functions are extrapolated using a third-order accurate explicit formula. This concept leads to the proposed method \cite{xie2025novel}:

\begin{equation}\label{proposed}
\left\{ 
\begin{aligned}
&\frac{\frac{11}{6} \tilde{\m}_h^{n+3} - 3 {\m}_h^{n+2} + \frac32 {\m}_h^{n+1}-\frac13 {\m}_h^n}{k}
=  - \hat{\m}_h^{n+3} \times \left(\epsilon \Delta_{h,(4)} \tilde{{\m}}_h^{n+3} +\hat{\f}_h^{n+3}\right) \\
&\quad -\alpha\hat{\m}_h^{n+3}\times (\hat{\m}_h^{n+3}\times (\epsilon\Delta_{h,(4)}\tilde{{\m}}_h^{n+3}+\hat{\f}_h^{n+3}))\\ 
&\m_h^{n+3} = \frac{\tilde{\m}_h^{n+3}}{ |\tilde{\m}_h^{n+3}| } ,
\end{aligned}
\right.
\end{equation}
where
\begin{align*}
\hat{\m}_h^{n+3} &= 3 \m_h^{n+2} - 3\m_h^{n+1}+\m_h^{n},\\
\hat{\f}_h^{n+3} &= 3 \f_h^{n+2} - 3\f_h^{n+1}+{\f}_h^n.
\end{align*}

The proposed method \cref{proposed} gives the linear system of equations
\begin{align*}
	&\left(\frac{11}{6}I+k\epsilon \hat{\m}_h^{n+3}\times \Delta_h+\alpha \epsilon k\hat{\m}_h^{n+3}\times (\hat{\m}_h^{n+3}\times \Delta_h)\right)\tilde{{\m}}_h^{n+3}=3 {\m}_h^{n+2} - \frac32 {\m}_h^{n+1}+\frac13 {\m}_h^n\\
	&-k\hat{\m}_h^{n+3}\times \hat{\f}_h^{n+3}-k\alpha \hat{\m}_h^{n+3}\times (\hat{\m}_h^{n+3}\times \hat{\f}_h^{n+3}).
\end{align*}

Here, we notice that
\begin{align*}
\hat{\m}\times \Delta_h \tilde{{\m}}=\begin{pmatrix}
0 &-\hat{w}\Delta_h & \hat{v} \Delta_h\\
\hat{w}\Delta_h &0 & -\hat{u} \Delta_h\\
-\hat{v}\Delta_h & \hat{u}\Delta_h &0
\end{pmatrix}\begin{pmatrix}
\tilde{u}\\
\tilde{v}\\
\tilde{w}
\end{pmatrix}
\end{align*}
and 
\begin{align*}
\hat{\m}\times(\hat{\m}\times \Delta_h \tilde{{\m}})&=\begin{pmatrix}
0 &-\hat{w} & \hat{v} \\
\hat{w} &0 & -\hat{u} \\
-\hat{v} & \hat{u} &0
\end{pmatrix}\begin{pmatrix}
0 &-\hat{w}\Delta_h & \hat{v} \Delta_h\\
\hat{w}\Delta_h &0 & -\hat{u} \Delta_h\\
-\hat{v}\Delta_h & \hat{u}\Delta_h &0
\end{pmatrix}\begin{pmatrix}
\tilde{u}\\
\tilde{v}\\
\tilde{w}
\end{pmatrix}\\
&=\begin{pmatrix}
-(\hat{v}^2+\hat{w}^2)\Delta_h& \hat{u}\hat{v}\Delta_h & \hat{u}\hat{w}\Delta_h\\
\hat{u}\hat{v}\Delta_h & -(\hat{u}^2+\hat{w}^2)\Delta_h & \hat{u}\hat{v} \Delta_h\\
\hat{u}\hat{w}\Delta_h & \hat{v}\hat{w}\Delta_h & -(\hat{u}^2+\hat{v}^2)\Delta_h
\end{pmatrix}\begin{pmatrix}
\tilde{u}\\
\tilde{v}\\
\tilde{w}
\end{pmatrix}
\end{align*}
Thus,
\begin{align*}
	B&=\frac{11}{6}I+k\epsilon \hat{\m}_h^{n+3}\times \Delta_h+\alpha \epsilon k\hat{\m}_h^{n+3}\times (\hat{\m}_h^{n+3}\times \Delta_h)\\
	&=\begin{pmatrix}
	\frac{11}{6}-k\alpha \epsilon(\hat{v}^2+\hat{w}^2)\Delta_h& -k\epsilon \hat{w}\Delta_h+ k\alpha \epsilon\hat{u}\hat{v}\Delta_h &k\epsilon \hat{v}\Delta_h+ k\alpha \epsilon\hat{u}\hat{w}\Delta_h\\
	k\epsilon\hat{w}\Delta_h+ k\alpha \epsilon\hat{u}\hat{v}\Delta_h &\frac{11}{6} -k\alpha \epsilon(\hat{u}^2+\hat{w}^2)\Delta_h &-k\epsilon \hat{u}\Delta_h+ k\alpha \epsilon\hat{u}\hat{v} \Delta_h\\
	-k\epsilon \hat{v}\Delta_h+ k\alpha \epsilon\hat{u}\hat{w}\Delta_h &k\epsilon\hat{u}\Delta_h+ k\alpha \epsilon \hat{v}\hat{w}\Delta_h &\frac{11}{6} -k\alpha \epsilon(\hat{u}^2+\hat{v}^2)\Delta_h
	\end{pmatrix}
\end{align*}

Here, we simply denote 
\begin{align*}
	B=H_2+S_2,
\end{align*}
where $H_1$ is a symmetric matrix and $S_2$ is a skew-symmetric matrix. Specifically,
\begin{align*}
	H_2&=\begin{pmatrix}
	\frac{11}{6}-k\alpha \epsilon(\hat{v}^2+\hat{w}^2)\Delta_h& k\alpha \epsilon\hat{u}\hat{v}\Delta_h & k\alpha \epsilon\hat{u}\hat{w}\Delta_h\\
 k\alpha \epsilon\hat{u}\hat{v}\Delta_h &\frac{11}{6} -k\alpha \epsilon(\hat{u}^2+\hat{w}^2)\Delta_h &k\alpha \epsilon\hat{u}\hat{v} \Delta_h \\
	k\alpha \epsilon\hat{u}\hat{w}\Delta_h &k\alpha \epsilon \hat{v}\hat{w}\Delta_h &\frac{11}{6} -k\alpha \epsilon(\hat{u}^2+\hat{v}^2)\Delta_h
	\end{pmatrix}\\
		S_2&=\begin{pmatrix}
0& -k\epsilon \hat{w}\Delta_h&k\epsilon \hat{v}\Delta_h\\
	k\epsilon\hat{w}\Delta_h&0&-k\epsilon \hat{u}\Delta_h\\
	-k\epsilon \hat{v}\Delta_h &k\epsilon\hat{u}\Delta_h &0
	\end{pmatrix}
\end{align*}


\begin{remark}
Mathematically, while \cref{eq-5} and \cref{eq-model} are equivalent under the normalization condition $|\m|=1$, their numerical implementations in schemes \cref{scheme-third-order} and \cref{proposed}  yield different performances. The differences between the schemes are as follows.
\begin{itemize}
	\item The damping terms differ because $|\hat{\m}^{n+3}|\neq 1$ and $\hat{\m}^{m+3}\cdot\Delta \tilde{{\m}}_h^{n+3}\neq |\nabla \hat{\m}^{n+3}|$;
	\item The coefficient matrices of the resulting linear systems differ: scheme \cref{scheme-third-order} yields $\frac{11}{6}I+\epsilon \hat{\m}^{n+3}\times \Delta_h+\alpha \epsilon \Delta_h$, whereas scheme \cref{proposed} yields $\frac{11}{6}I+\epsilon \hat{\m}^{n+3}\times \Delta_h+\alpha \hat{\m}^{n+3}\times (\hat{\m}^{n+3}\times \Delta_h)$. The right-hand side source terms also differ.
\end{itemize}
\end{remark}

\section{Linear Solver}
In this section, we will combine the matrix structure and sparsity of the linear system and develop a matrix splitting method for its convergence and efficiency.

\subsection{Hermitian/skew-Hermitian splitting (HSS) method}\label{sec:HSS}

We can proceed a Hermitian/skew-Hermitian splitting (HSS) by \cite{bai2003hermitian,bai2004preconditioned}
\begin{equation}
B = H + S
\end{equation}
where
\begin{equation}
H = \frac{1}{2}(B + B^*) \quad \text{and} \quad S = \frac{1}{2}(B - B^*)
\end{equation}

Note that the eigenvalues of $S$ are either $0$ or imaginary number and of $H$ are all real parts. We have known the real parts of complex eigenvalues for $B$ are positive. Thus, the eigenvalues of Hermitian $H$ are positive which means the matrix $B$ is positive definite, and non-Hermitian .

Consequently, we can list the feature of the matrix $B$:
\begin{enumerate}
	\item Large sparse;
	\item Non-Hermitian: $B \neq B^*$;
	\item Positive definite (or real positive): The real part $\mathcal{R}(\boldsymbol{x}^* B \boldsymbol{x}) > 0$ for any nonzero $\boldsymbol{x} \in \mathbb{C}^n$.
\end{enumerate}
The HSS iteration method. Given an initial guess $\boldsymbol{x}^{(0)}$. For $\ell = 0, 1, \dots$ until the constructed sequences $\{\boldsymbol{x}^{(\ell)}\}$ converges, compute
\[
\begin{cases}
(\beta I + H)\boldsymbol{x}^{(\ell+1/2)} = (\beta I - S)\boldsymbol{x}^{(\ell)} + \boldsymbol{b}, \\
\quad (\beta I + S)\boldsymbol{x}^{(\ell+1)} = (\beta I - H)\boldsymbol{x}^{(\ell+1/2)} + \boldsymbol{b},
\end{cases}
\]

\begin{thm}[Convergence analysis of HSS] \label{thm:HSS-convergence}
	Let $B \in \mathbb{C}^{n \times n}$ be a positive definite matrix, let $H = \frac{1}{2}(B+B^*)$ and $S = \frac{1}{2}(B-B^*)$ be its Hermitian and skew-Hermitian parts, and let $\beta$ be a positive constant. Then the iteration matrix $M(\beta)$ of the HSS iteration is given by
	\begin{equation}
	M(\beta) = (\beta I + S)^{-1}(\beta I - H)(\beta I + H)^{-1}(\beta I - S),
	\end{equation}
	and its spectral radius $\rho(M(\beta))$ is bounded by
	\begin{equation}
	\sigma(\beta) = \max_{\lambda_i \in \lambda(H)} \left| \frac{\beta - \lambda_i}{\beta + \lambda_i} \right|,
	\end{equation}
	where $\lambda(H)$ is the spectral set of the matrix $H$. Therefore, it holds that
	\begin{equation}
	\rho(M(\beta)) \leq \sigma(\beta) < 1 \quad \forall \beta > 0;
	\end{equation}
	i.e., the HSS iteration converges to the unique solution $x^* \in \mathbb{C}^n$ of the system of linear equations
\end{thm}

\begin{corollary} \label{cor:HSS-opt-beta}
	Let $B \in \mathbb{C}^{n \times n}$ be a positive definite matrix, let $H = \frac{1}{2}(B+B^*)$ and $S = \frac{1}{2}(B-B^*)$ be its Hermitian and skew-Hermitian parts, and let $\beta$ be a positive constant and let $\gamma_{\min}$ and $\gamma_{\max}$ be the minimum and the maximum eigenvalues of the matrix $H$, respectively. Then
	\begin{equation}
	\beta^* = \arg\min_{\beta} \left\{ \max_{\gamma_{\min} \leq \lambda \leq \gamma_{\max}} \left| \frac{\beta - \lambda}{\beta + \lambda} \right| \right\} = \sqrt{\gamma_{\min}\gamma_{\max}},
	\end{equation}
	and
	\begin{equation}
	\sigma(\beta^*) = \frac{\sqrt{\gamma_{\max}} - \sqrt{\gamma_{\min}}}{\sqrt{\gamma_{\max}} + \sqrt{\gamma_{\min}}} = \frac{\sqrt{\kappa(H)} - 1}{\sqrt{\kappa(H)} + 1},
	\end{equation}
	where $\kappa(H)$ is the spectral condition number of $H$.
\end{corollary}

We define the optimal parameter:
\begin{align*}
	\tilde{\beta} = \arg\min_{\beta} \{\rho(M(\beta))\},
\end{align*}

For the Scheme A and B, since the variable-coefficient matrix for the asymmetric linear system evolved, the HSS method can not be solved explicitly. The inexact HSS (IHSS) method should be considered.

The IHSS iteration method. Given an initial guess $\boldsymbol{x}^{(0)}$, for $\ell = 0,1,\ldots$, until $\{\boldsymbol{x}^{(\ell)}\}$ converges.

(1) approximate the solution of $(\beta I + H)\boldsymbol{z}^{(\ell+\frac{1}{2})} = \boldsymbol{r}^{(\ell)}$ with $\boldsymbol{r}^{(\ell)} = \boldsymbol{b} - A\boldsymbol{x}^{(\ell)}$
by iterating until $\boldsymbol{z}^{(\ell+\frac{1}{2})}$ is such that the residual
\begin{equation}
\boldsymbol{p}^{(\ell+\frac{1}{2})} = \boldsymbol{r}^{(\ell)} - (\beta I + H)\boldsymbol{z}^{(\ell+\frac{1}{2})}
\end{equation}
satisfies
\[
\|\boldsymbol{p}^{(\ell+\frac{1}{2})}\| \leq \varepsilon_\ell\|\boldsymbol{r}^{(\ell)}\|,
\]
and then compute $\boldsymbol{x}^{(\ell+\frac{1}{2})} = \boldsymbol{x}^{(\ell)} + \boldsymbol{z}^{(\ell+\frac{1}{2})}$;

(2) approximate the solution of $(\beta I + S)\boldsymbol{z}^{(\ell+1)} = \boldsymbol{r}^{(\ell+\frac{1}{2})} = \boldsymbol{b} - A\boldsymbol{x}^{(\ell+\frac{1}{2})}$ by
iterating until $\boldsymbol{z}^{(\ell+1)}$ is such that the residual
\begin{equation}
\boldsymbol{q}^{(\ell+1)} = \boldsymbol{r}^{(\ell+\frac{1}{2})} - (\beta I + S)\boldsymbol{z}^{(\ell+1)}
\end{equation}
satisfies
\[
\|\boldsymbol{q}^{(\ell+1)}\| \leq \eta_\ell\|\boldsymbol{r}^{(\ell+\frac{1}{2})}\|
\]
and then compute $\boldsymbol{x}^{(\ell+1)} = \boldsymbol{x}^{(\ell+\frac{1}{2})} + \boldsymbol{z}^{(\ell+1)}$. Here $\|\cdot\|$ is a norm of a vector.
\section{Numerical experiments}
\label{sec:experiments}

The sparsity and non-Hermitian property can be easily verified numerically in 1D and 3D. For the positiveness, the results are presented in

\begin{table}[htbp]
    \centering
    \begin{tabular}{c|c|c|c}
    \hline
      $N_x$   &$N_t$& $\min \textbf{eig}\left(\frac{A+A^{*}}{2}\right)$ & $\min \textbf{eig}\left(\frac{B+B^{*}}{2}\right)$  \\
         \hline
         64&100 &1.831359910960193 &1.825702675166990 \\
         128&100&1.831359132456233&1.824550880717313\\
         500&100&1.831359080012343&1.823691838263642\\
         500&10&1.813590359902322&1.736916484841441\\
         \hline
    \end{tabular}
    \begin{tabular}{c|c|c|c}
    \hline
      $N_x$   &$N_t$& $\min \textbf{eig}\left(\frac{A+A^{*}}{2}\right)$ & $\min \textbf{eig}\left(\frac{B+B^{*}}{2}\right)$  \\
         \hline
         64&100 &1.833308730408240 &1.832885749784061 \\
         128&100&1.833308730413283&1.832885747243268 \\
         500&10&1.833087299230090&1.828857415801380\\
         \hline
    \end{tabular}
    \caption{The positiveness for Matrix A and B over time with the damping parameter $\alpha=0.01$ (top) and $\alpha=0.5$ (bottom) in 1D up to the final time $T=0.1$. The smallest eigenvalue is positive. Since that the dominant part of the matrices is $11/6=1.833$ when $\alpha$, $\Delta t/h^2$ are smaller.}
    \label{tab:placeh}
\end{table}

\begin{table}[htbp]
    \centering
    \begin{tabular}{c|c|c|c|c|c}
    \hline
    $N_x$&$N_x$  &$N_x$   &$N_t$& $\min \textbf{eig}\left(\frac{A+A^{*}}{2}\right)$ & $\min \textbf{eig}\left(\frac{B+B^{*}}{2}\right)$  \\
         \hline
         4&4&4&40 &1.831606958090378&1.830877241360574 \\
         6&6&6&40 &1.831095310243323&1.829711010777179 \\
         8&8&8&40&1.830753045464913&1.828639160121246\\
         \hline
    \end{tabular}
    \begin{tabular}{c|c|c|c|c|c}
    \hline
    $N_x$&$N_x$  &$N_x$   &$N_t$& $\min \textbf{eig}\left(\frac{A+A^{*}}{2}\right)$ & $\min \textbf{eig}\left(\frac{B+B^{*}}{2}\right)$  \\
         \hline
         4&4&4&40 &1.833286773195681&1.831933551302752 \\
         6&6&6&40 &1.833287121445087&1.831920706743620 \\
         8&8&8&40&1.833287185184424&1.831916528114020\\
         \hline
    \end{tabular}
    \caption{The positiveness for Matrix A and B over time with the damping parameter $\alpha=0.01$ (top) and $\alpha=0.5$ (bottom) in 3D up to the final time $T=0.1$. The smallest eigenvalue is positive. Since that the dominant part of the matrices is $11/6=1.833$ when $\alpha$, $\Delta t/h^2$ are smaller.}
    \label{tab:placeh}
\end{table}

\begin{figure}[htbp]
	\centering
	\subfloat[1D ]{\label{cputime_ED_1D_time}\includegraphics[width=1.5in]{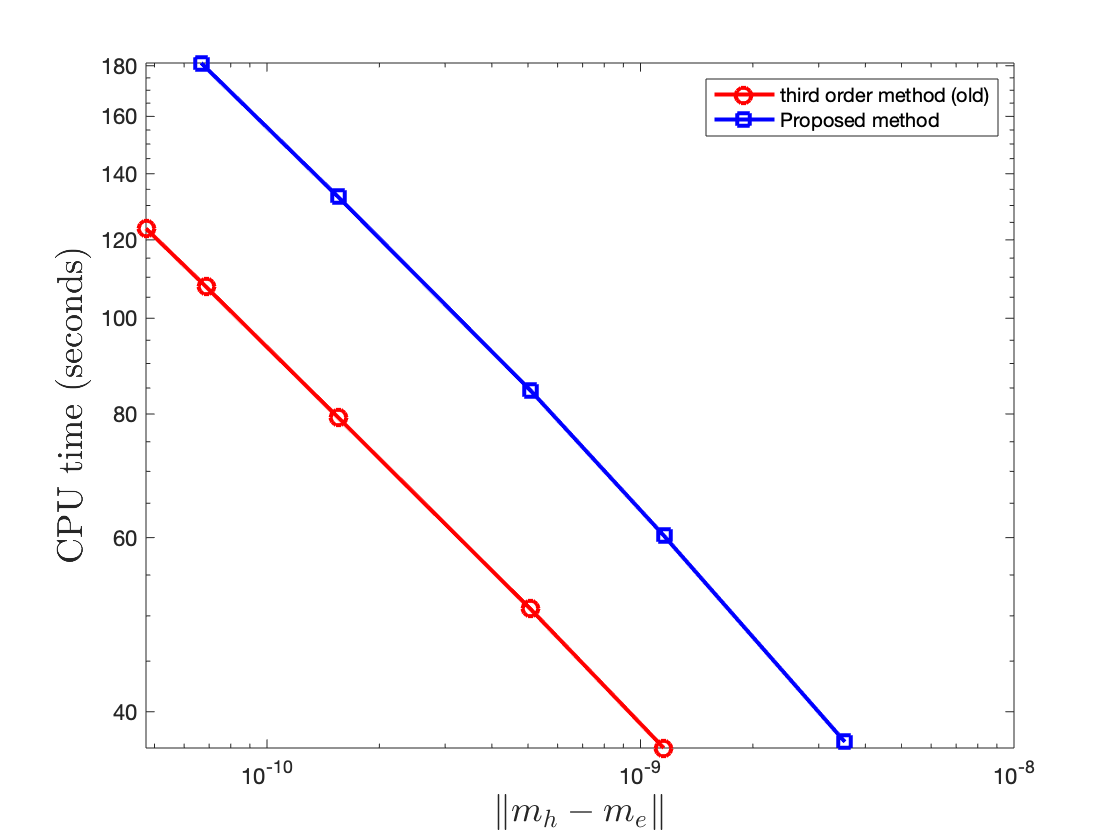}}
		\subfloat[3D ]{\label{cputime_ED_3D_time}\includegraphics[width=1.5in]{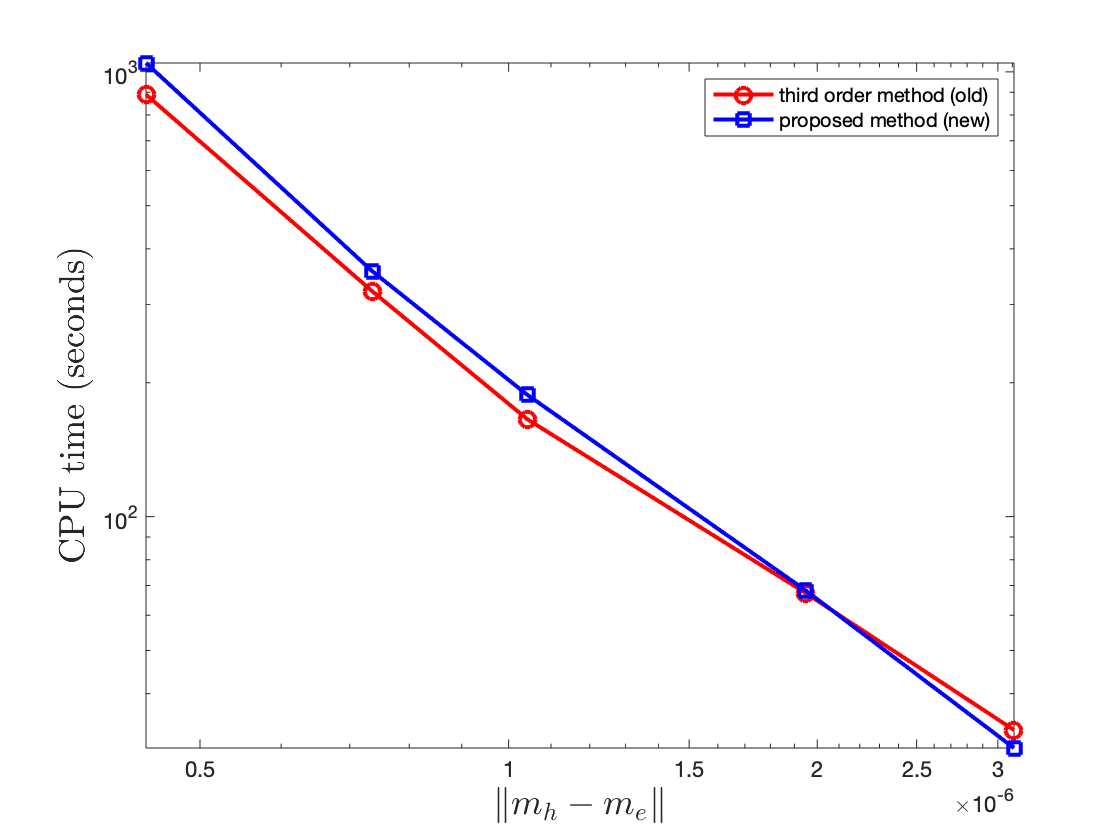}}
	\subfloat[1D ]{\label{cputime_ED_1D_space}\includegraphics[width=1.5in]{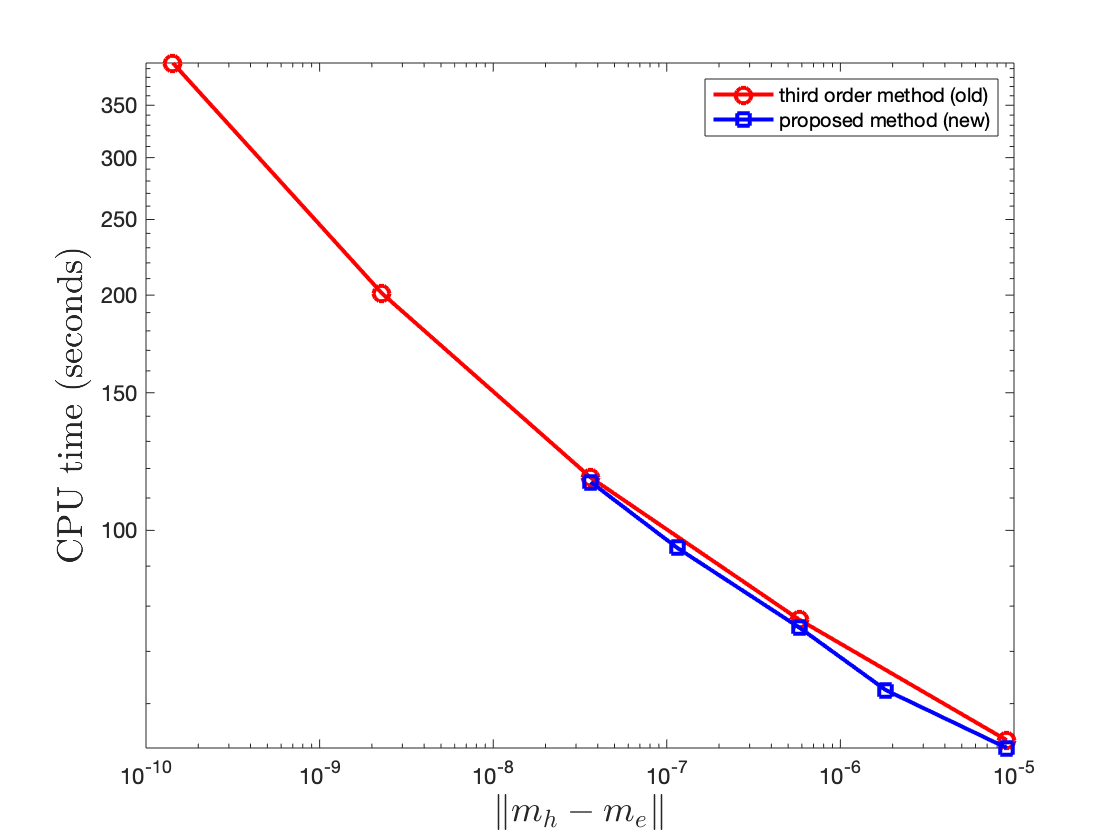}}
		\subfloat[3D ]{\label{cputime_ED_3D_space}\includegraphics[width=1.5in]{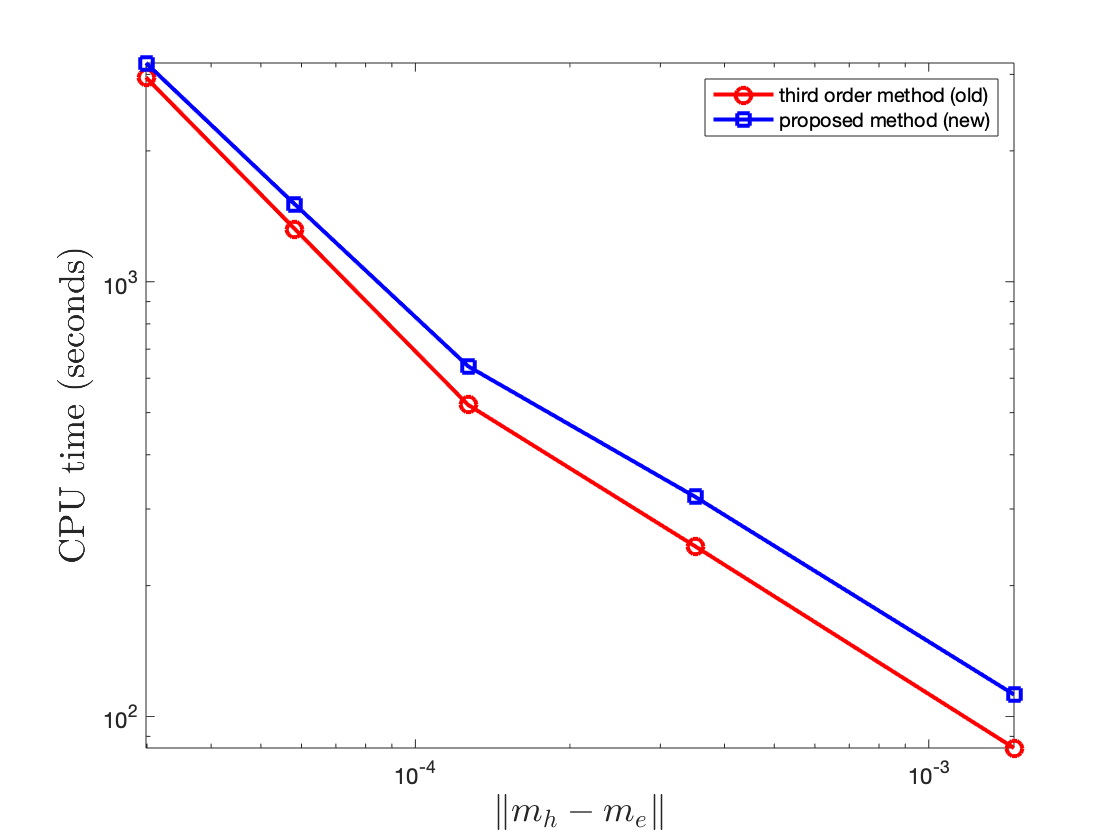}}
		\hspace{0.1in}
		\subfloat[1D ]{\label{cputime_ED_BDF2_1D_time}\includegraphics[width=1.5in]{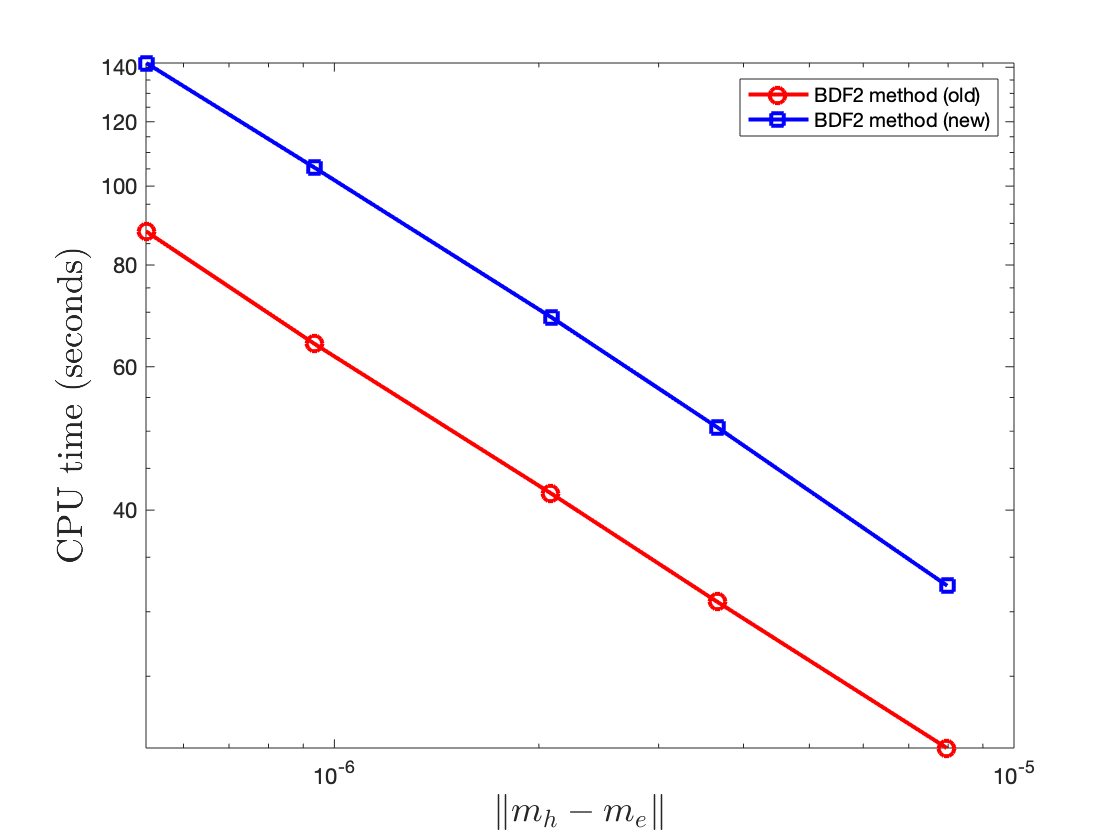}}
		\subfloat[3D ]{\label{cputime_ED_BDF2_3D_time}\includegraphics[width=1.5in]{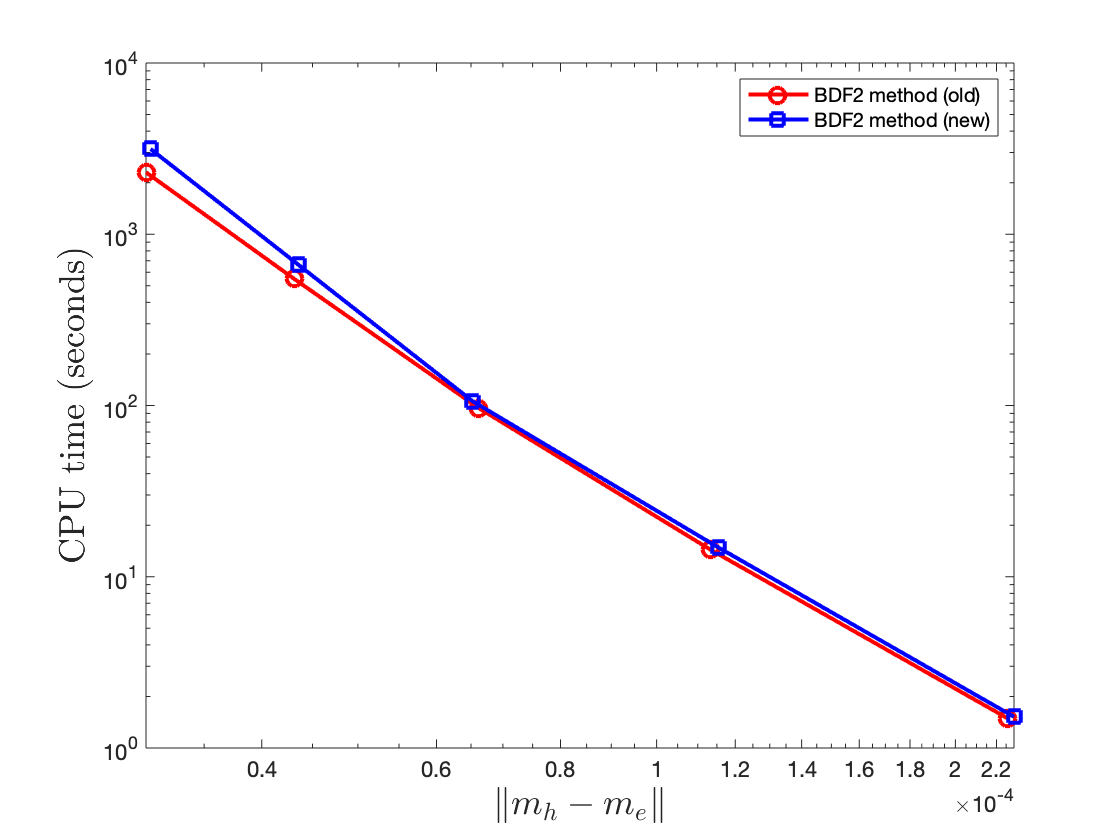}}
		\subfloat[1D ]{\label{cputime_ED_BDF2_1D_space}\includegraphics[width=1.5in]{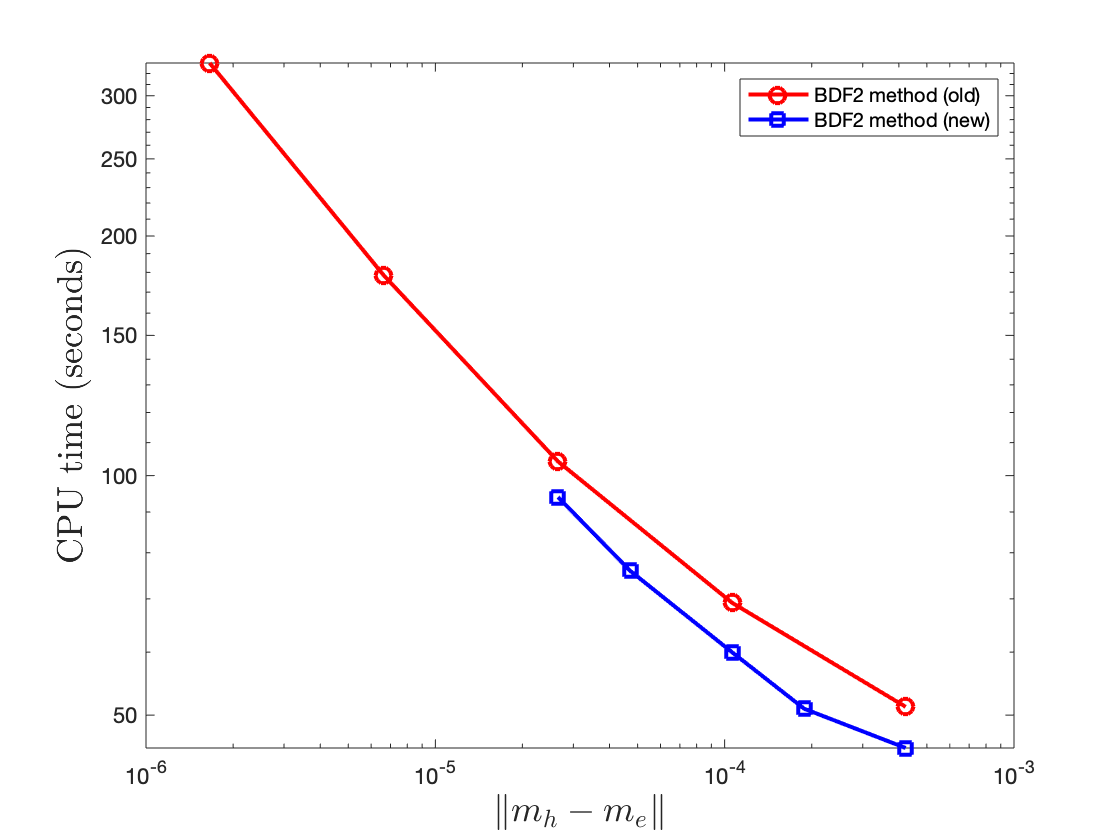}}
		\subfloat[3D ]{\label{cputime_ED_BDF2_3D_space}\includegraphics[width=1.5in]{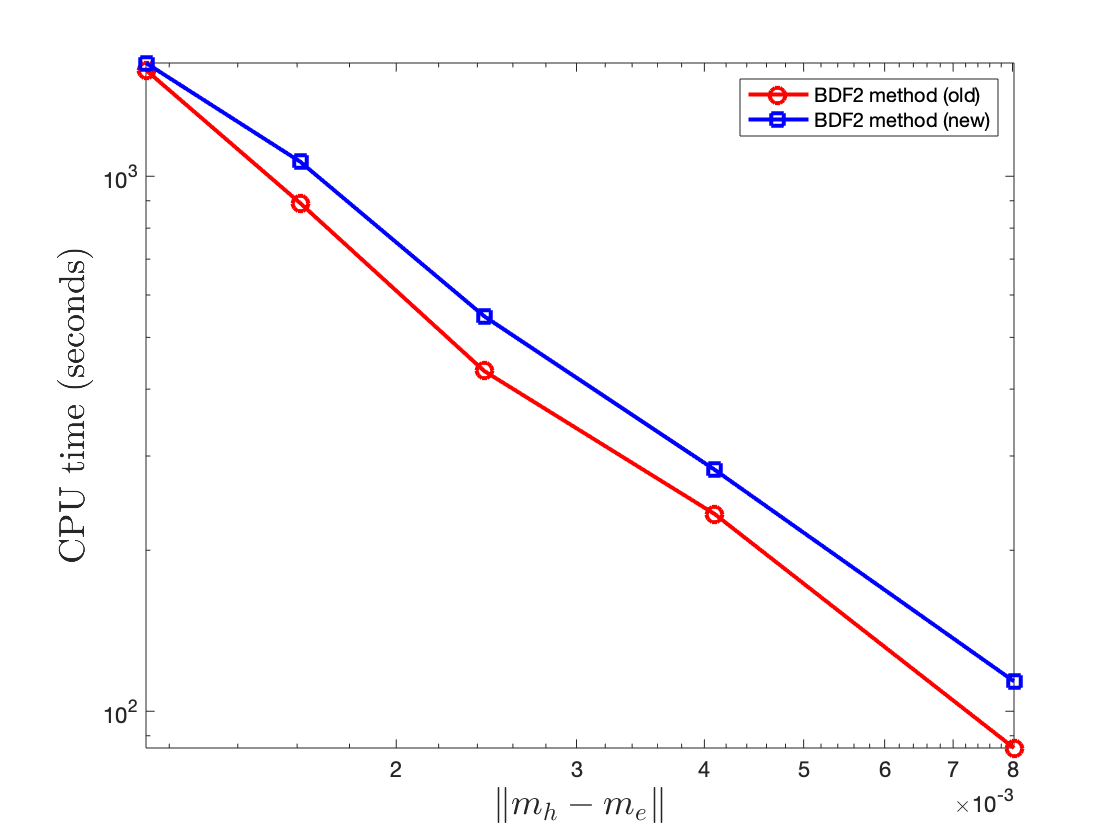}}
		\hspace{0.1in}
		\subfloat[1D ]{\label{cputime_ED_BDF1_1D_time}\includegraphics[width=1.5in]{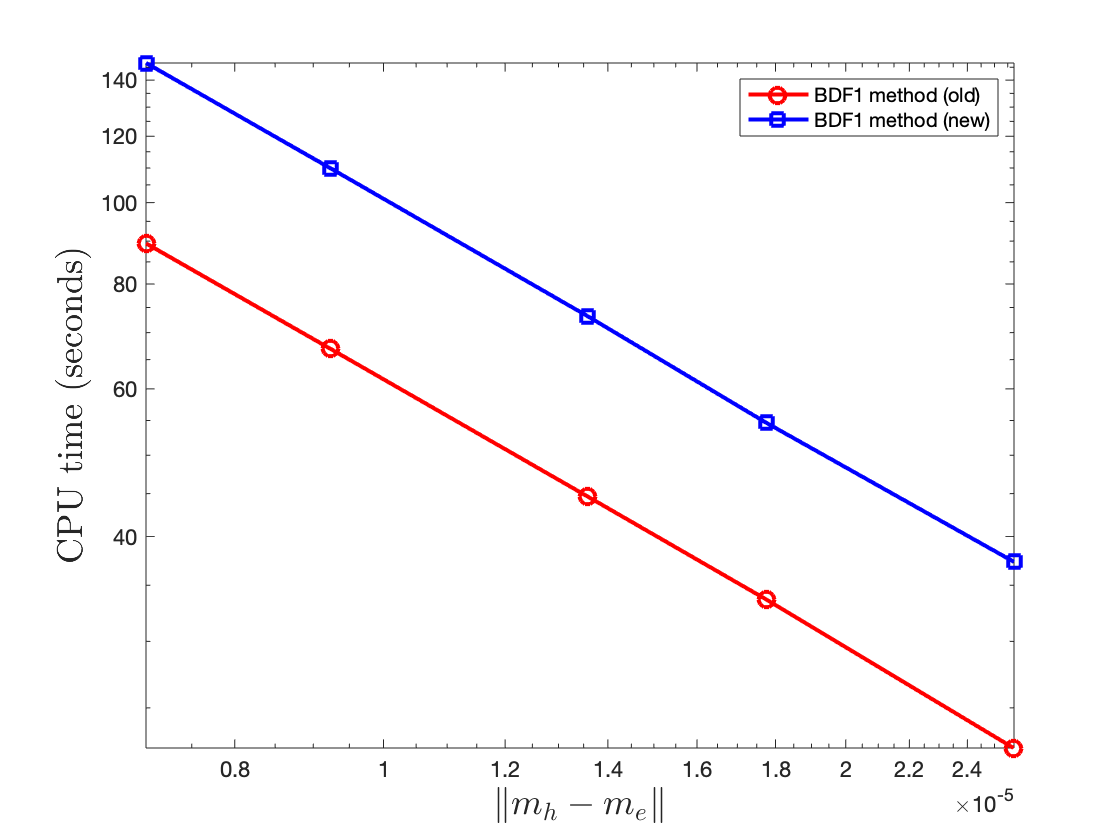}}
		\subfloat[3D ]{\label{cputime_ED_BDF1_3D_time}\includegraphics[width=1.5in]{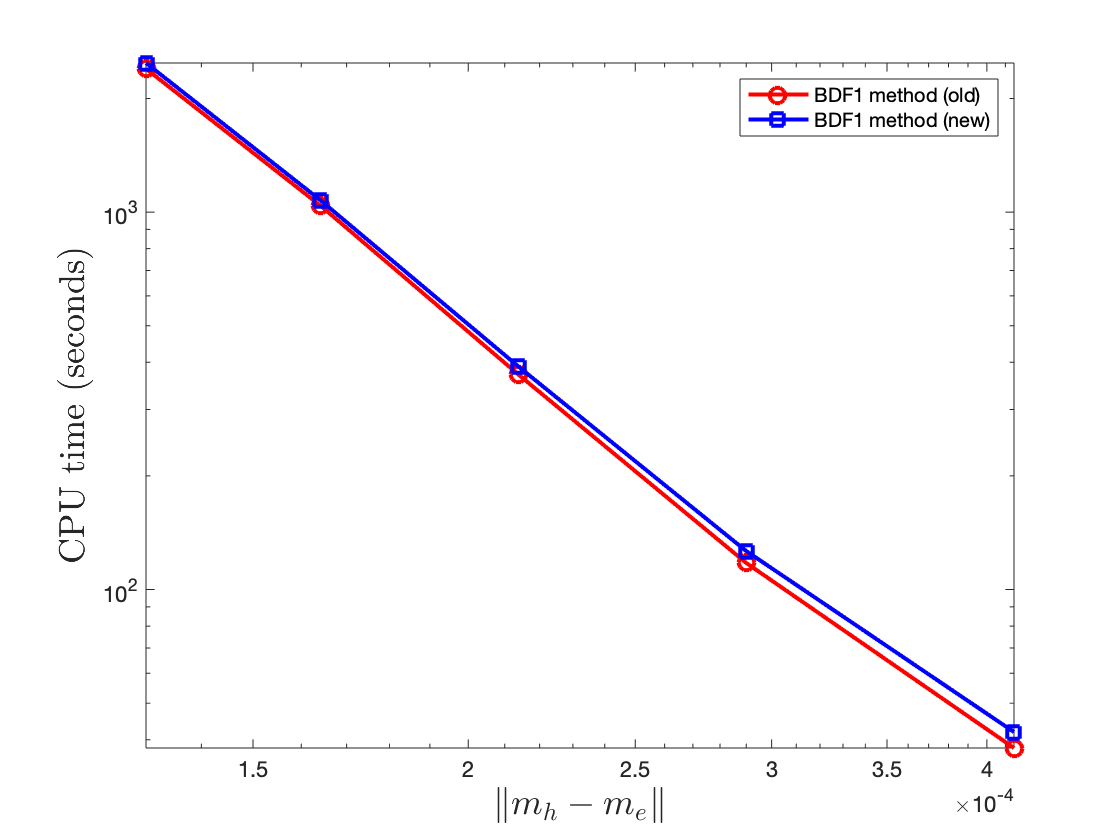}}
		\subfloat[1D ]{\label{cputime_ED_BDF1_1D_space}\includegraphics[width=1.5in]{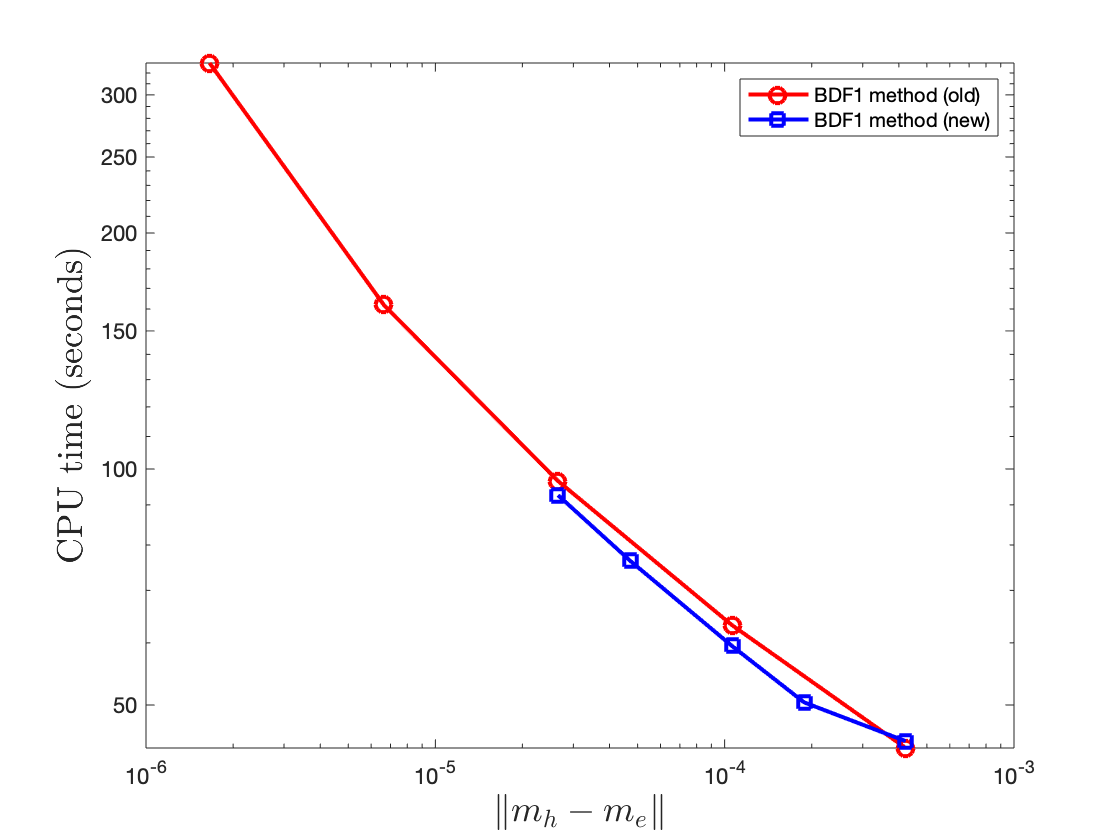}}
		\subfloat[3D ]{\label{cputime_ED_BDF1_3D_space}\includegraphics[width=1.5in]{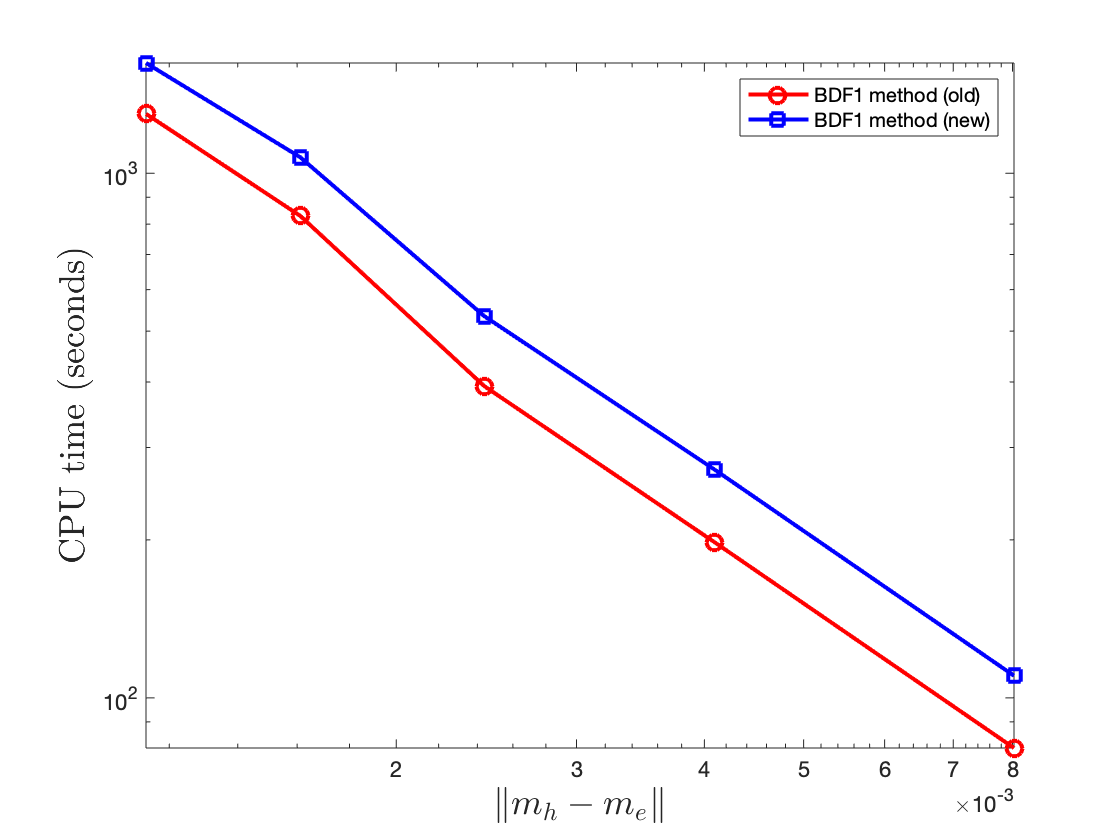}}
	\caption{Comparison of CPU time required to achieve the desired numerical accuracy for the proposed method for \cref{eq-5} and the third order semi-implicit scheme for \cref{eq-model}, along with their BDF2 and BDF1 counterparts. Top row: BDF3 method; Middle row: BDF2 method; Bottom row: BDF1 method. The observation is that the proposed method consumes more time than previous methods to achieve the same level of accuracy, indicating that its linear system is more difficult to solve. Numerical evidence reveals that GMRES solver convergence slows for larger temporal step size $k$ or smaller spatial grid-size $h$, increasing computational challenge. Thus, the proposed method requires a more efficient solver. Here the results in original code using \text{diag}, while \texttt{spdiags} will be faster. The conclusion is the same. From left to right: first two with varying the time-step size $k$ and last two with varying the grid-size $h$.}\label{figure:3order}
\end{figure}

Additionally, to compare numerical efficiency between the proposed method \cref{proposed} and the revised third-order semi-implicit scheme \cref{scheme-third-order}, we plot CPU time (seconds) versus the $\ell^{\infty}$ norm. Results are presented in \cref{figure:3order}. We observe that the scheme B requires more time than previous methods to achieve the same accuracy level. This indicates that the linear system for the scheme B is more challenging to solve. 

\subsection{Spectral radius}\label{sec:radii}

\begin{figure}[htbp]
    \centering
    \subfloat{\includegraphics[width=0.25\linewidth]{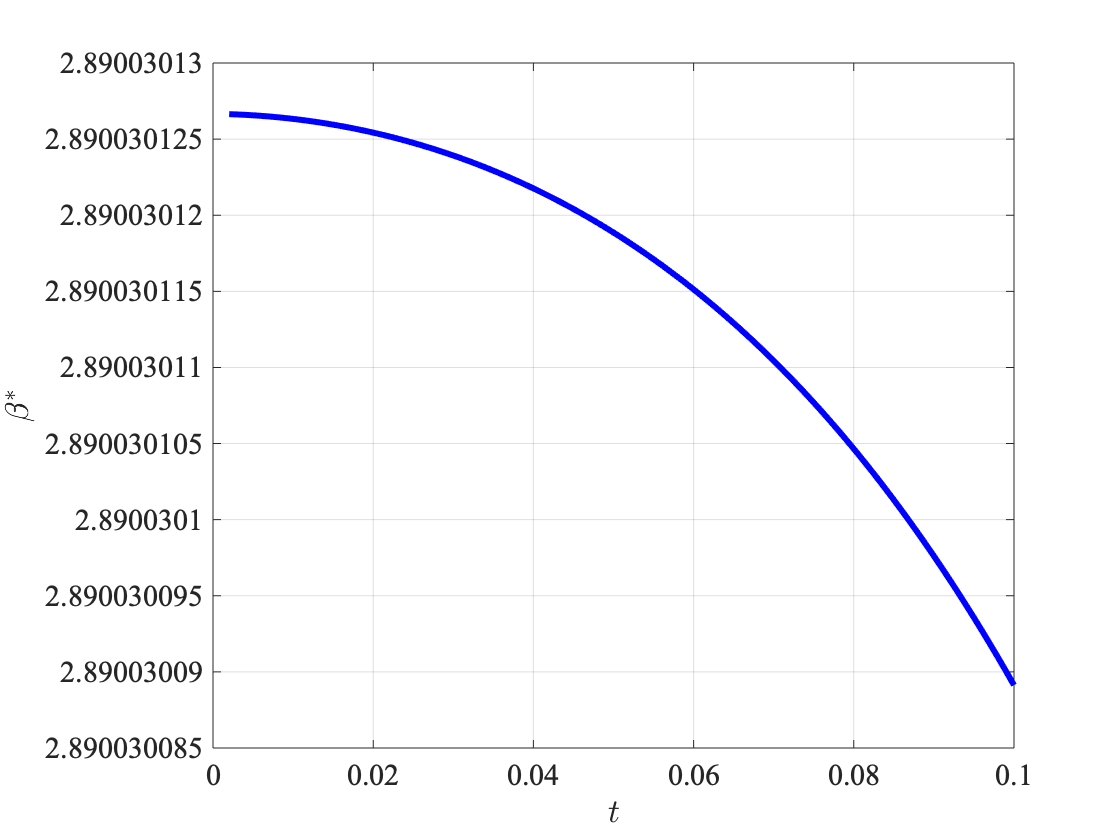}}
    \subfloat{\includegraphics[width=0.25\linewidth]{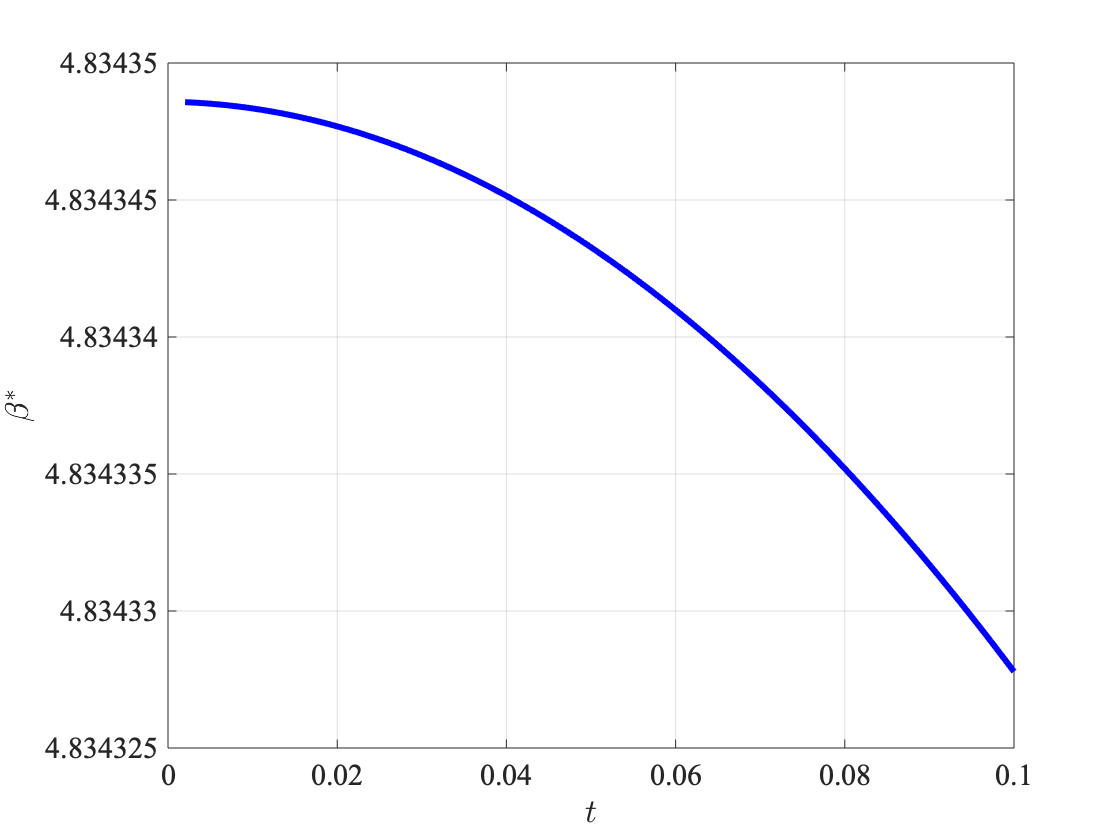}}
     \subfloat{\includegraphics[width=0.25\linewidth]{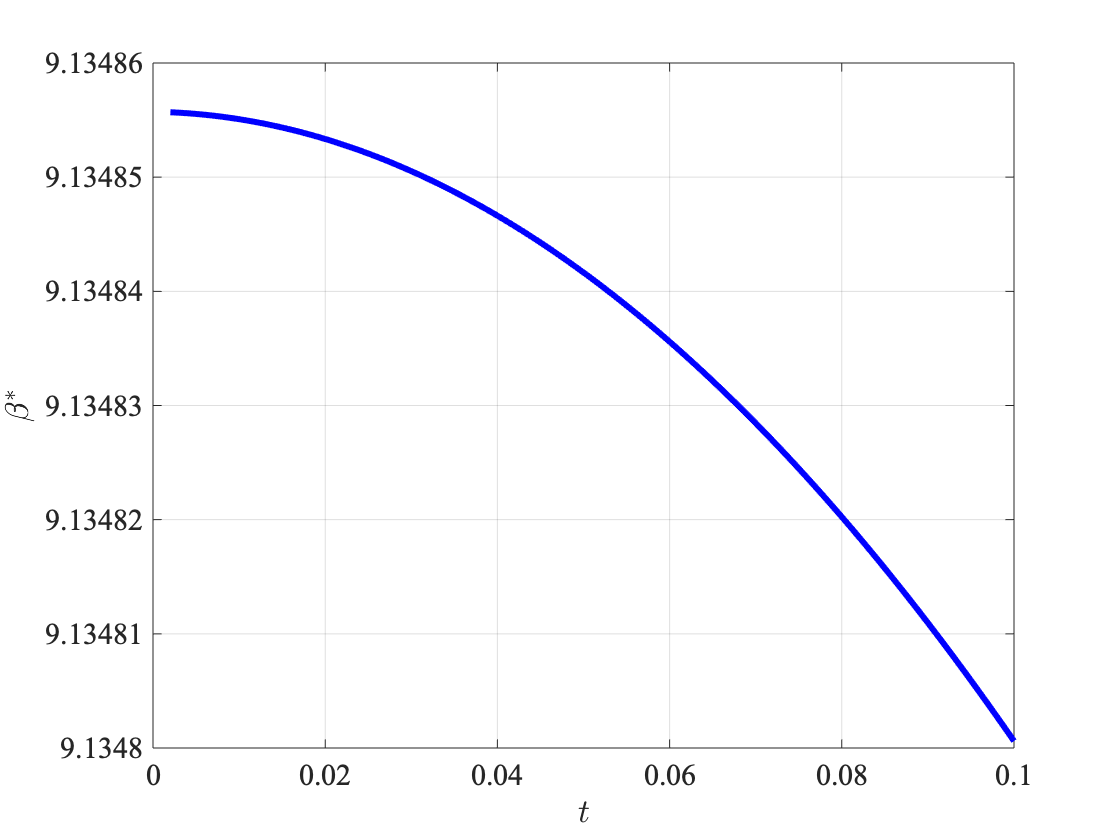}}
     \hspace{0.1in}
     \subfloat{\includegraphics[width=0.25\linewidth]{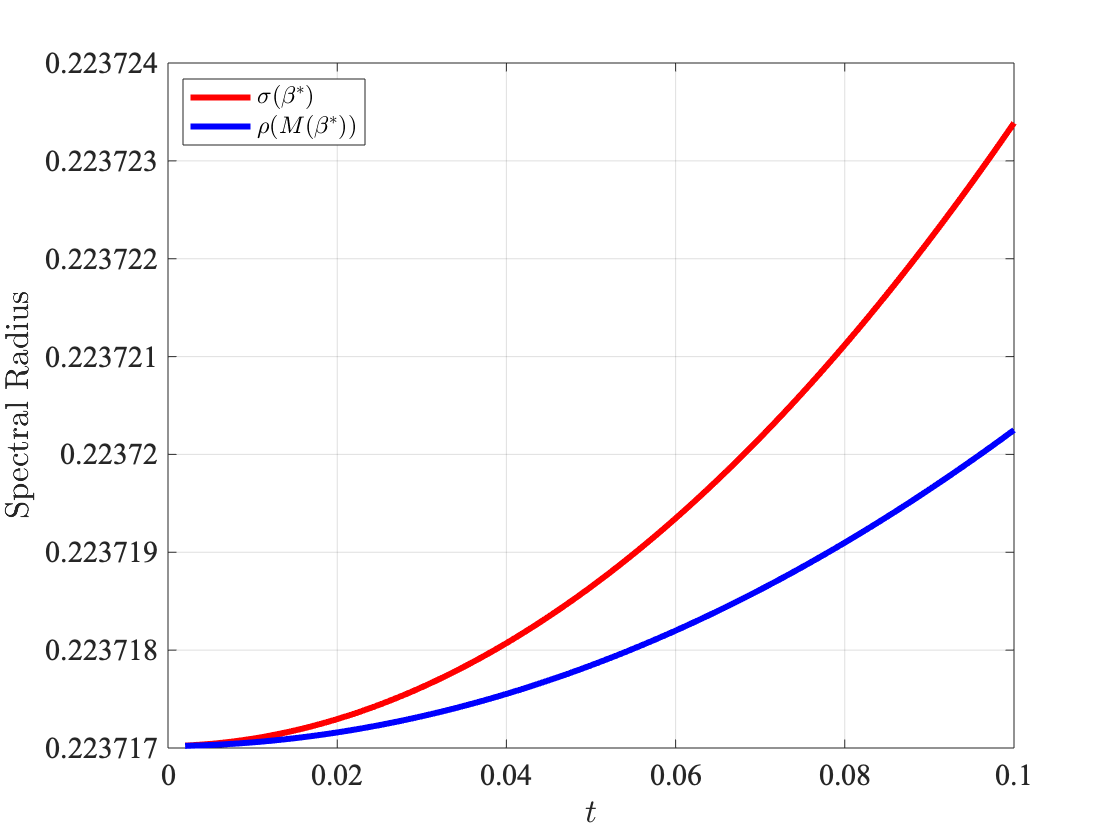}}
    \subfloat{\includegraphics[width=0.25\linewidth]{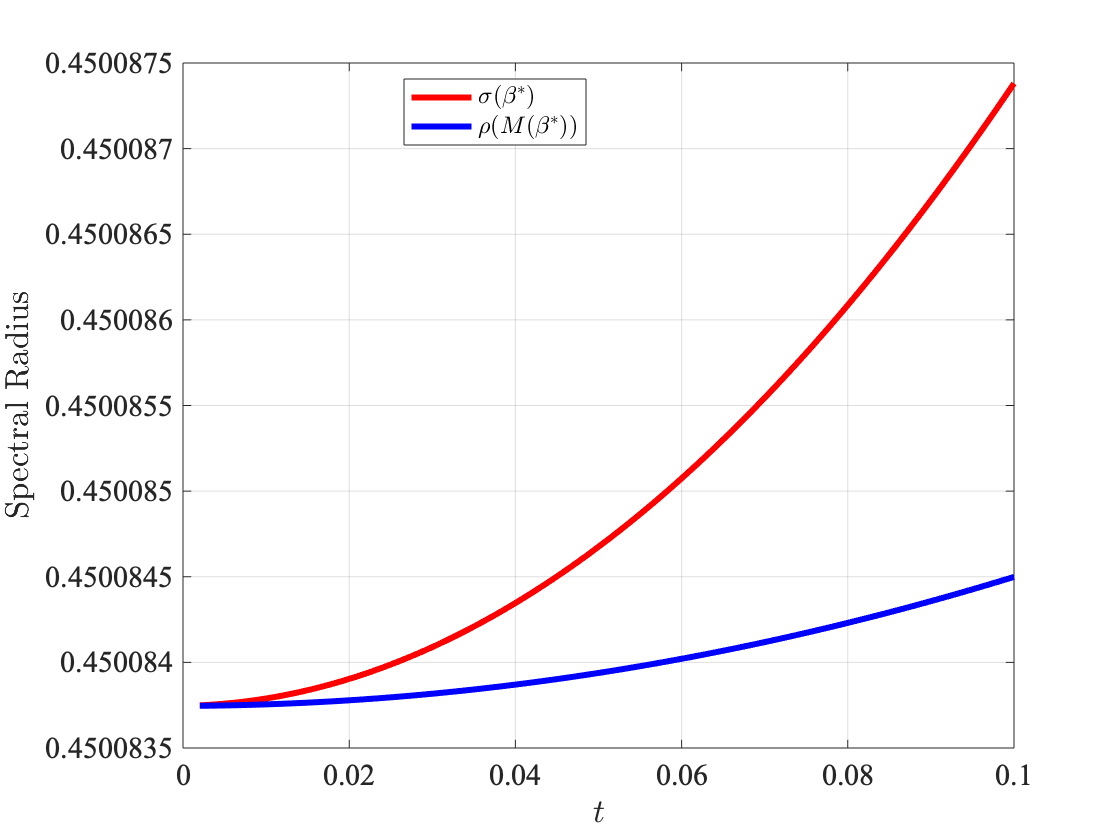}}
     \subfloat{\includegraphics[width=0.25\linewidth]{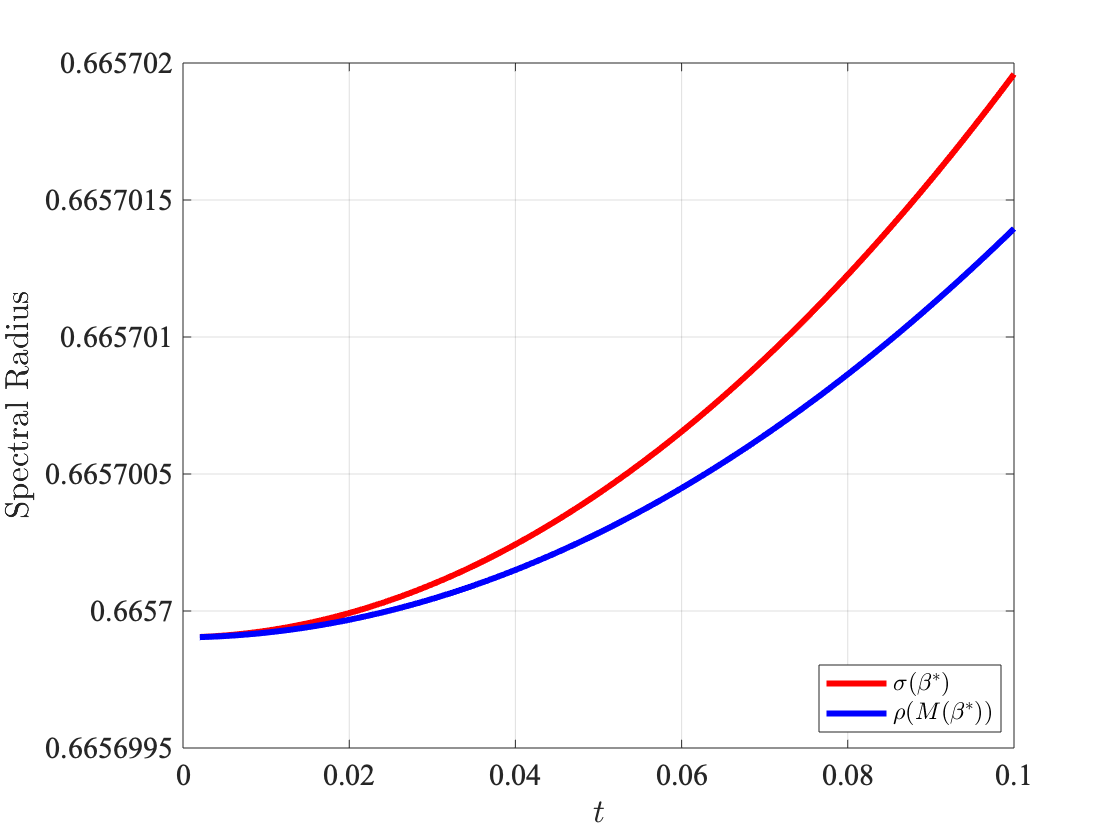}}
     \hspace{0.1in}
     \subfloat{\includegraphics[width=0.25\linewidth]{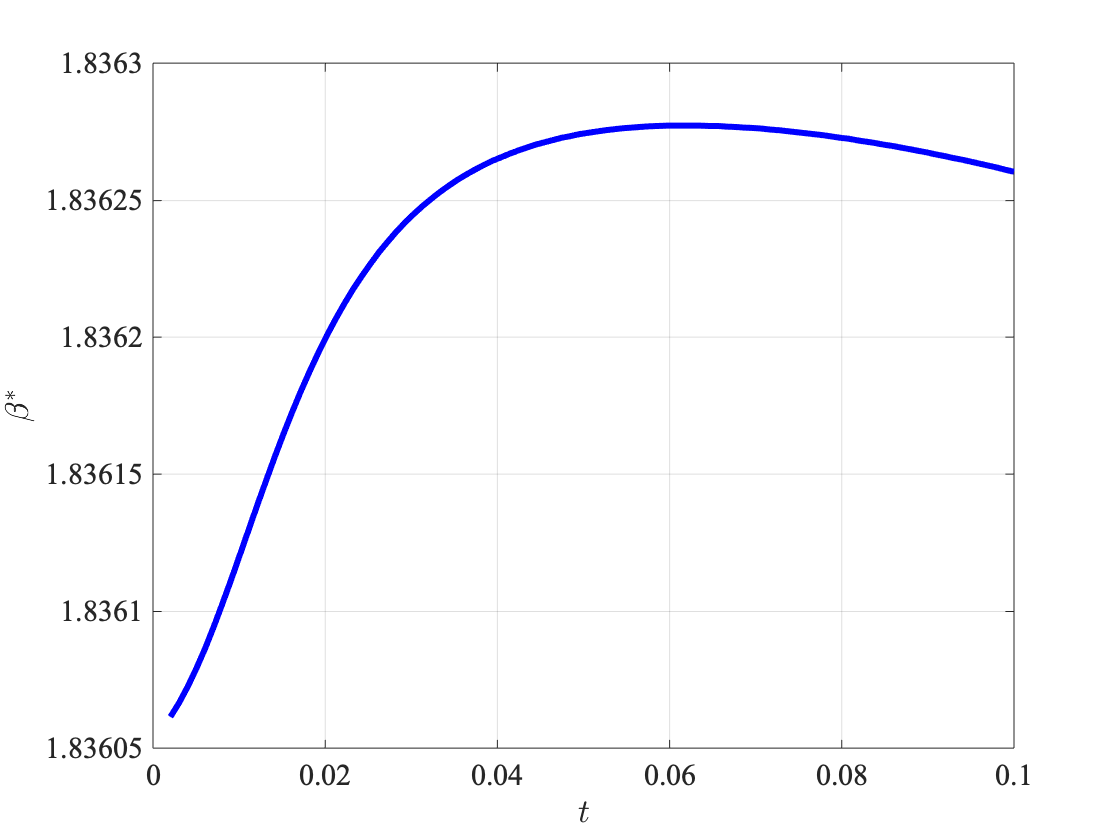}}
    \subfloat{\includegraphics[width=0.25\linewidth]{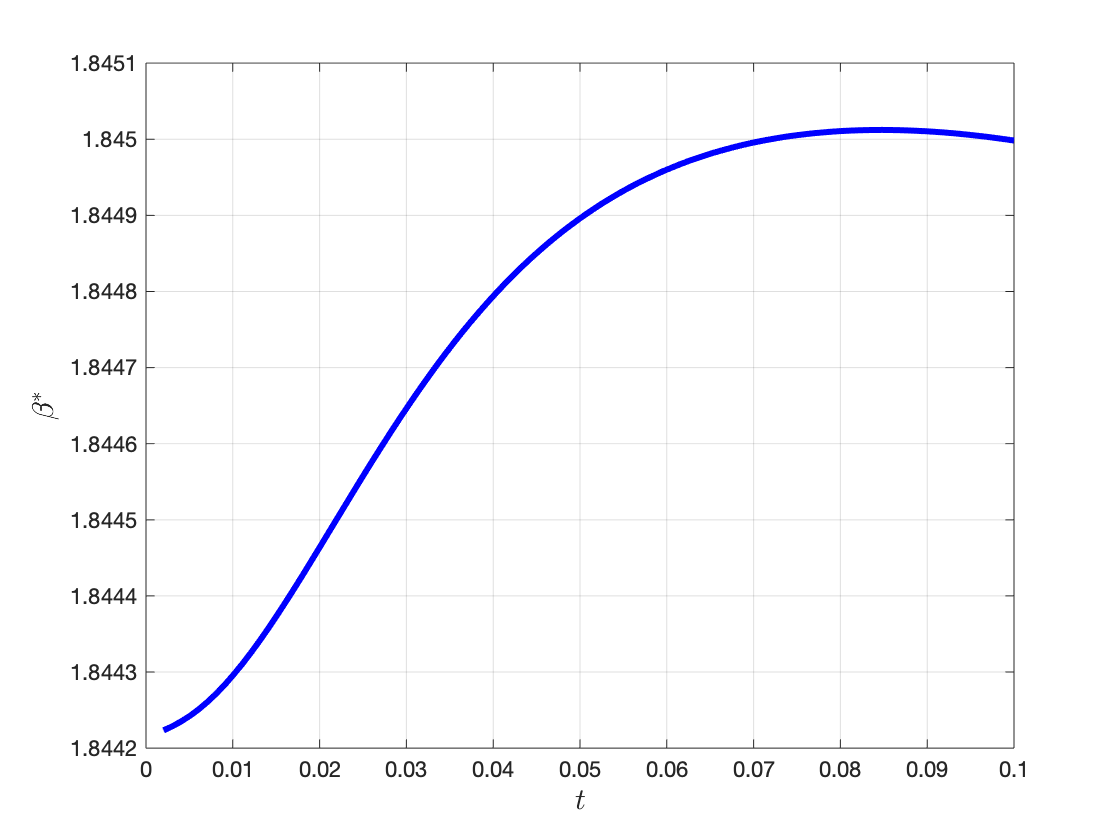}}
     \subfloat{\includegraphics[width=0.25\linewidth]{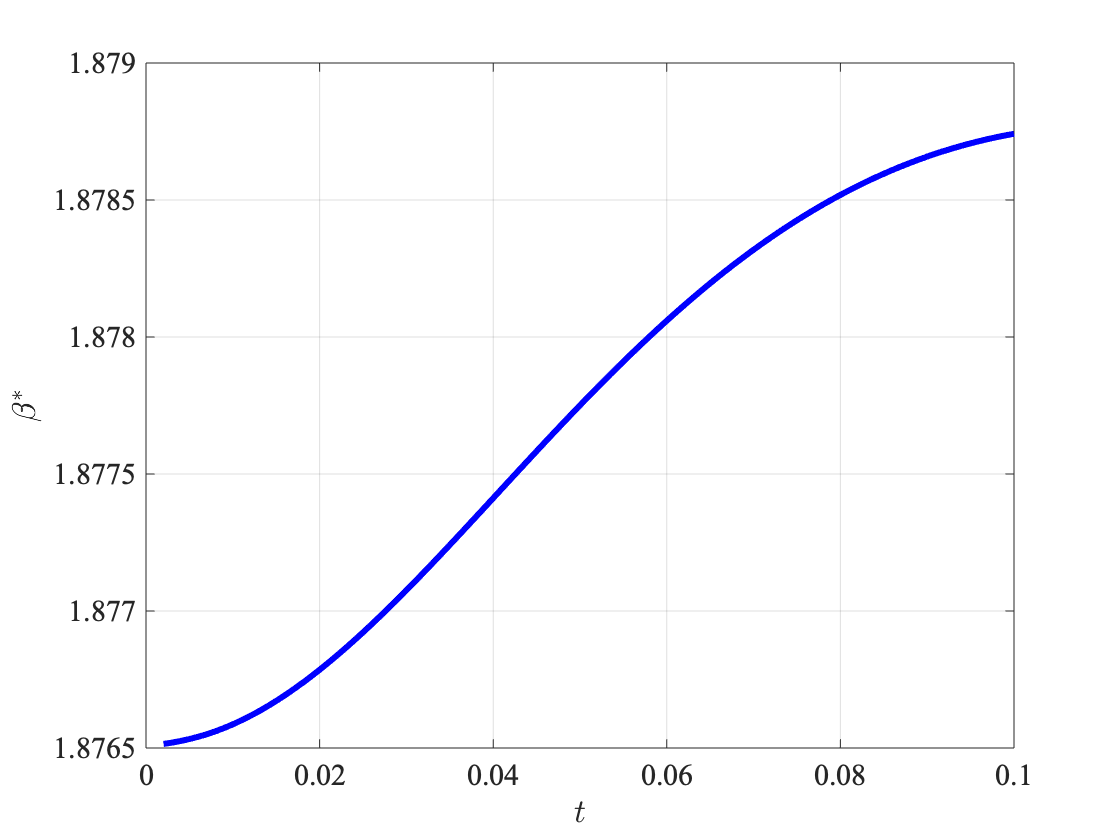}}
     \hspace{0.1in}
     \subfloat{\includegraphics[width=0.25\linewidth]{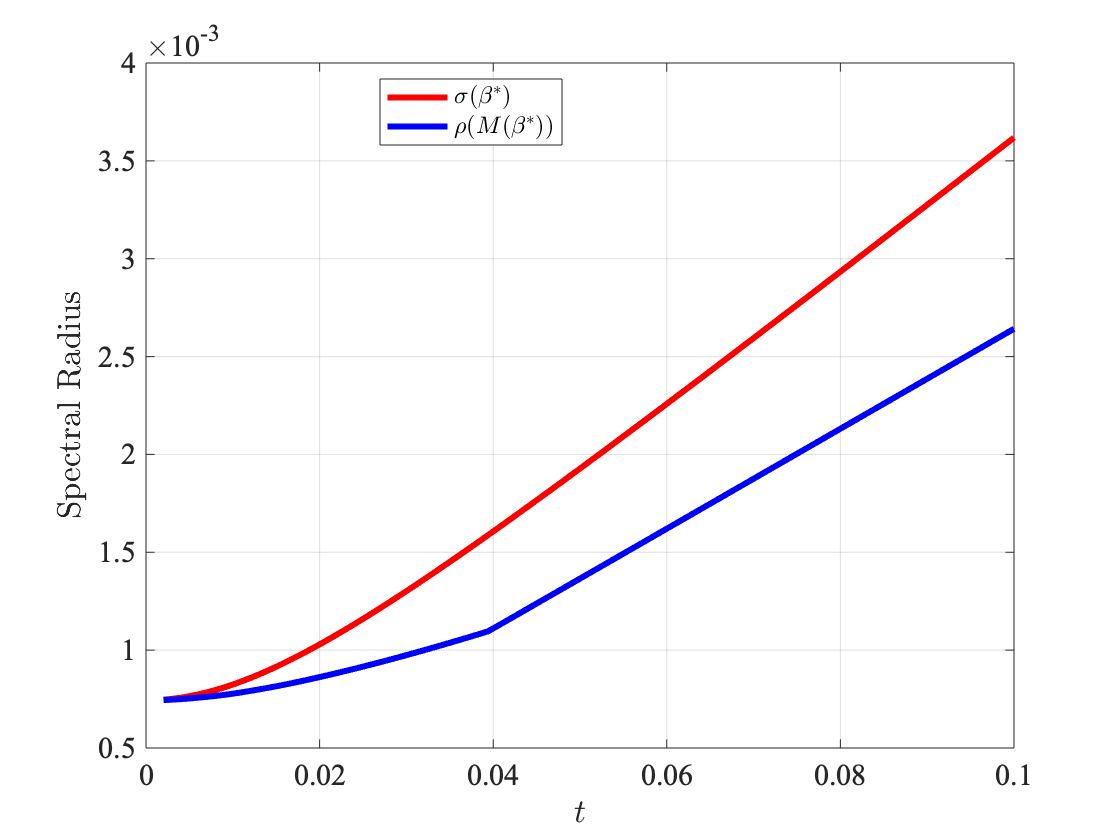}}
    \subfloat{\includegraphics[width=0.25\linewidth]{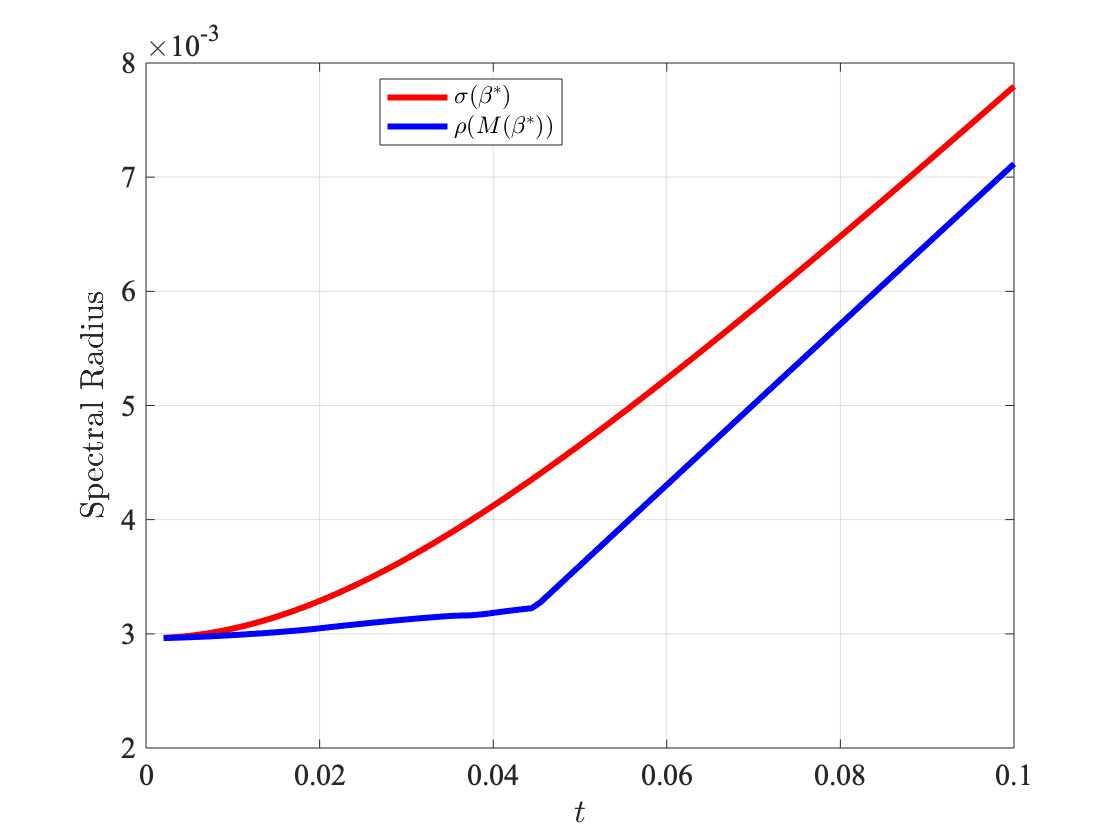}}
     \subfloat{\includegraphics[width=0.25\linewidth]{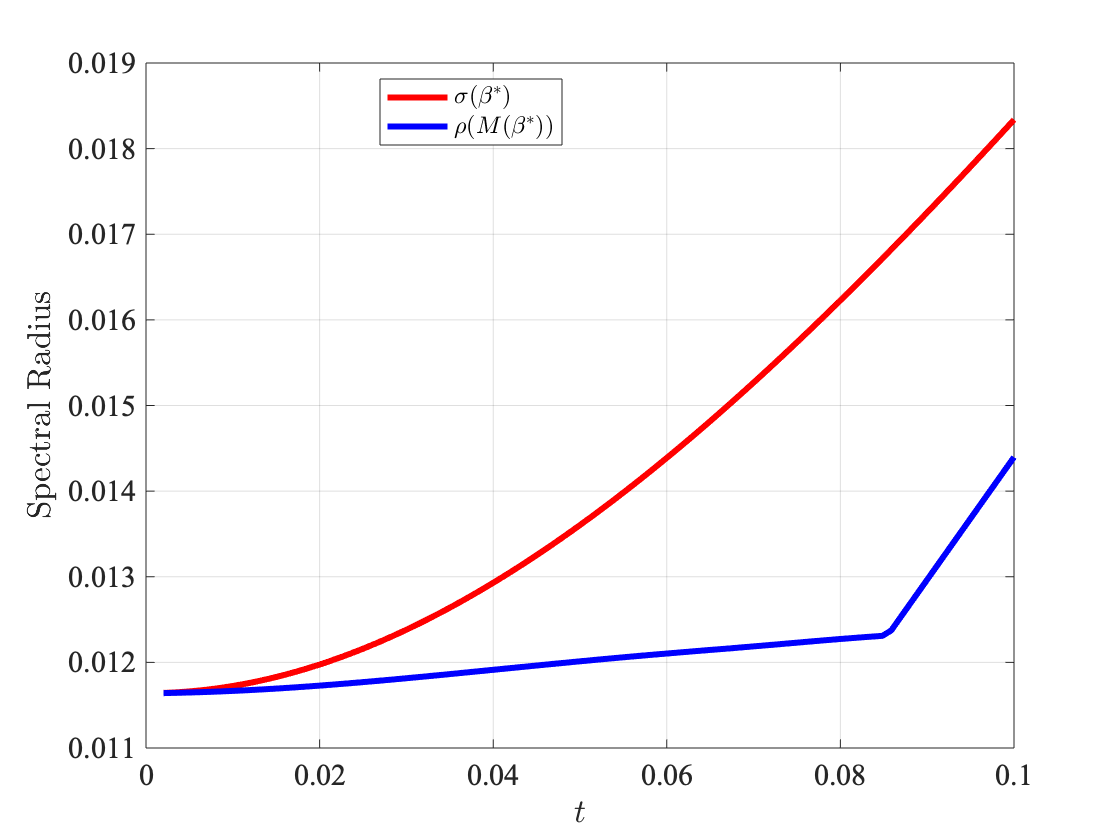}}
     \hspace{0.1in}
     \subfloat{\includegraphics[width=0.25\linewidth]{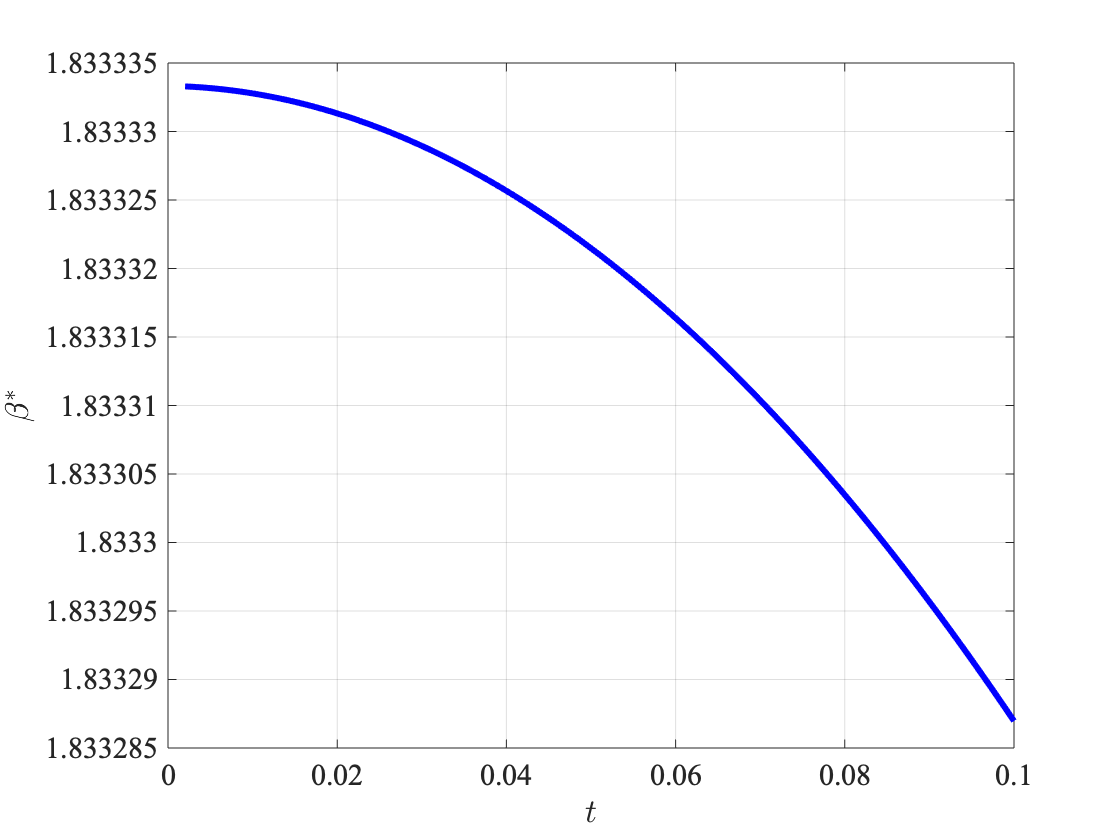}}
    \subfloat{\includegraphics[width=0.25\linewidth]{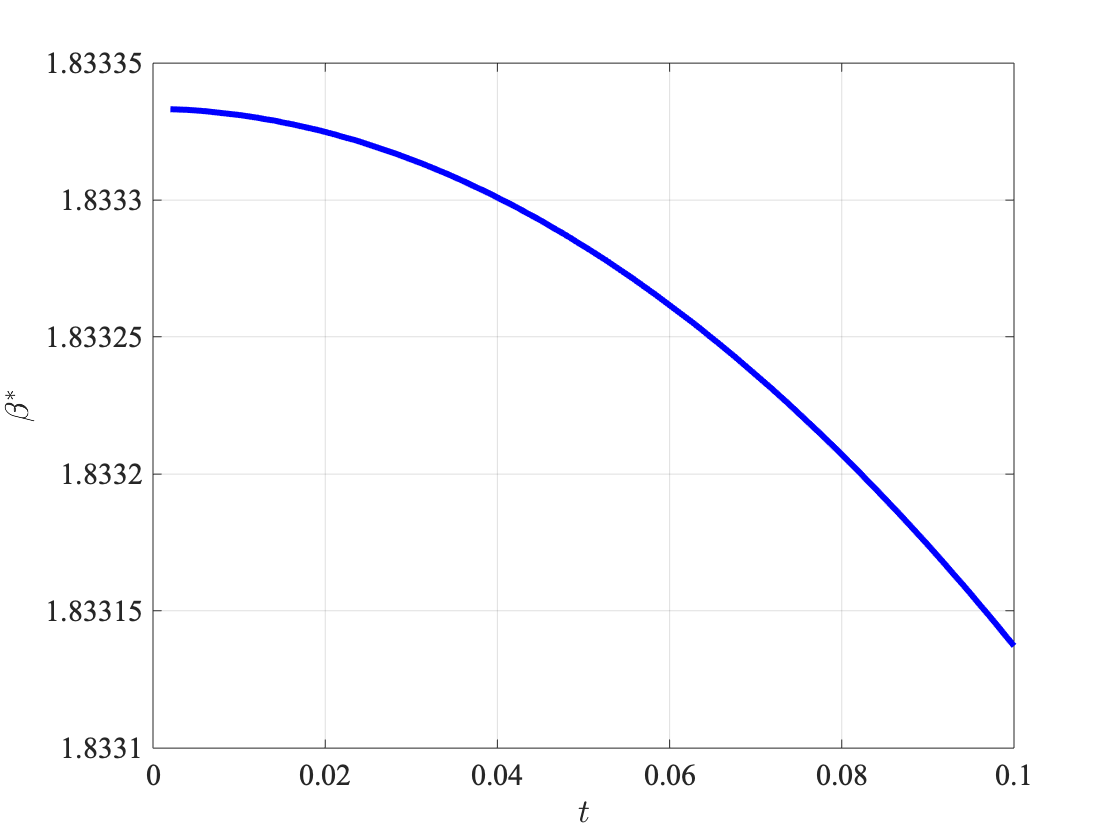}}
     \subfloat{\includegraphics[width=0.25\linewidth]{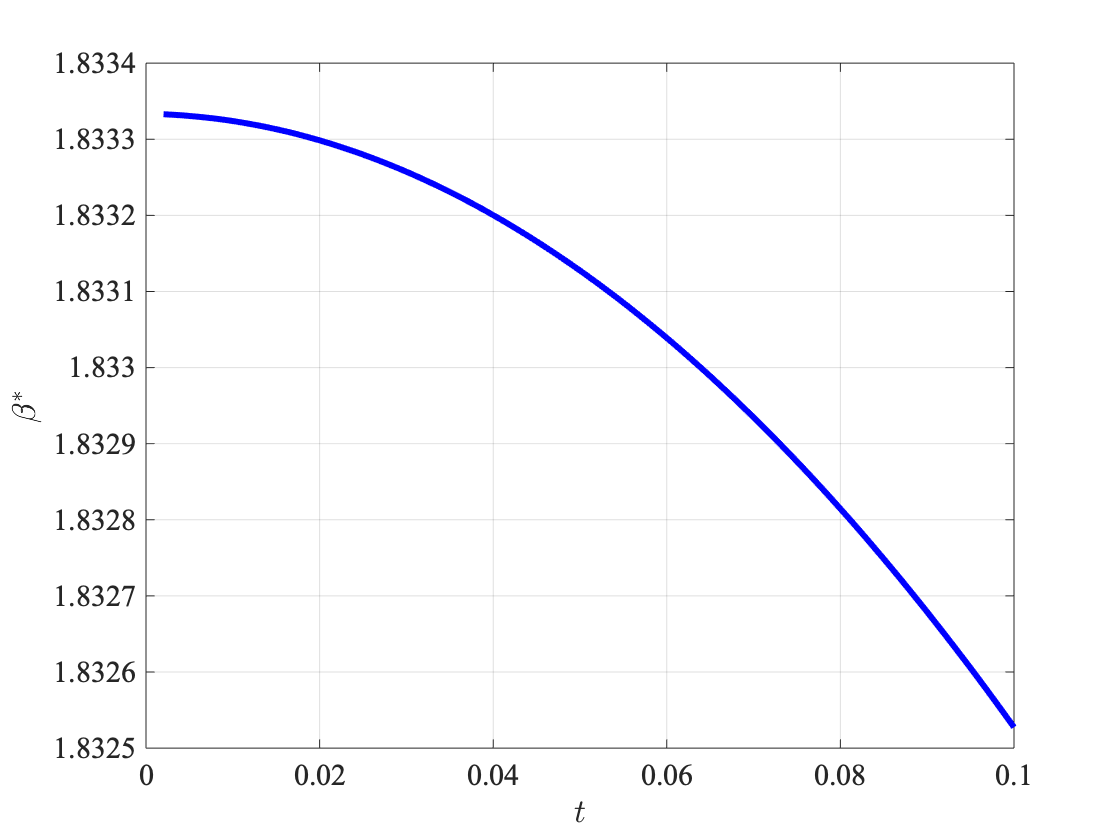}}
     \hspace{0.1in}
     \subfloat{\includegraphics[width=0.25\linewidth]{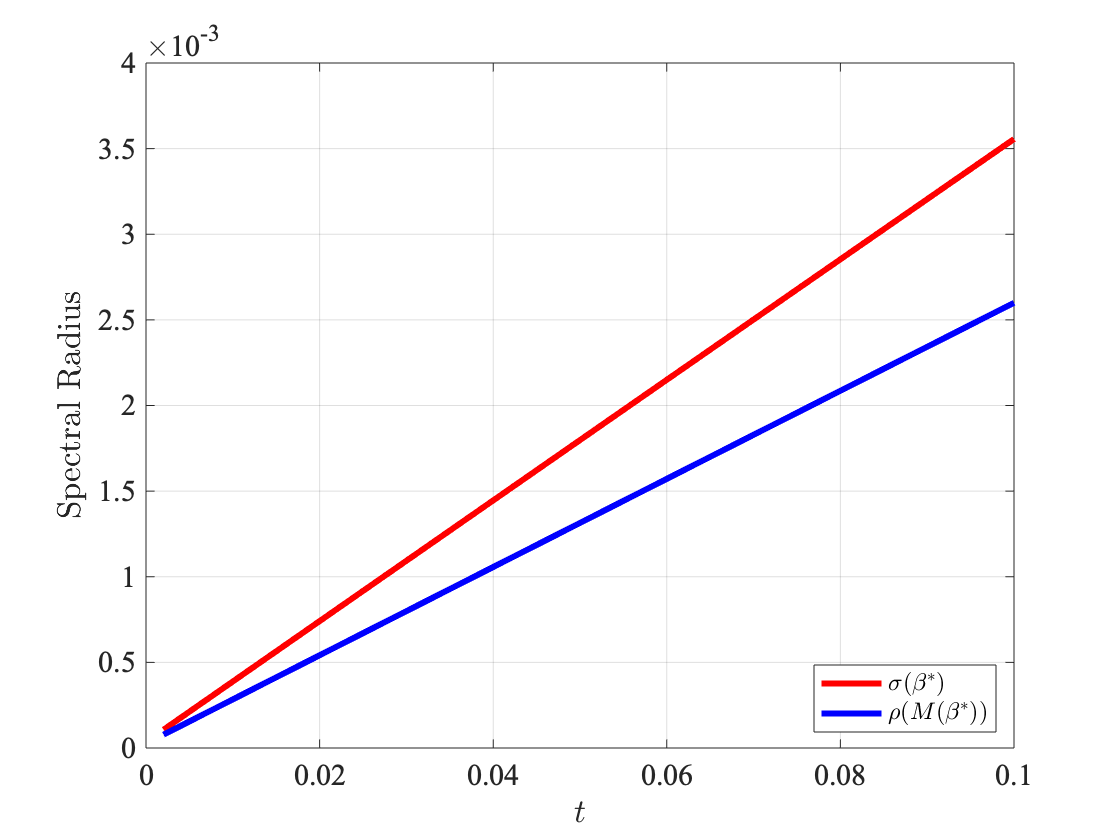}}
    \subfloat{\includegraphics[width=0.25\linewidth]{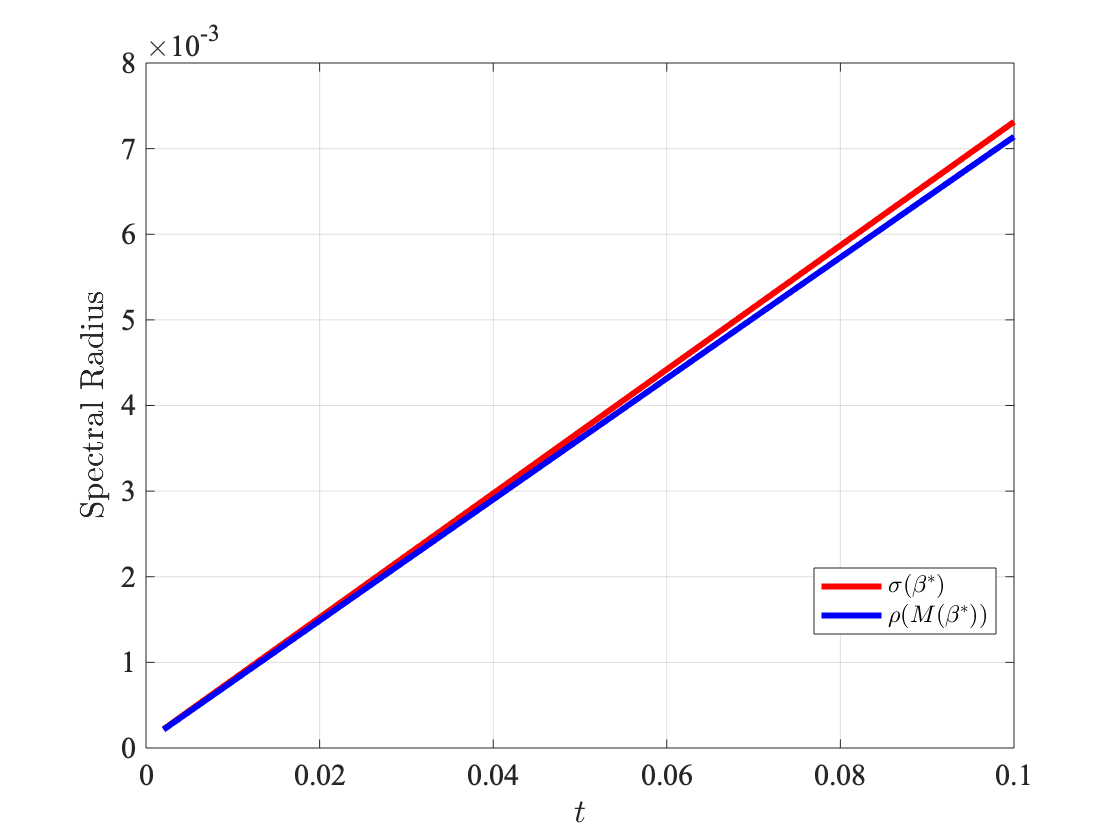}}
     \subfloat{\includegraphics[width=0.25\linewidth]{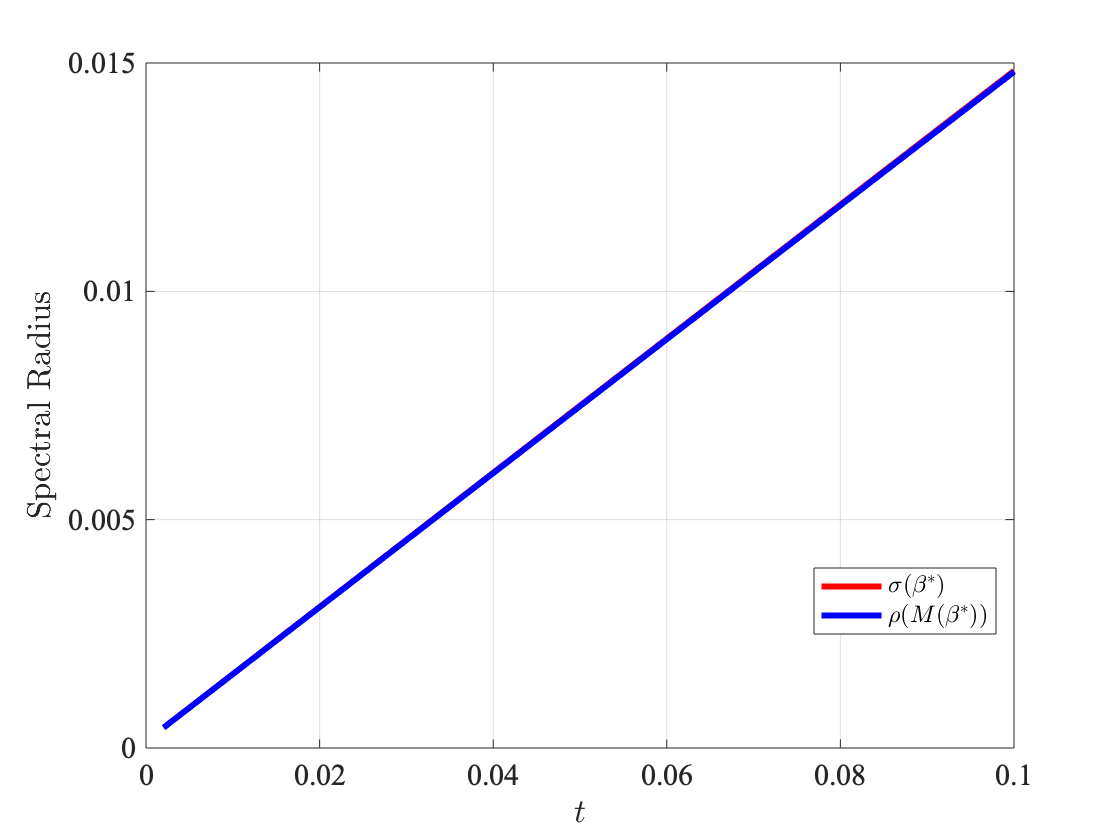}}
    \caption{The estimated $\beta^{*}$, the spectral radius of $\rho(M(\beta^{*}))$ and its upper bound $\sigma(\beta^*)$ over time for scheme A in 1D. From left to right $N_x=32,64,128$ and $N_t=100$. First two rows: $\alpha=0.5$; Middle two rows: $\alpha=0.001$; Bottom two rows: $\alpha=0$.}
    \label{fig:1}
\end{figure}

Figure~\ref{fig:1} shows that the spectral radius $\rho(M(\beta^*))$ remains below the corresponding upper bound $\sigma(\beta^*)$ for all tested spatial resolutions and damping parameters. Although the spectral radius increases gradually in some cases as time evolves, it stays strictly below one throughout the simulation. This numerically verifies the spectral-radius estimate derived in the previous analysis and confirms the convergence of the corresponding iteration. Similar behavior is observed for different values of $N_x$ and $\alpha$, indicating that the theoretical bound is robust with respect to the spatial resolution and damping parameter considered here.
\begin{figure}[htbp]
    \centering
    \subfloat{\includegraphics[width=0.25\linewidth]{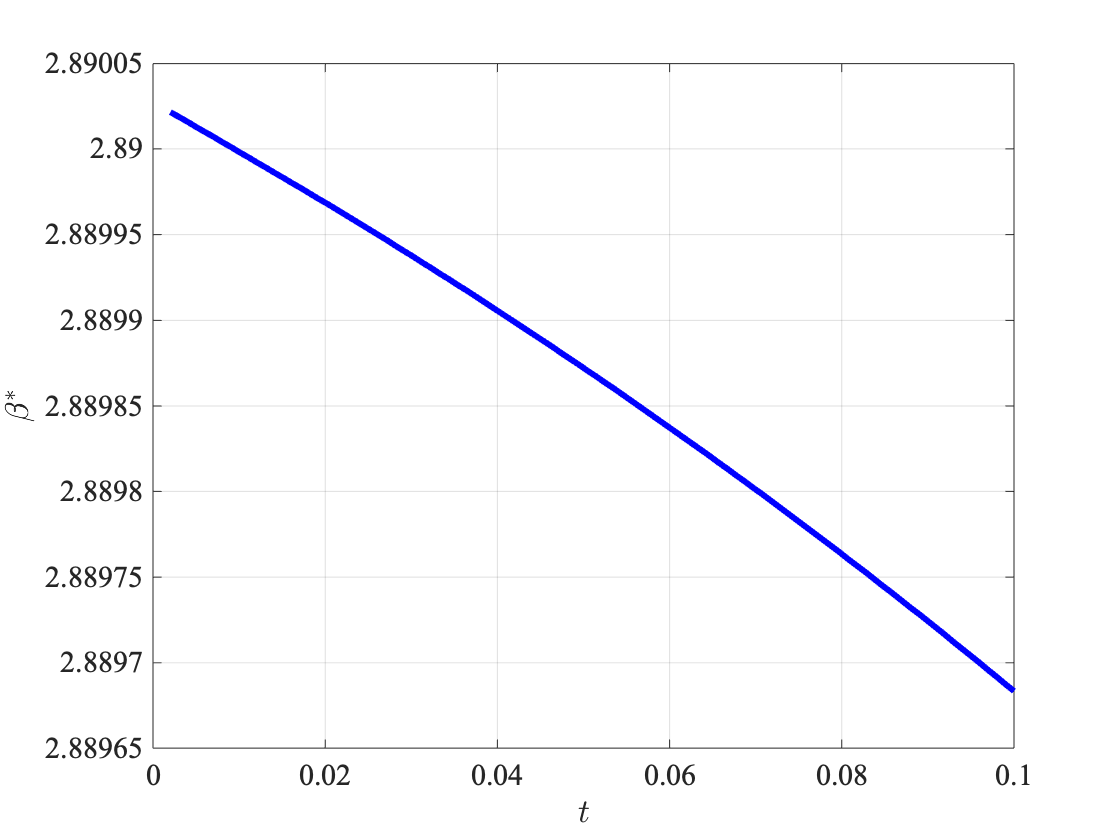}}
    \subfloat{\includegraphics[width=0.25\linewidth]{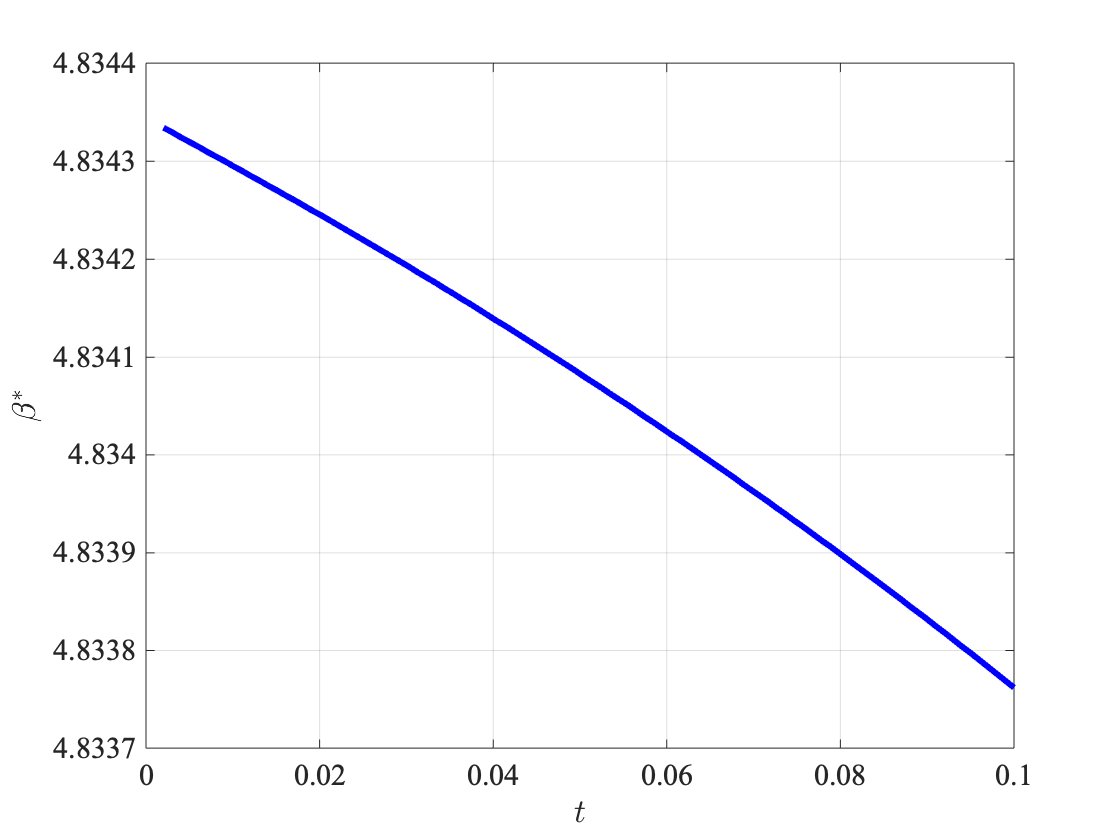}}
     \subfloat{\includegraphics[width=0.25\linewidth]{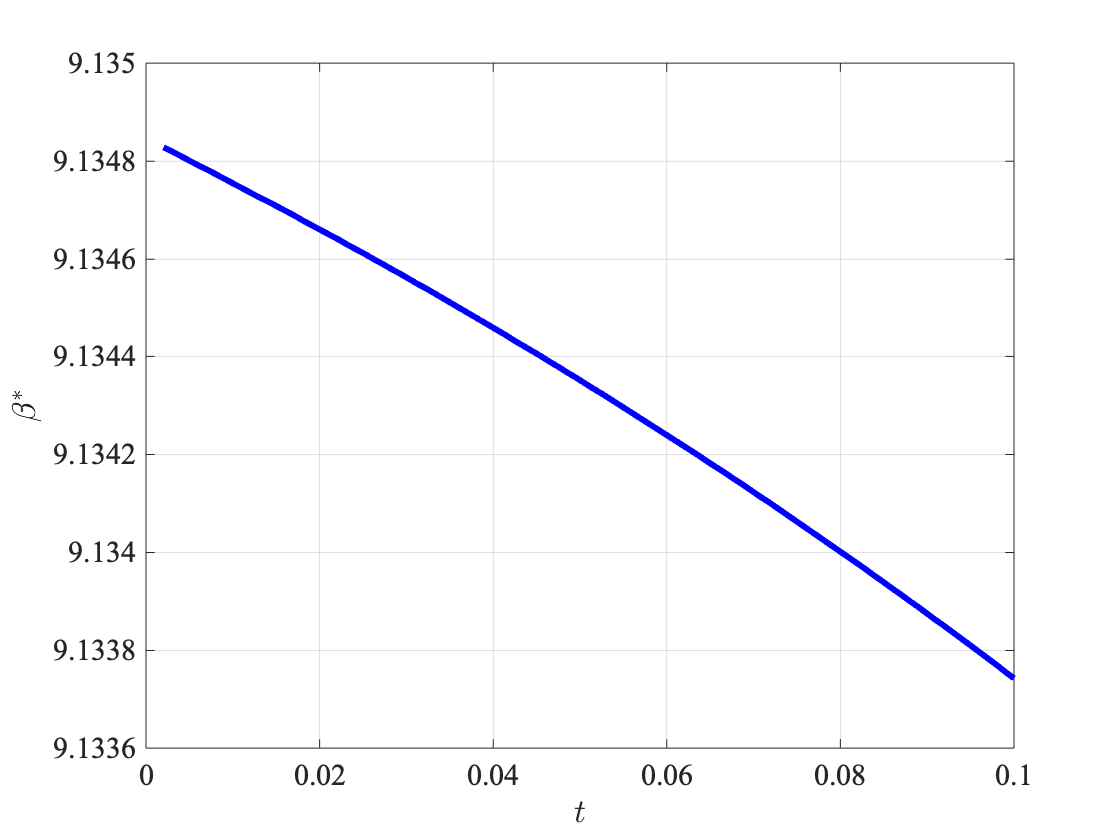}}
     \hspace{0.1in}
     \subfloat{\includegraphics[width=0.25\linewidth]{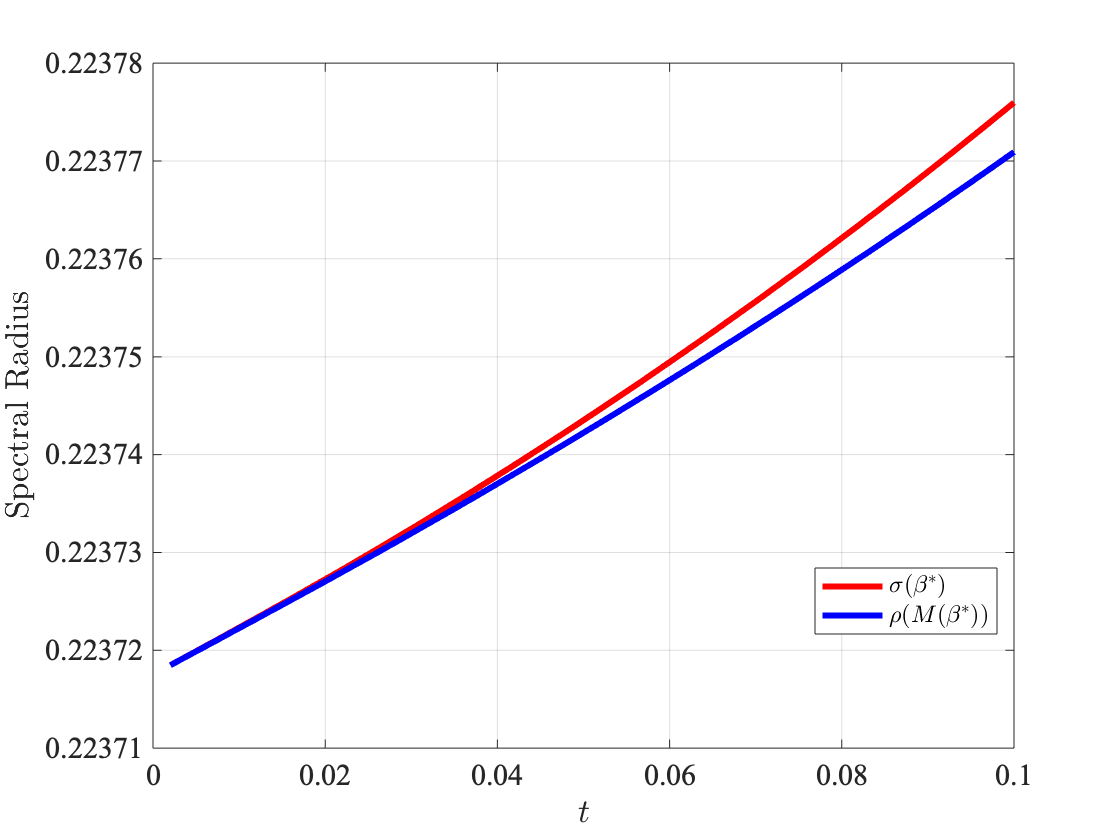}}
    \subfloat{\includegraphics[width=0.25\linewidth]{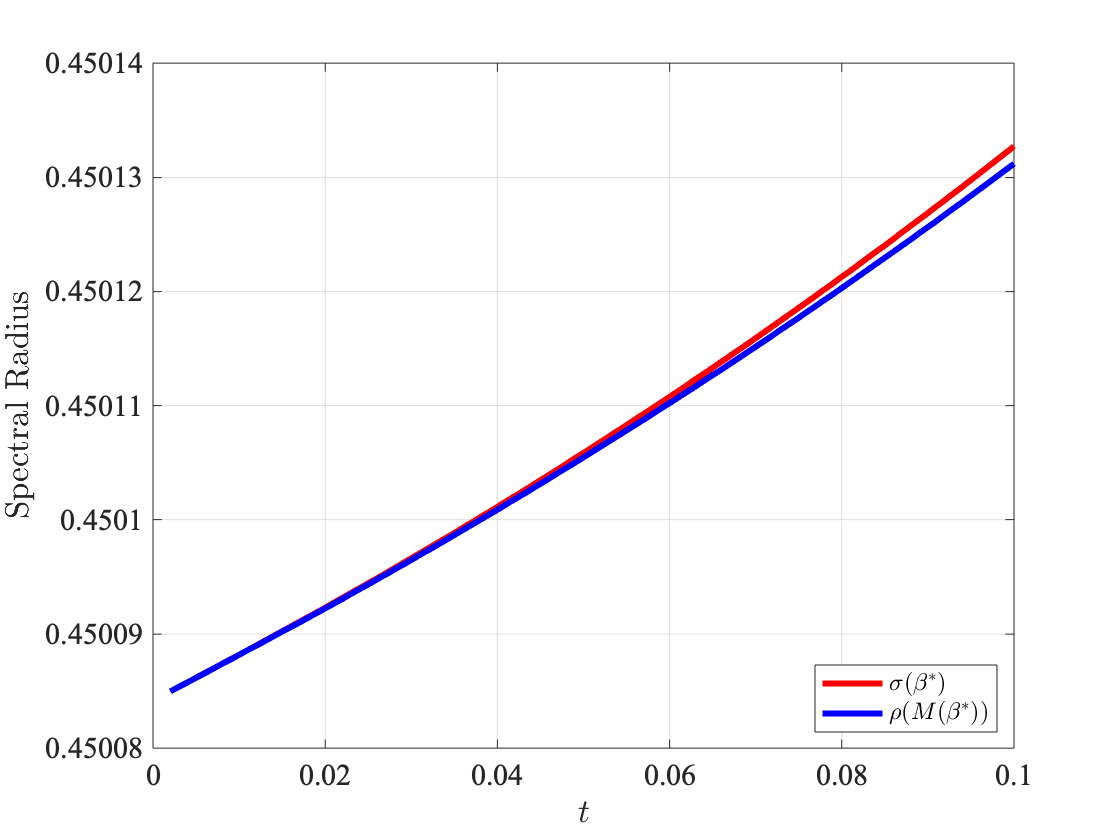}}
     \subfloat{\includegraphics[width=0.25\linewidth]{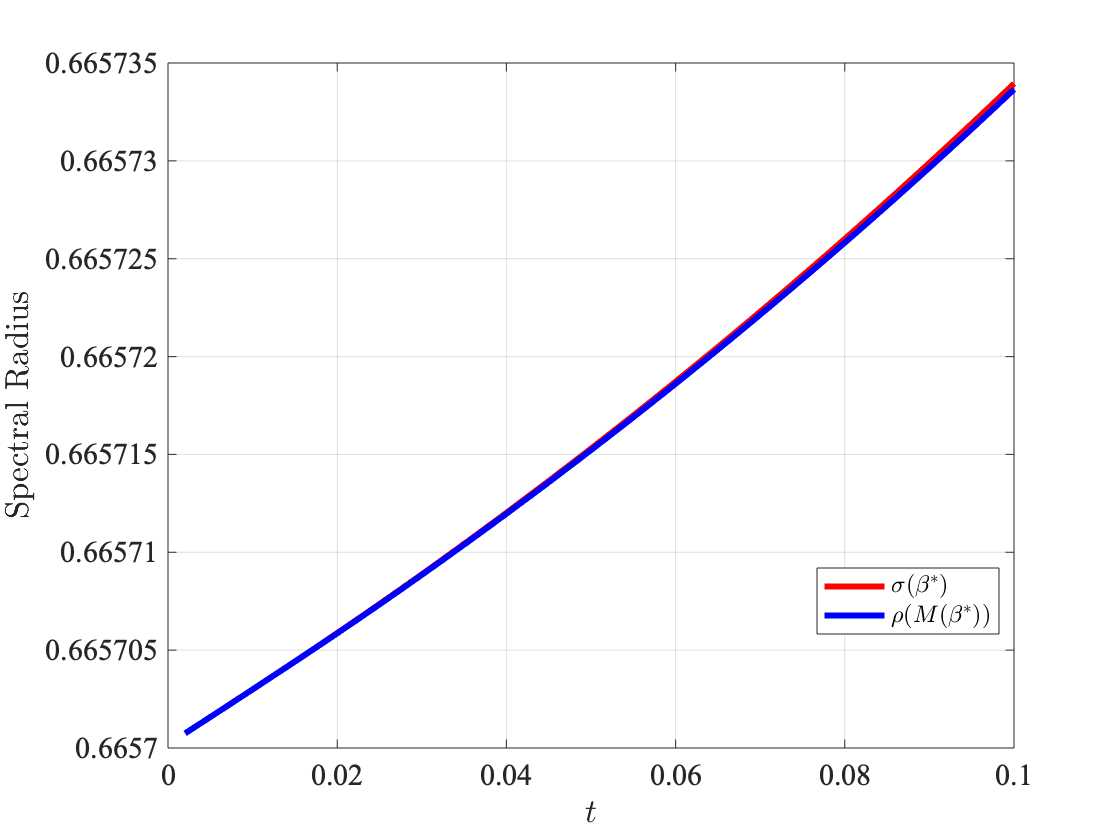}}
     \hspace{0.1in}
     \subfloat{\includegraphics[width=0.25\linewidth]{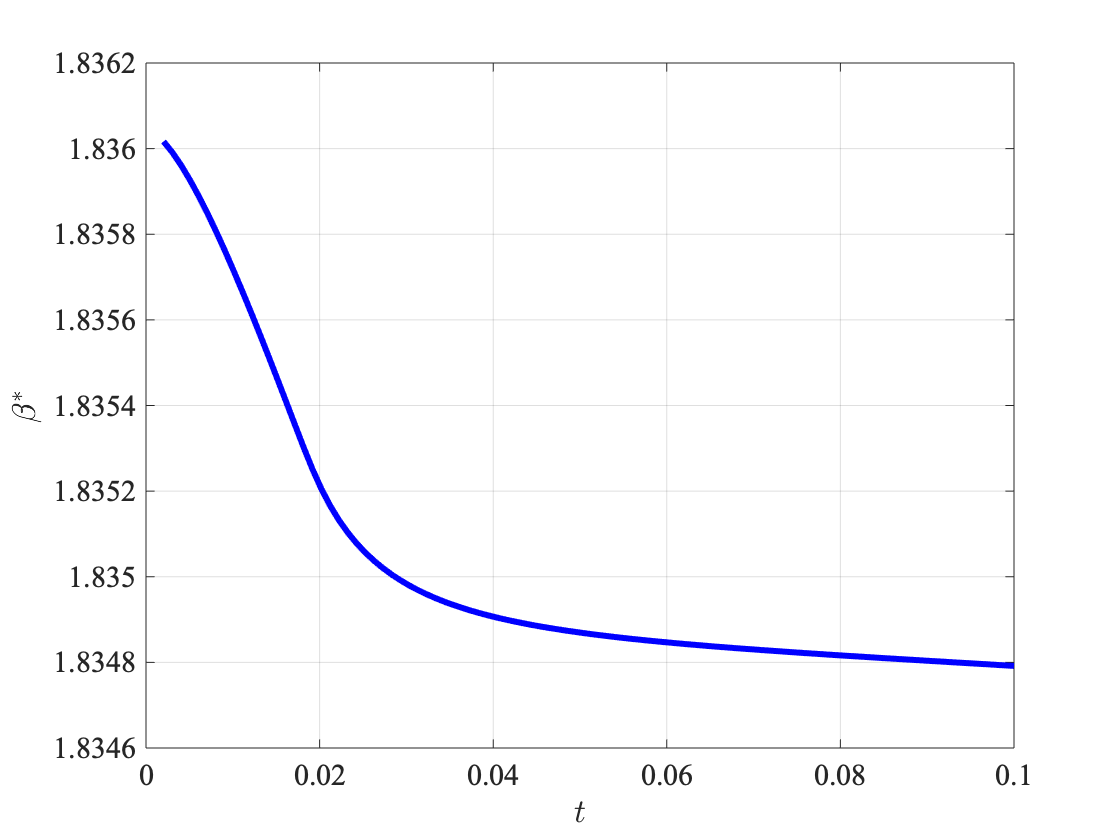}}
    \subfloat{\includegraphics[width=0.25\linewidth]{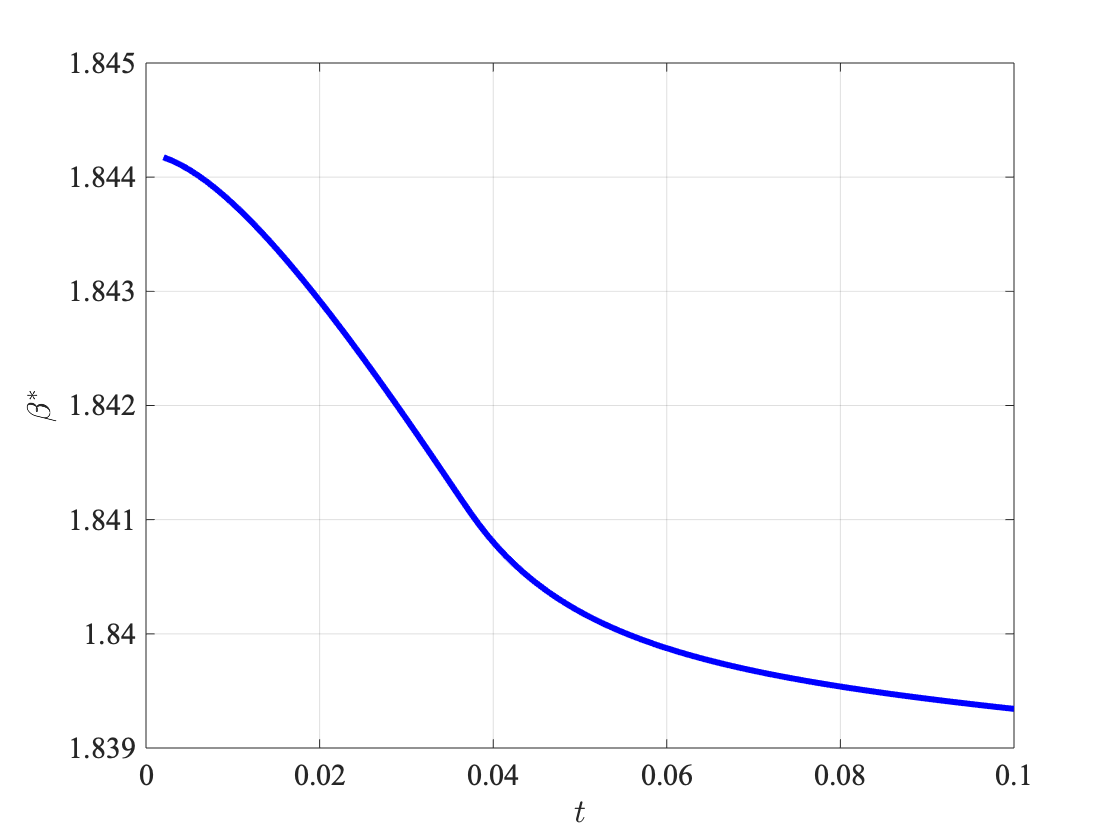}}
     \subfloat{\includegraphics[width=0.25\linewidth]{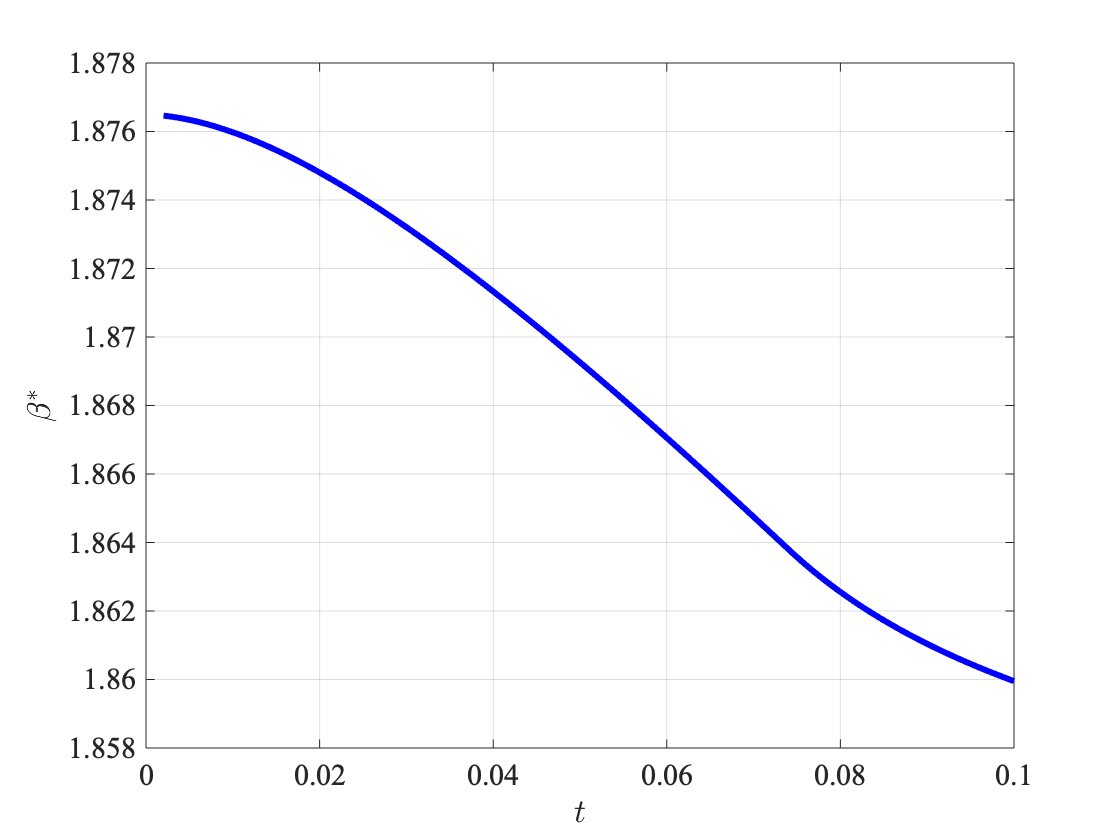}}
     \hspace{0.1in}
     \subfloat{\includegraphics[width=0.25\linewidth]{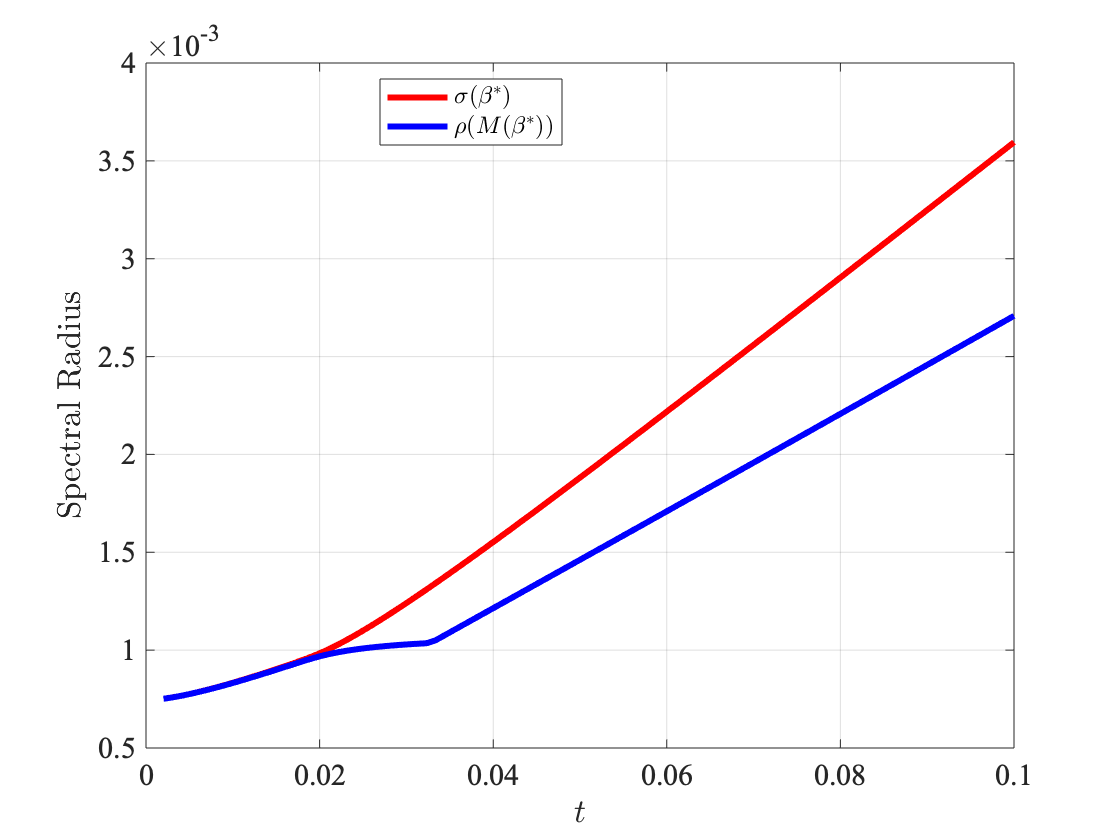}}
    \subfloat{\includegraphics[width=0.25\linewidth]{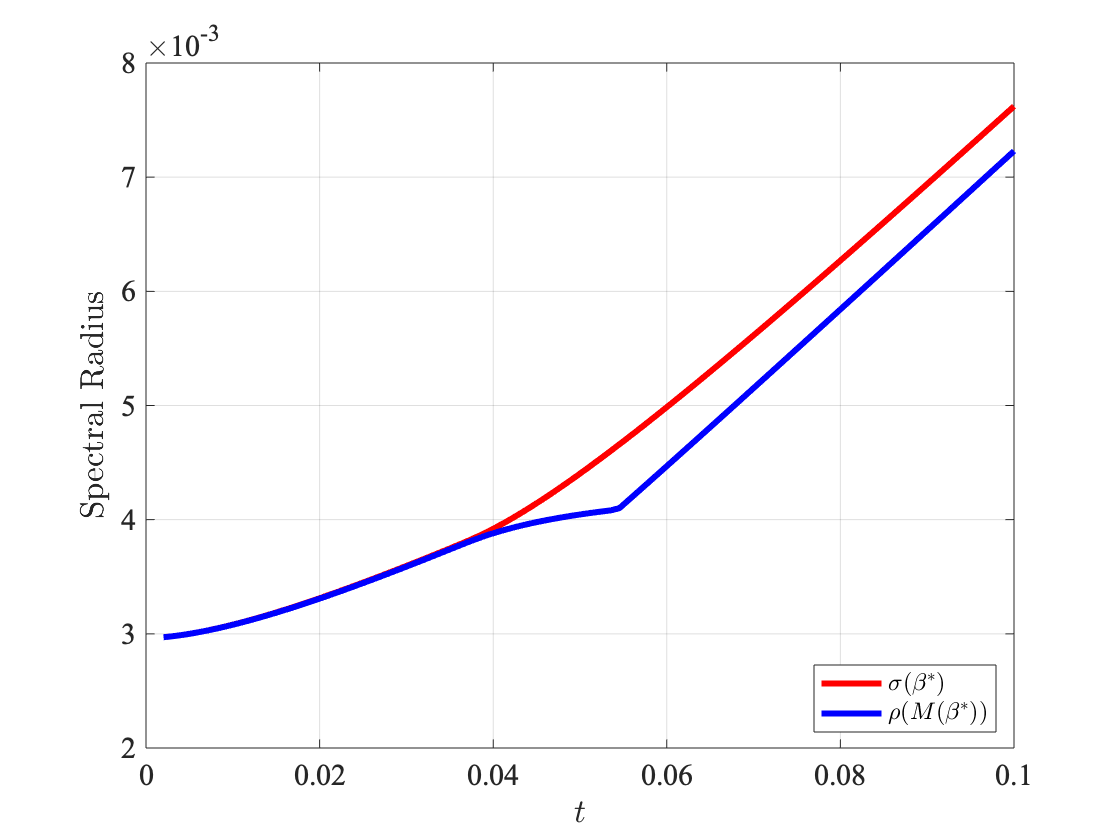}}
     \subfloat{\includegraphics[width=0.25\linewidth]{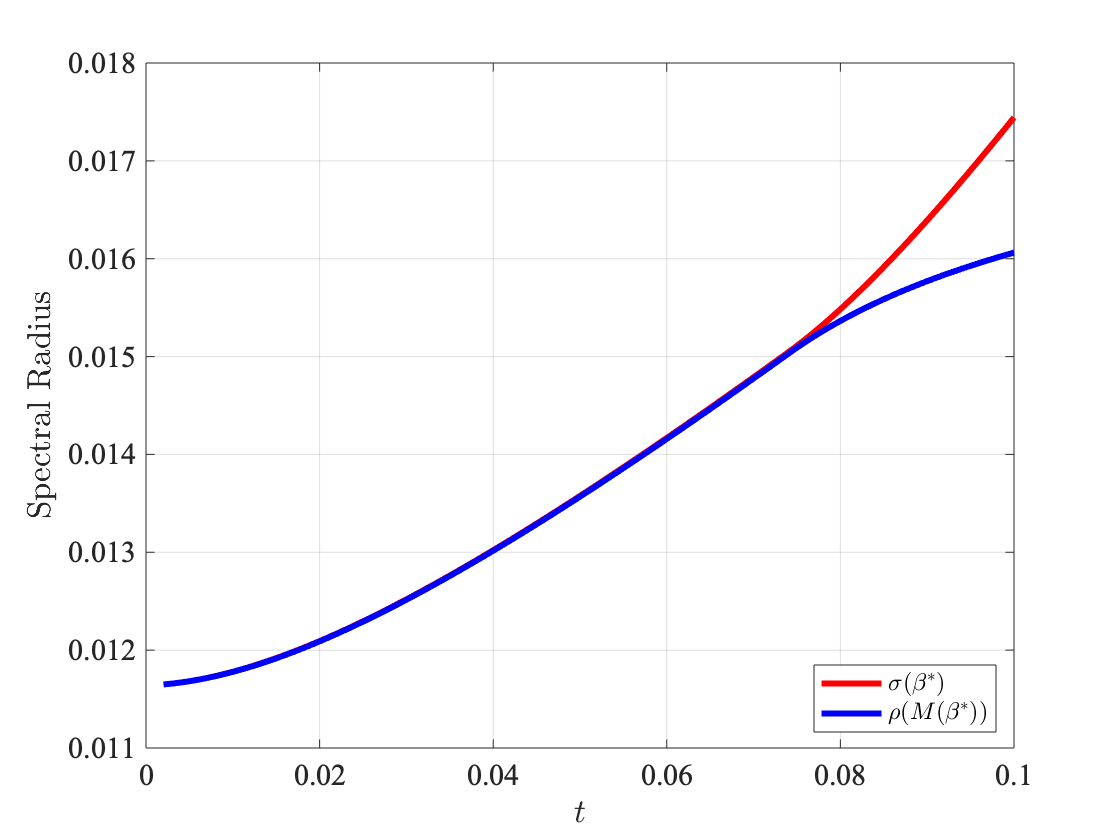}}
     \hspace{0.1in}
     \subfloat{\includegraphics[width=0.25\linewidth]{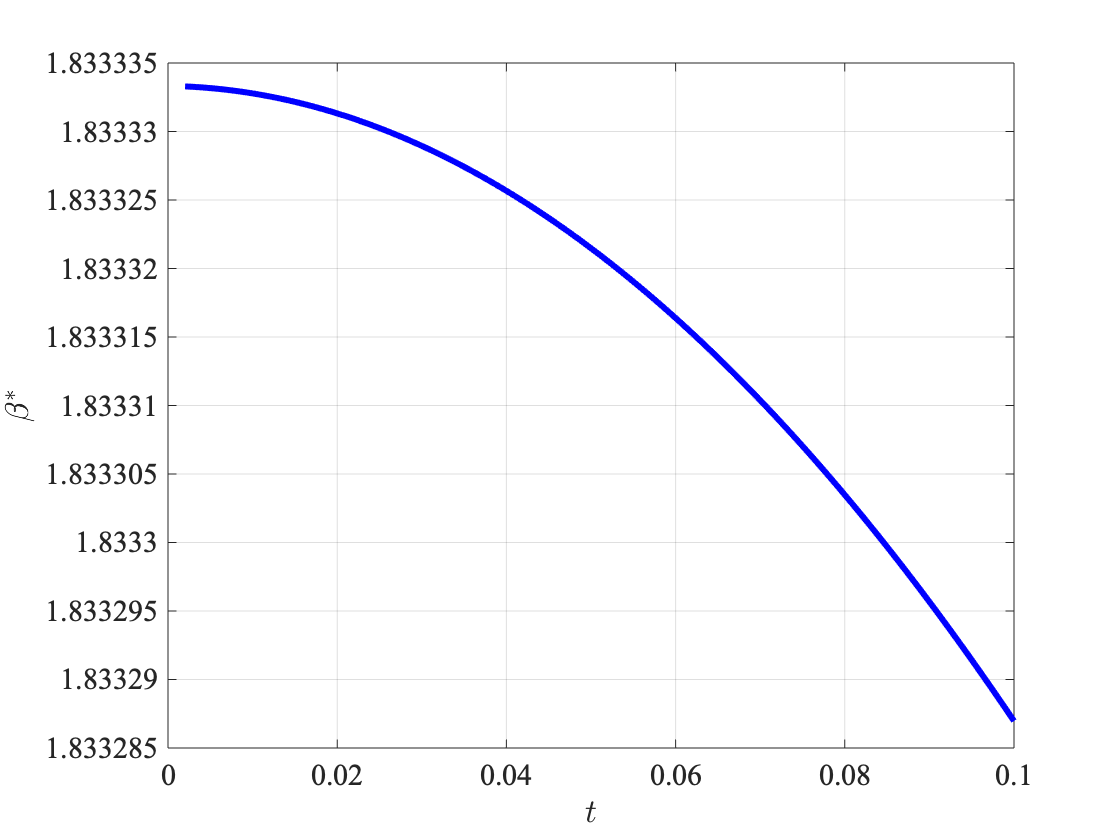}}
    \subfloat{\includegraphics[width=0.25\linewidth]{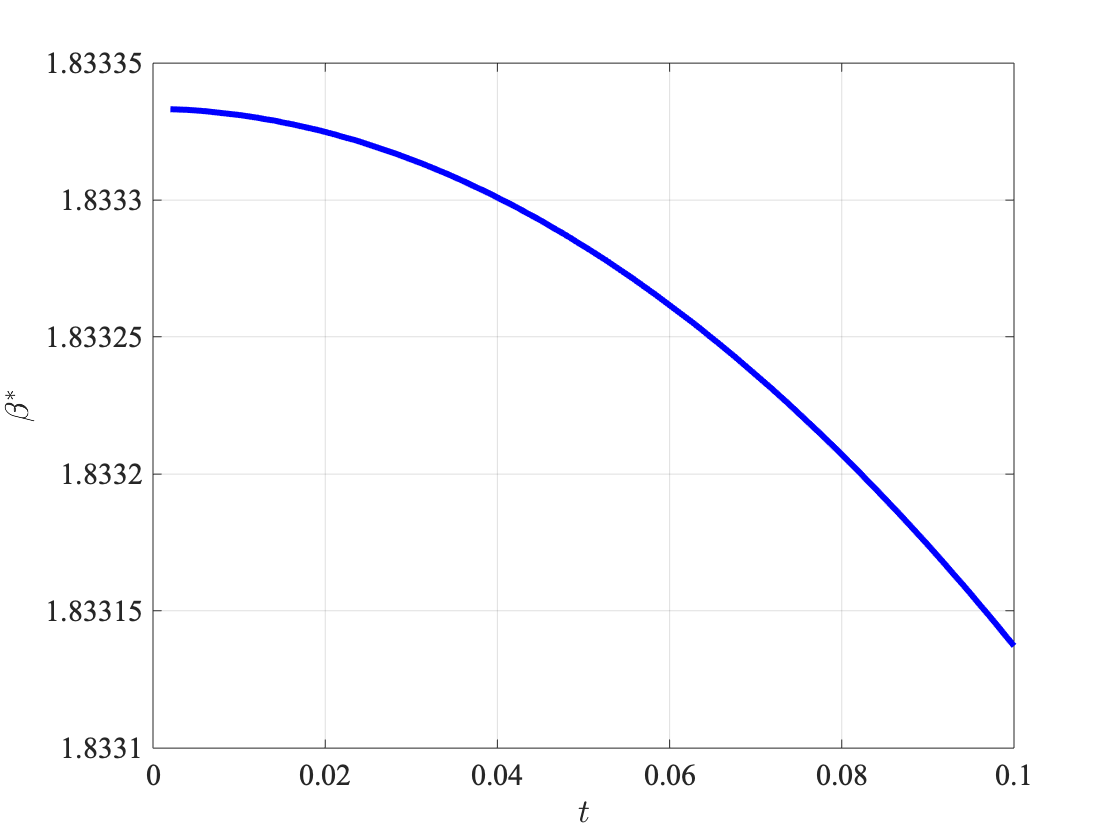}}
     \subfloat{\includegraphics[width=0.25\linewidth]{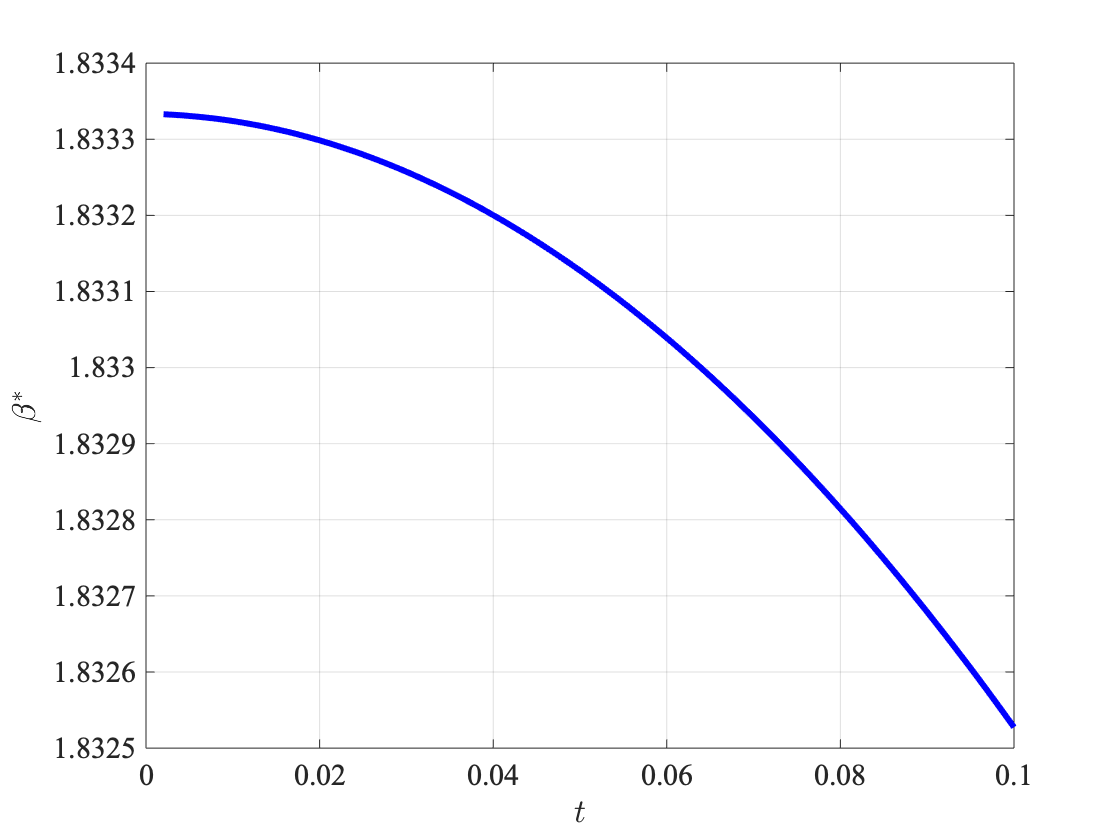}}
     \hspace{0.1in}
     \subfloat{\includegraphics[width=0.25\linewidth]{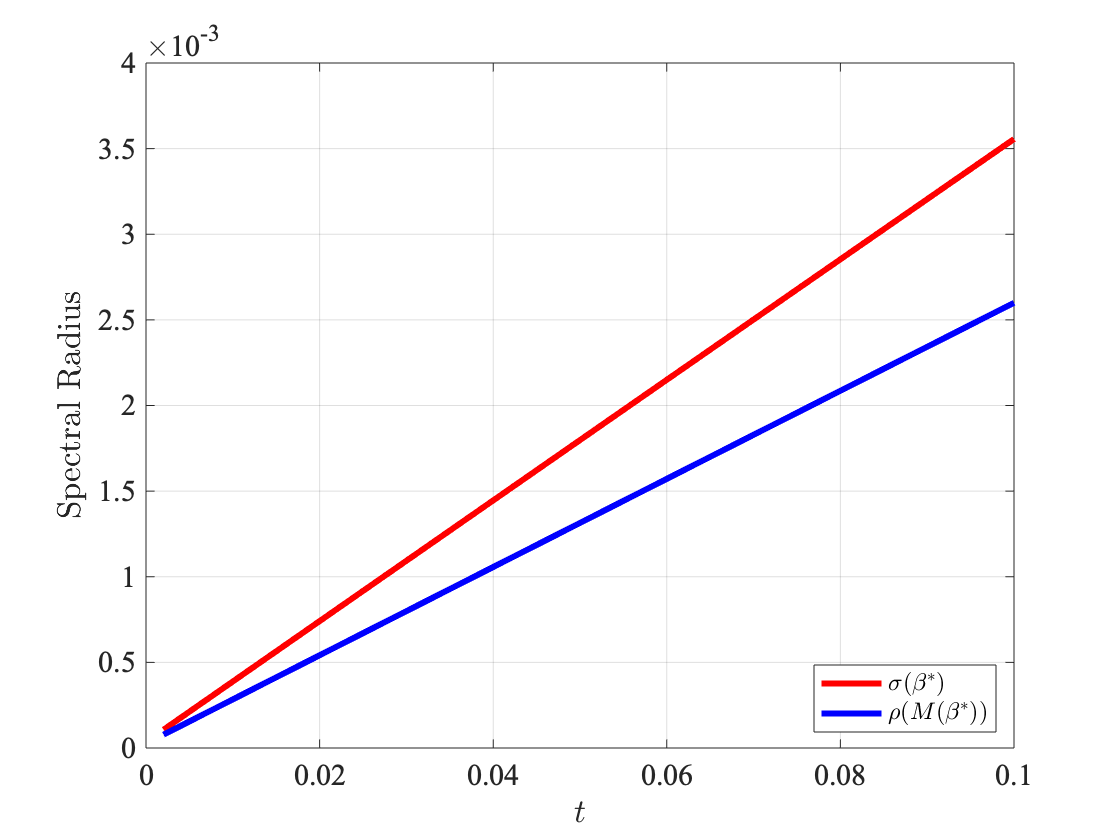}}
    \subfloat{\includegraphics[width=0.25\linewidth]{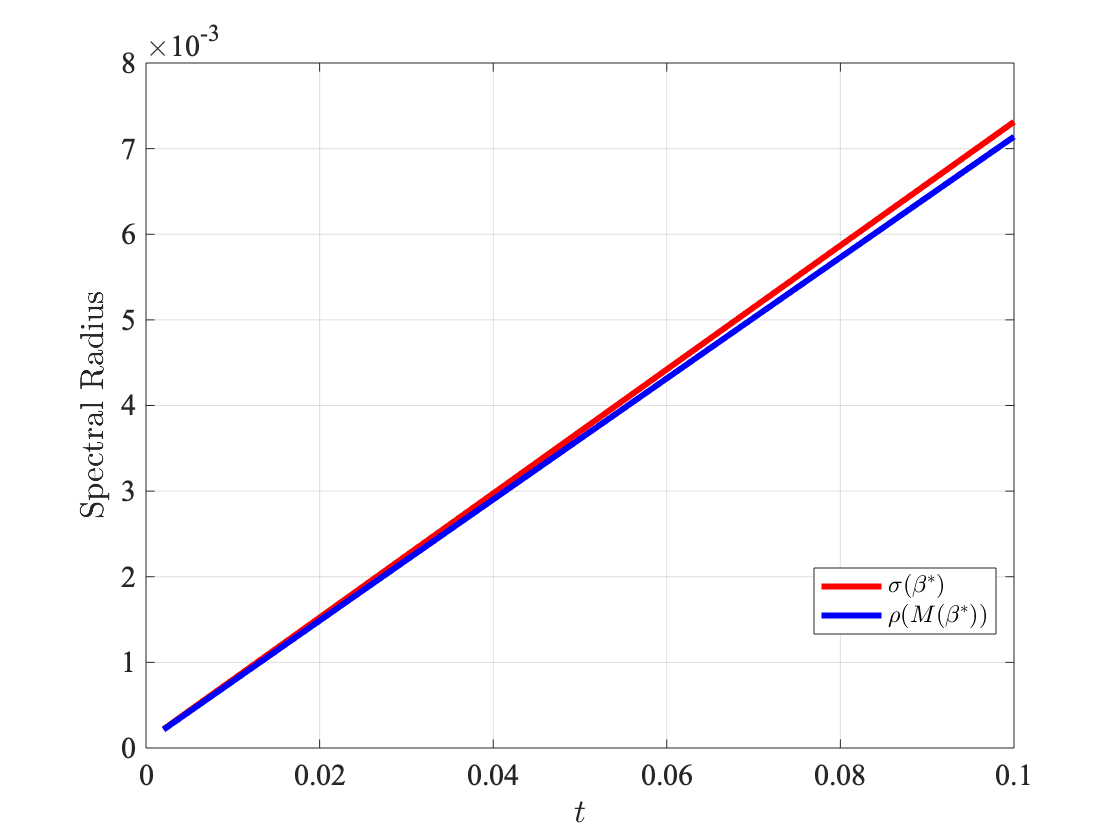}}
     \subfloat{\includegraphics[width=0.25\linewidth]{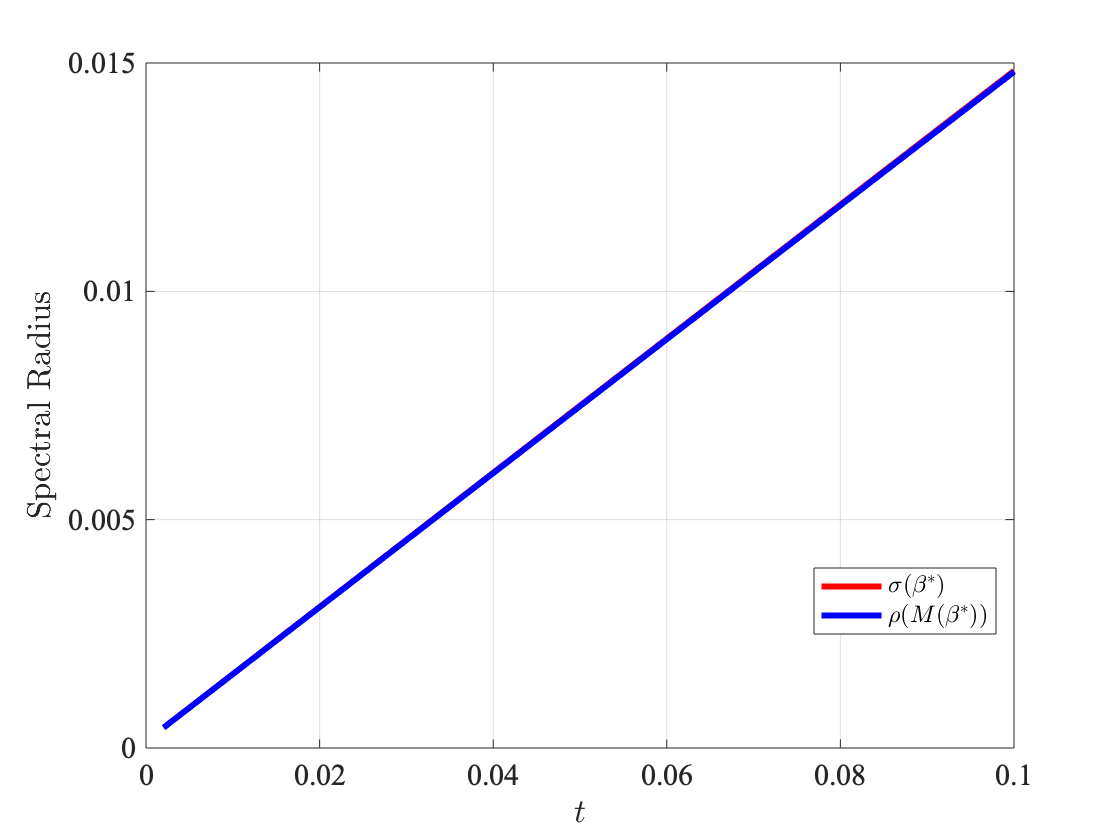}}
    \caption{The estimated $\beta^{*}$, the spectral radius of $\rho(M(\beta^{*}))$ and its upper bound $\sigma(\beta^*)$ over time for scheme B in 1D. From left to right $N_x=32,64,128$ and $N_t=100$. First two rows: $\alpha=0.5$; Middle two rows: $\alpha=0.001$; Bottom two rows: $\alpha=0$}
    \label{fig:2}
\end{figure}
A similar behavior is observed for scheme B in Figure~\ref{fig:2}. The spectral radius remains below its corresponding upper bound for all tested values of $N_x$ and $\alpha$, and stays below one throughout the simulation, further supporting the theoretical spectral-radius estimate.

\begin{figure}[htbp]
    \centering
    \subfloat{\includegraphics[width=0.25\linewidth]{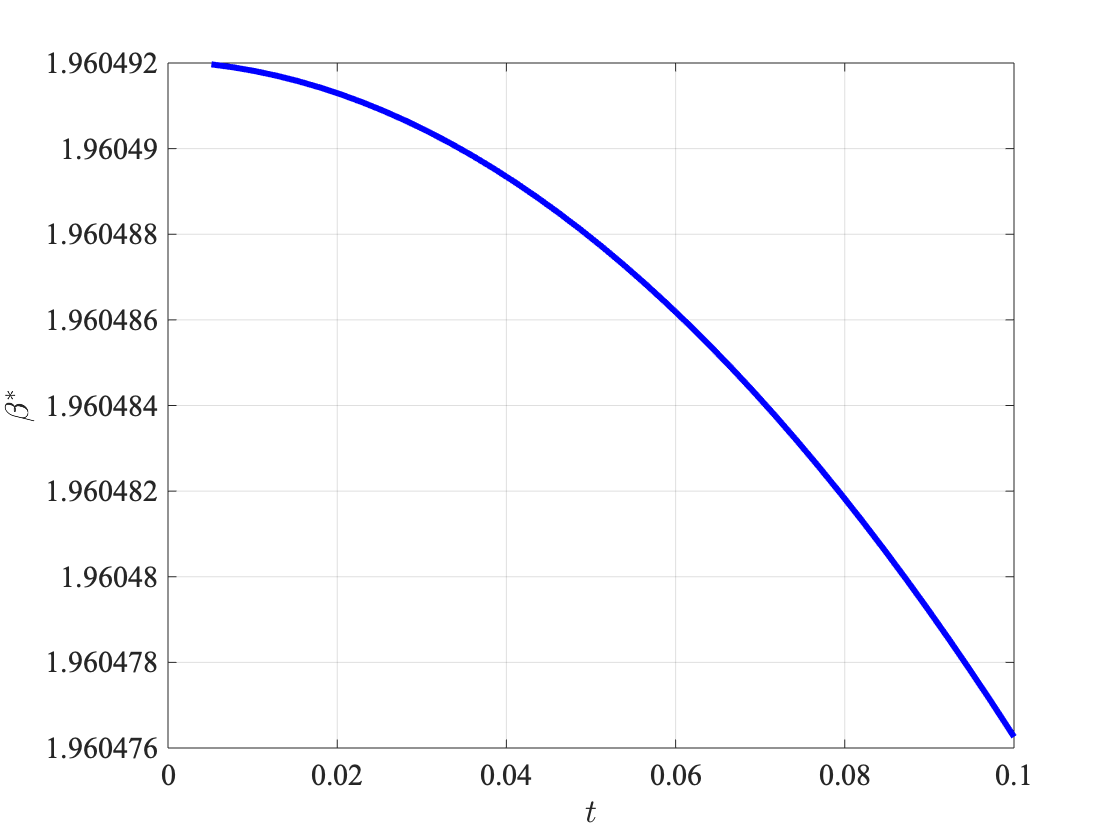}}
    \subfloat{\includegraphics[width=0.25\linewidth]{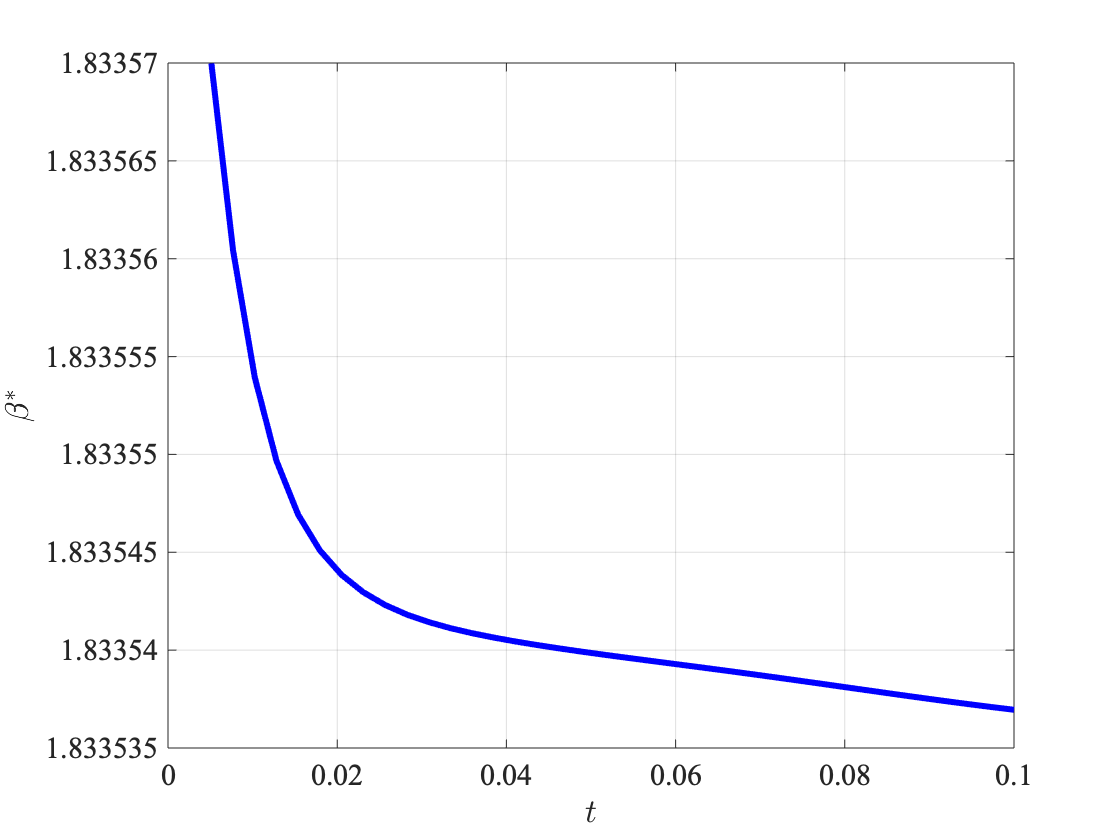}}
    \subfloat{\includegraphics[width=0.25\linewidth]{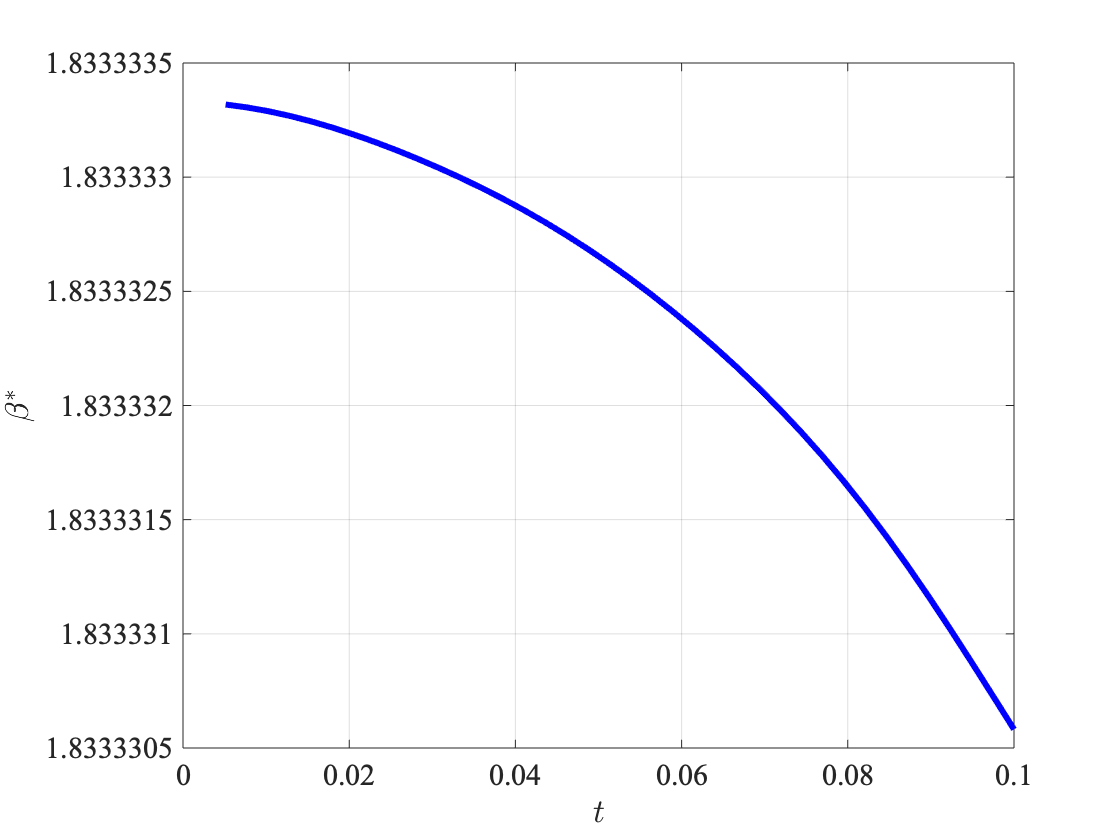}}
    \hspace{0.1in}
    \subfloat{\includegraphics[width=0.25\linewidth]{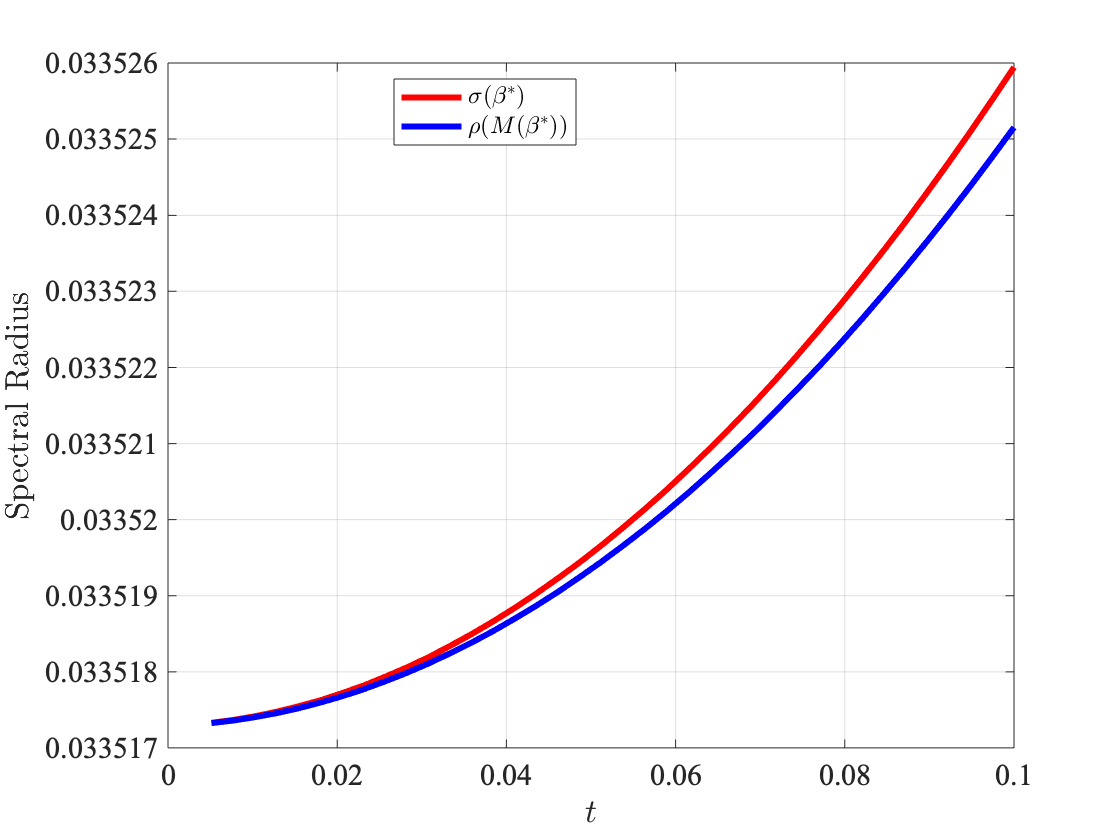}}
    \subfloat{\includegraphics[width=0.25\linewidth]{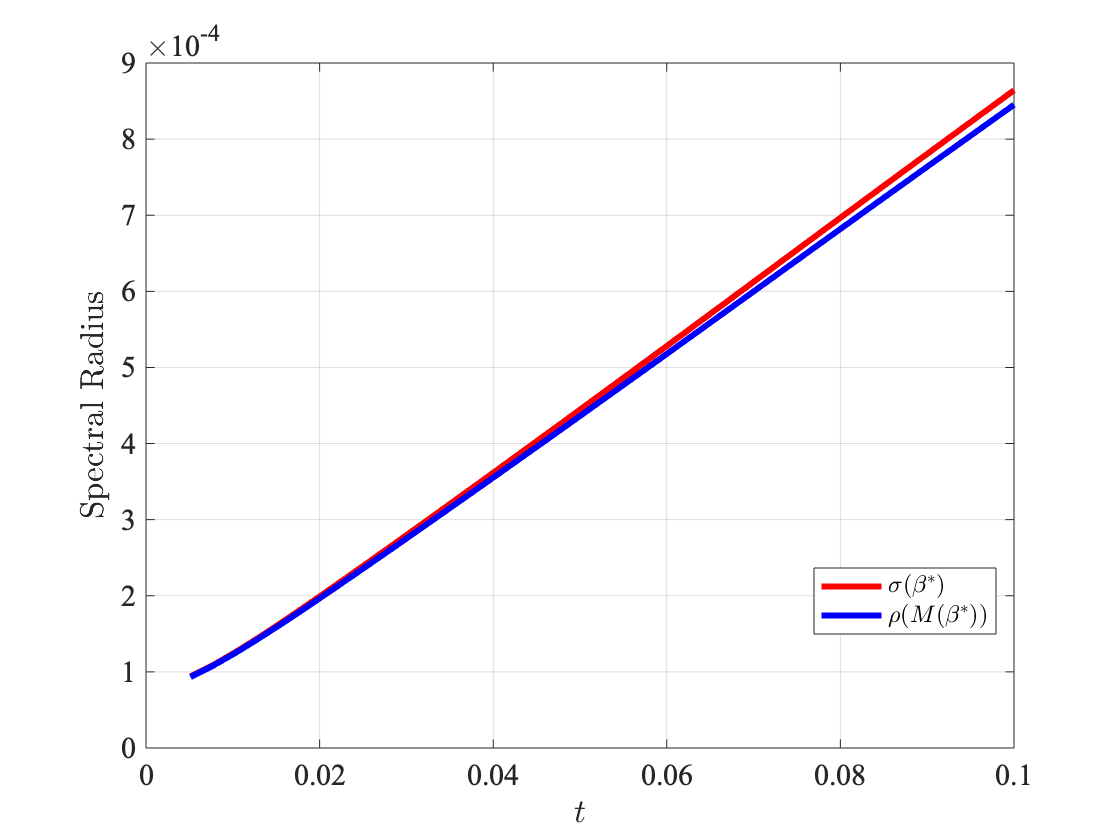}}
    \subfloat{\includegraphics[width=0.25\linewidth]{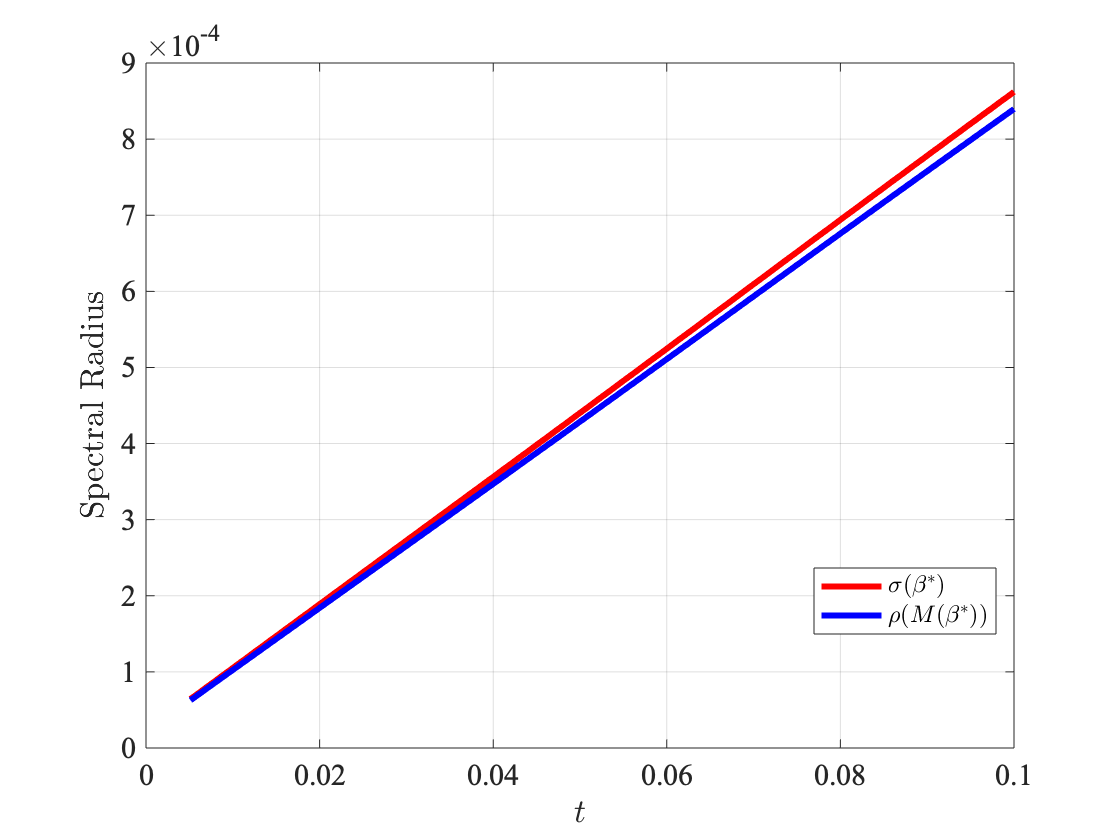}}
    \hspace{0.1in}
    \subfloat{\includegraphics[width=0.25\linewidth]{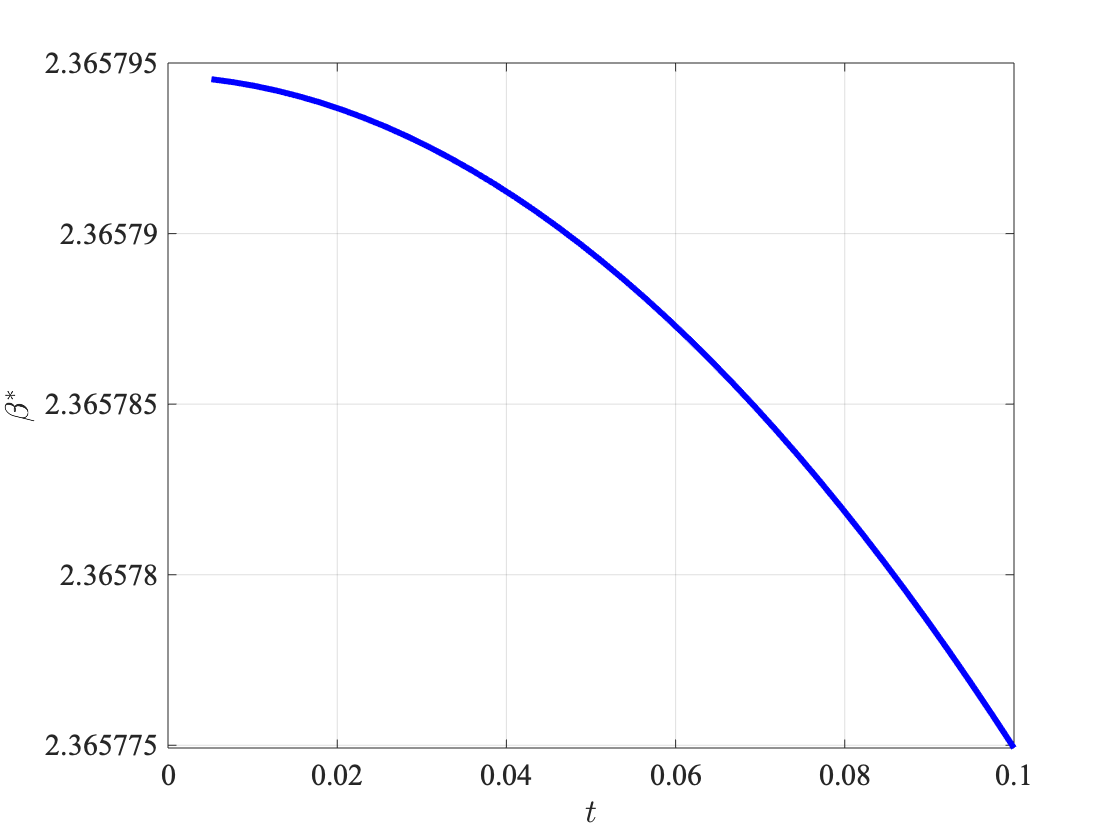}}
    \subfloat{\includegraphics[width=0.25\linewidth]{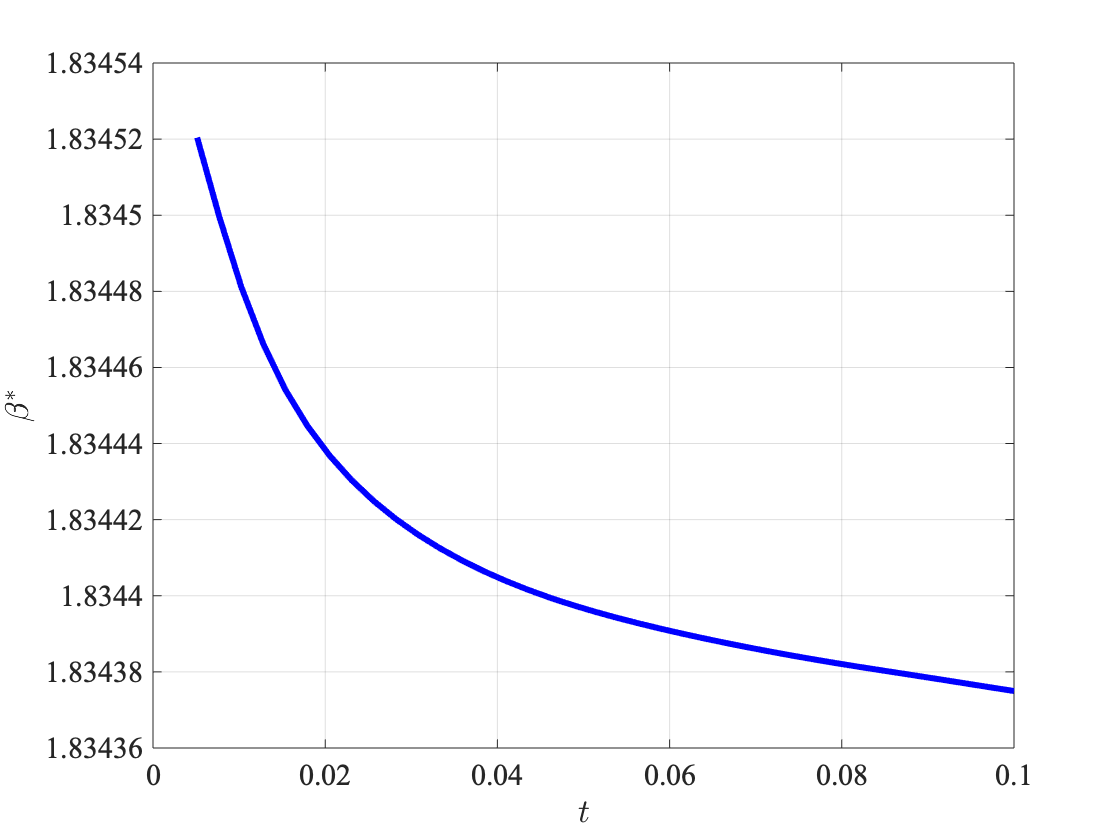}}
    \subfloat{\includegraphics[width=0.25\linewidth]{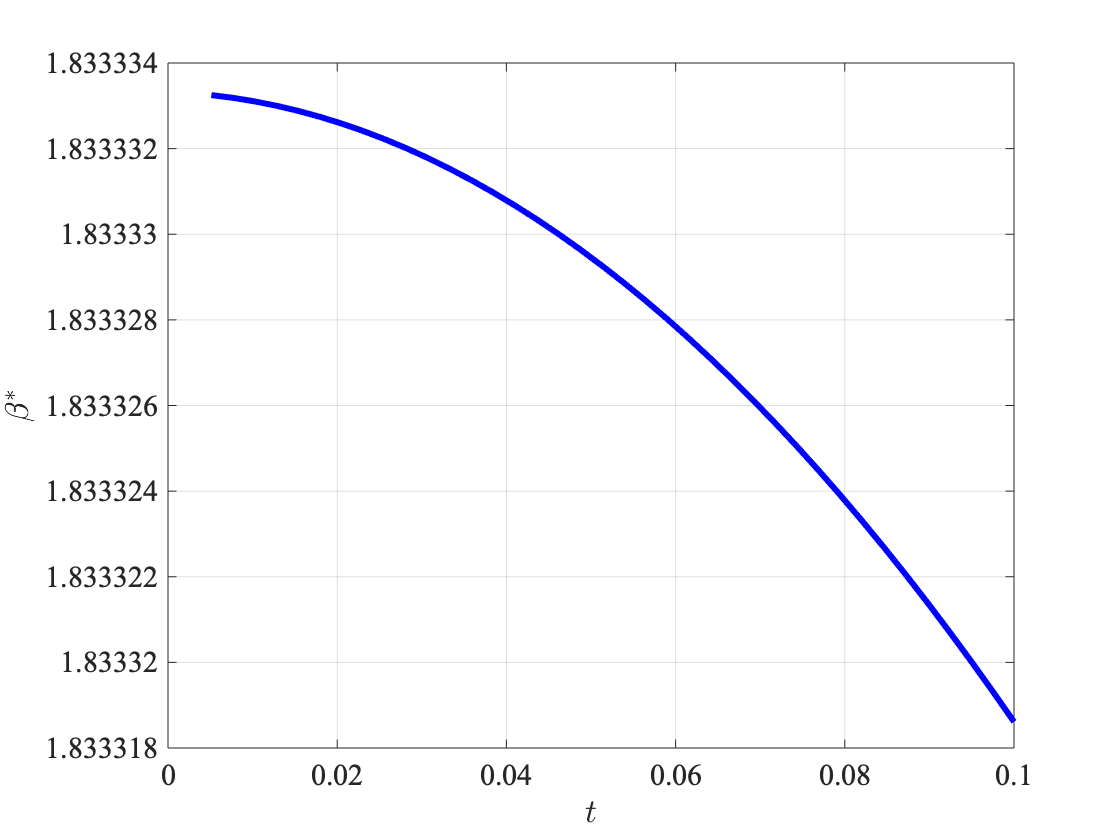}}
    \hspace{0.1in}
    \subfloat{\includegraphics[width=0.25\linewidth]{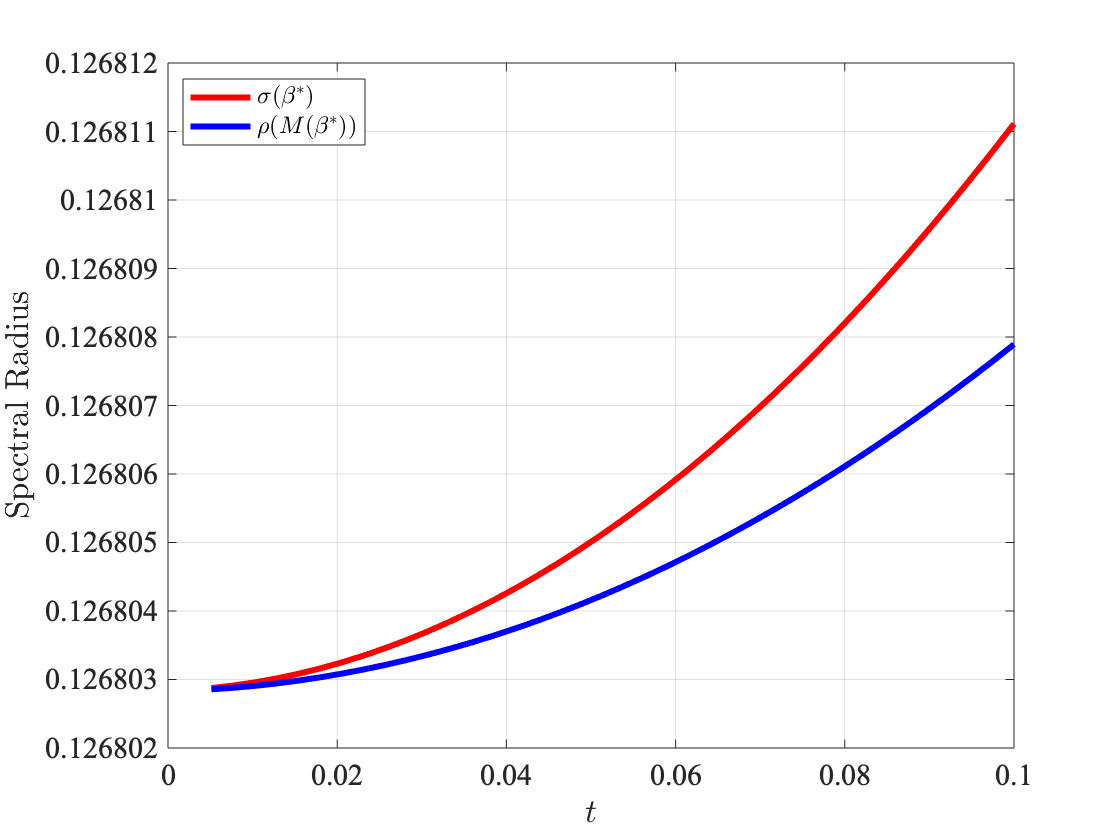}}
    \subfloat{\includegraphics[width=0.25\linewidth]{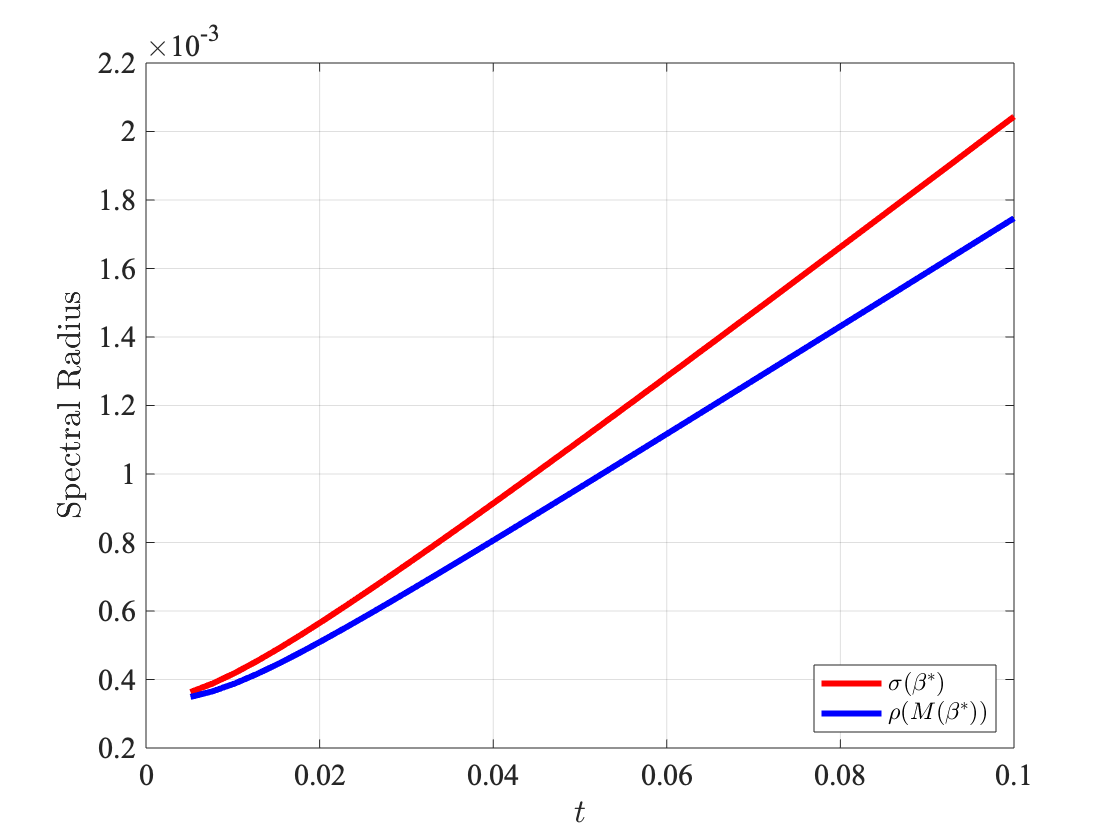}}
    \subfloat{\includegraphics[width=0.25\linewidth]{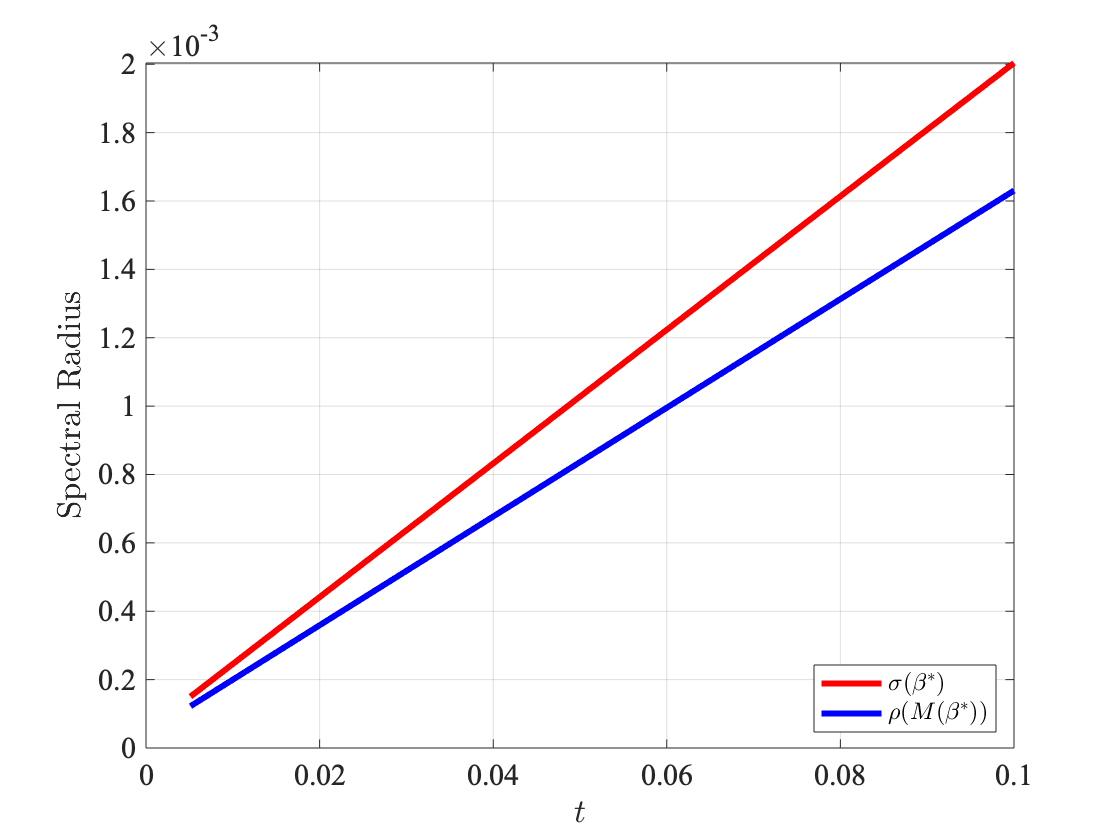}}
\caption{The estimated $\beta^{*}$, the spectral radius of $\rho(M(\beta^{*}))$ and its upper bound $\sigma(\beta^*)$ over time for scheme A in 3D. Top two rows: $N_x=N_y=N_z=4$, $N_t=40$; Bottom two rows: $N_x=N_y=N_z=8$, $N_t=40$. From left to right: $\alpha=0.5,0.001,0$.}
    \label{fig:3}
\end{figure}
\begin{figure}[htbp]
    \centering
    \subfloat{\includegraphics[width=0.25\linewidth]{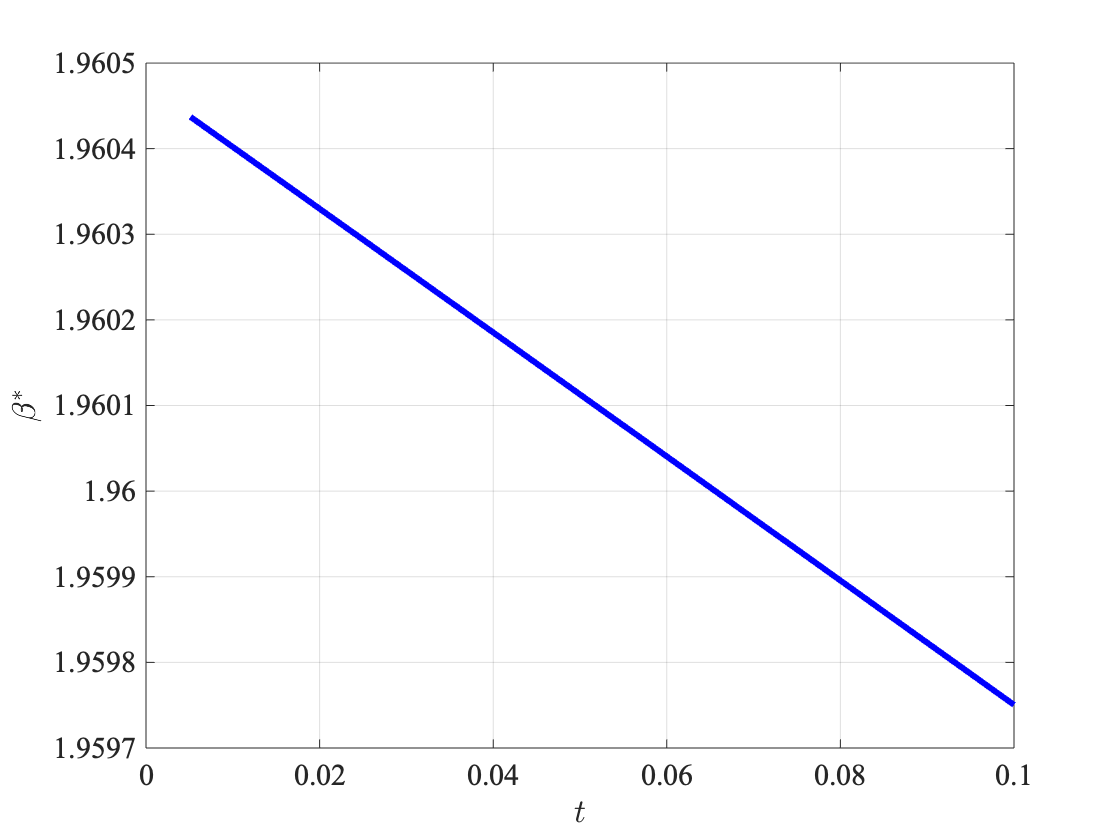}}
    \subfloat{\includegraphics[width=0.25\linewidth]{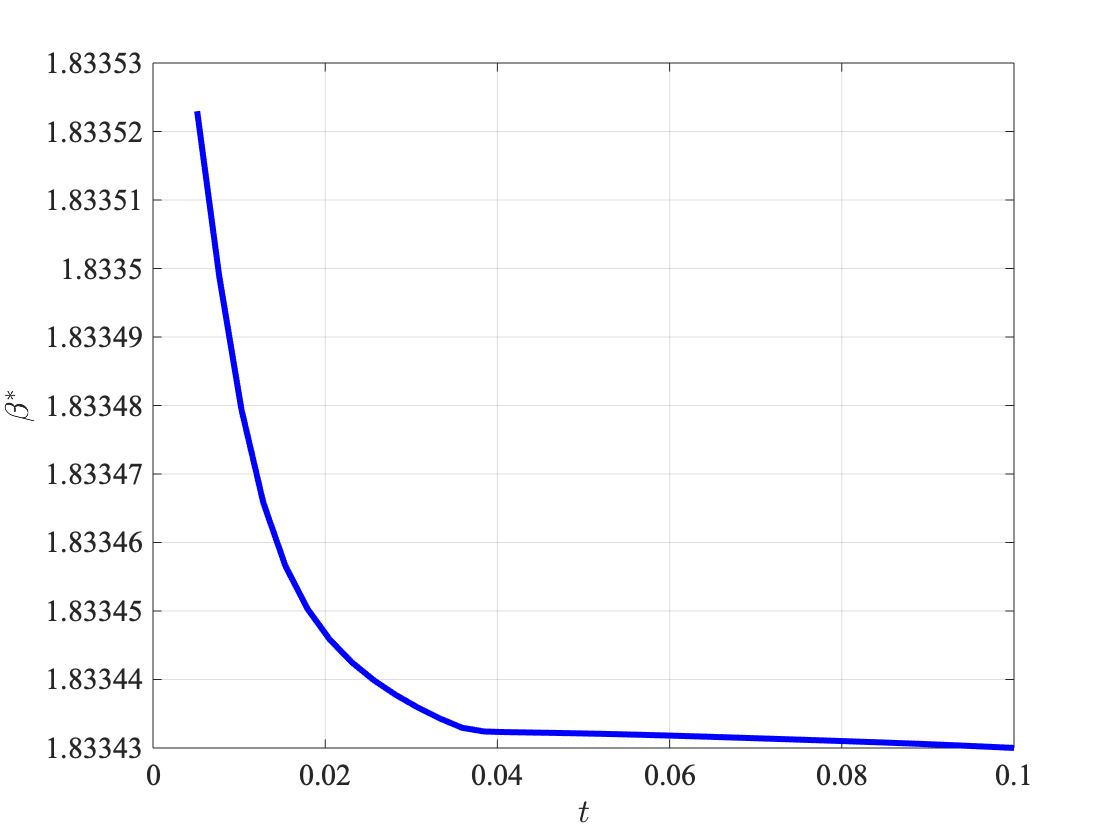}}
    \subfloat{\includegraphics[width=0.25\linewidth]{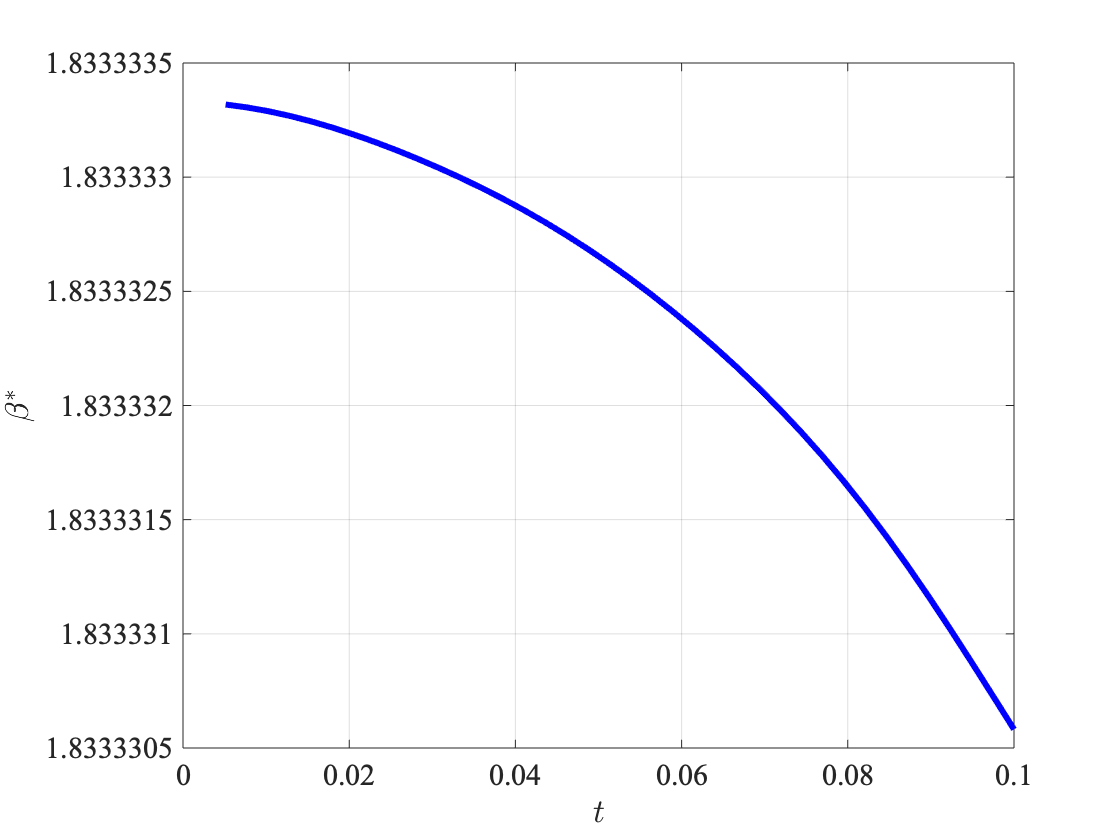}}
    \hspace{0.1in}
    \subfloat{\includegraphics[width=0.25\linewidth]{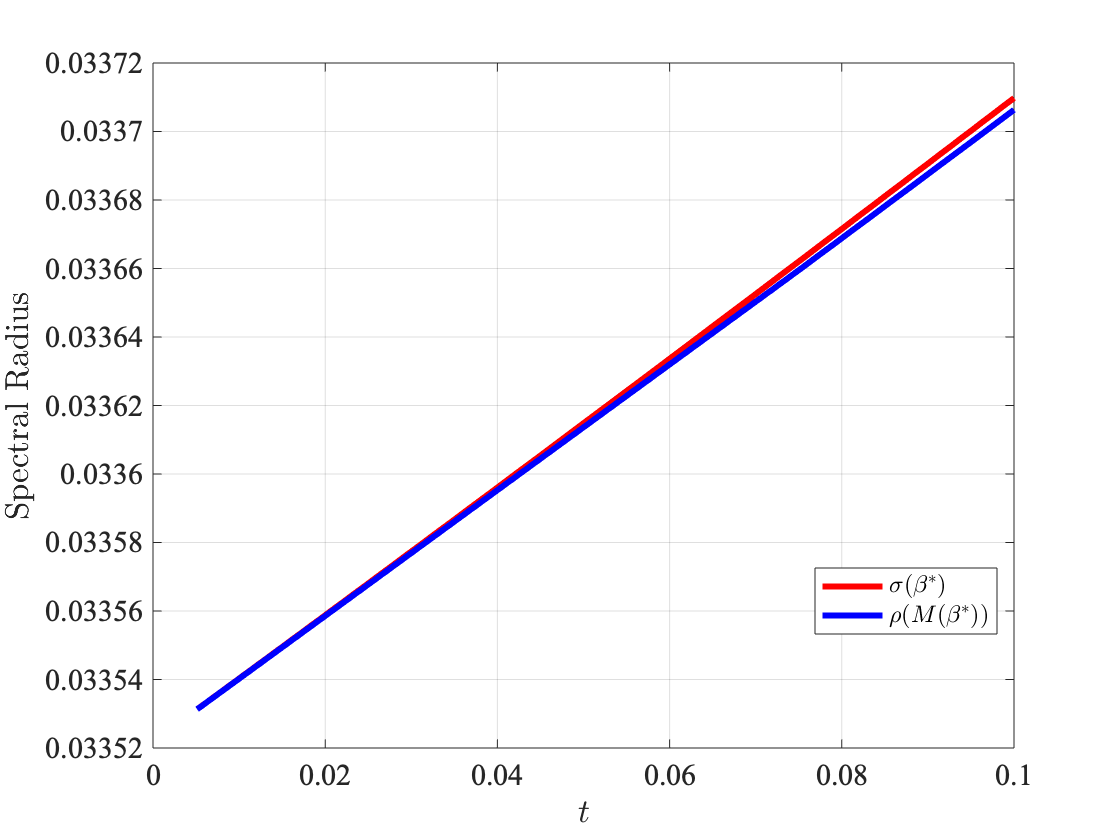}}
    \subfloat{\includegraphics[width=0.25\linewidth]{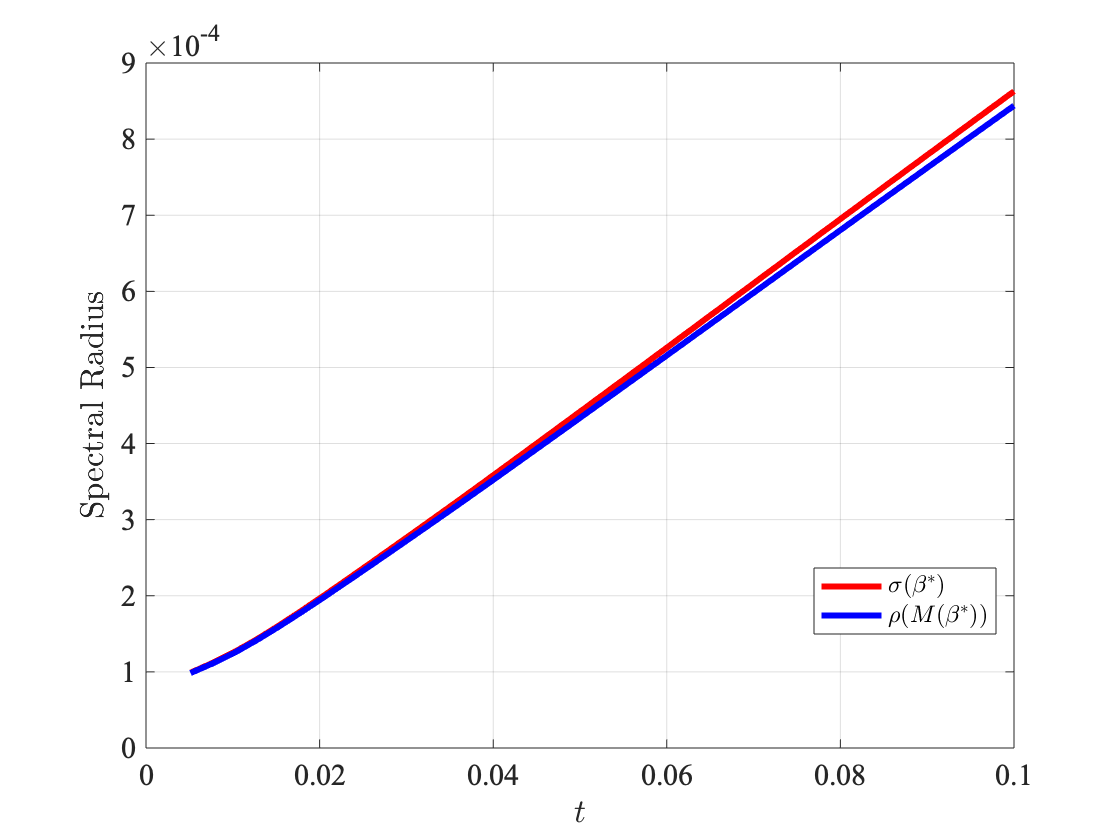}}
    \subfloat{\includegraphics[width=0.25\linewidth]{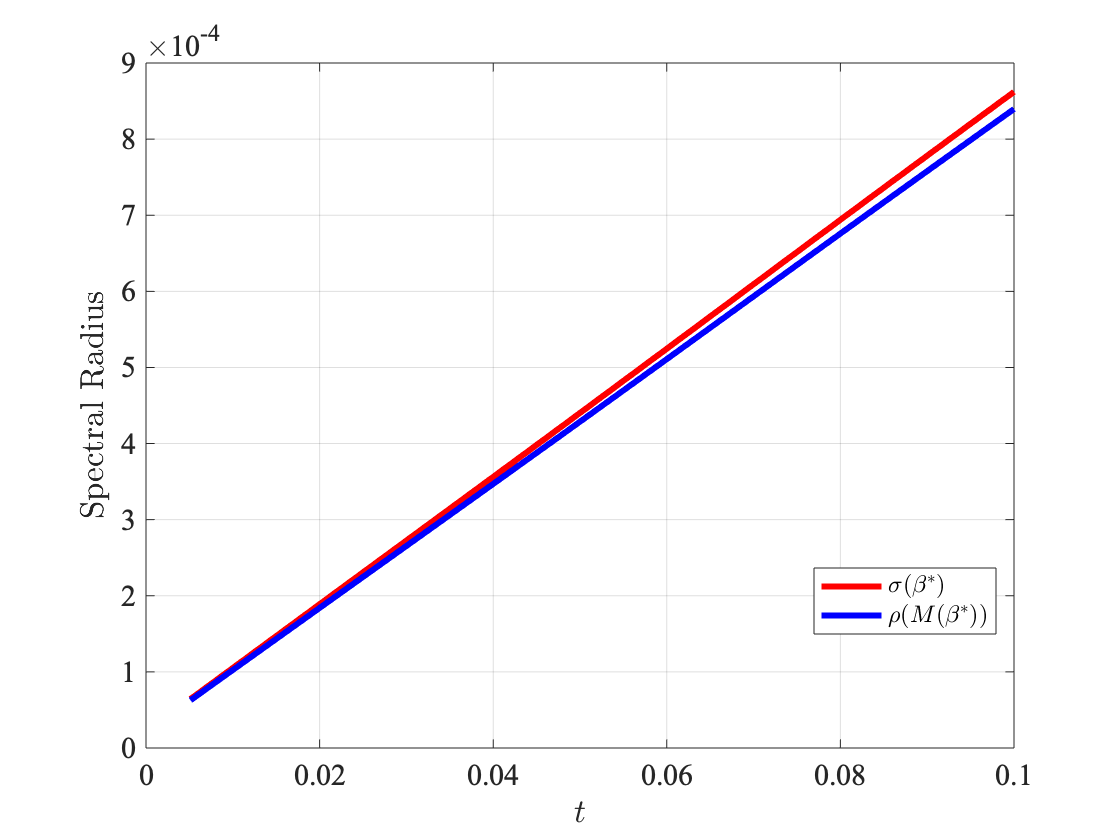}}
    \hspace{0.1in}
    \subfloat{\includegraphics[width=0.25\linewidth]{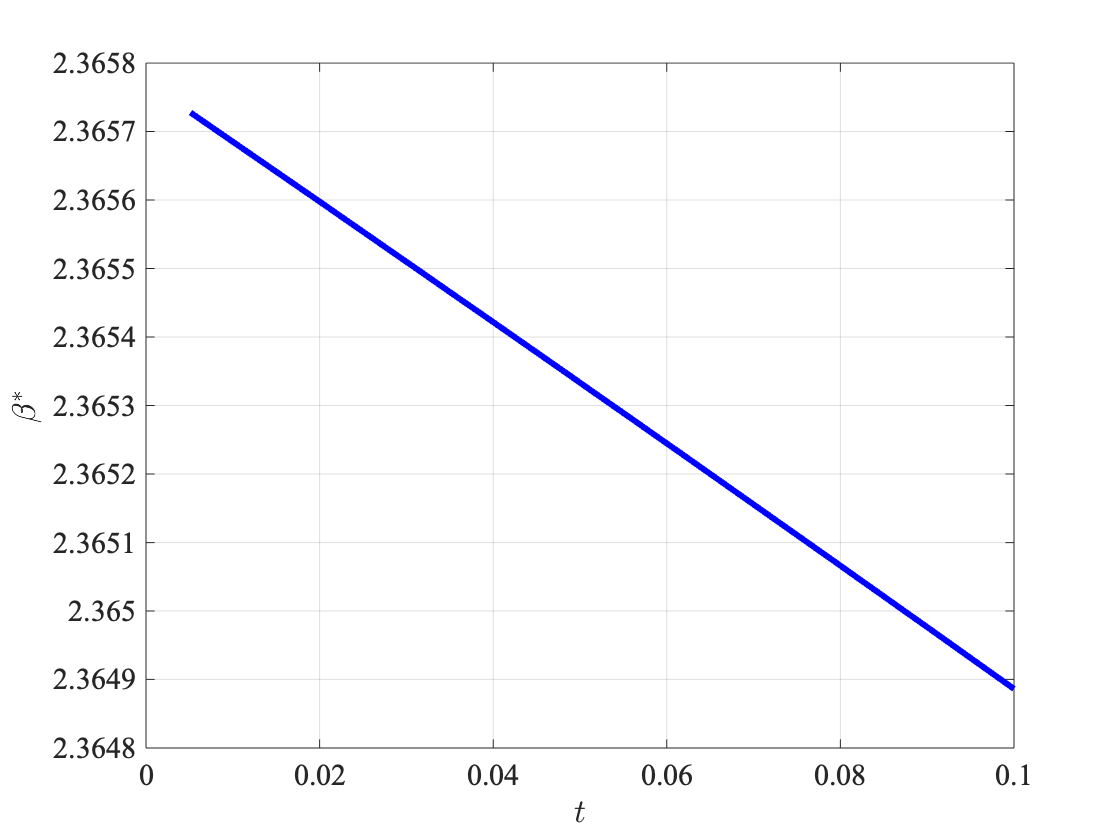}}
    \subfloat{\includegraphics[width=0.25\linewidth]{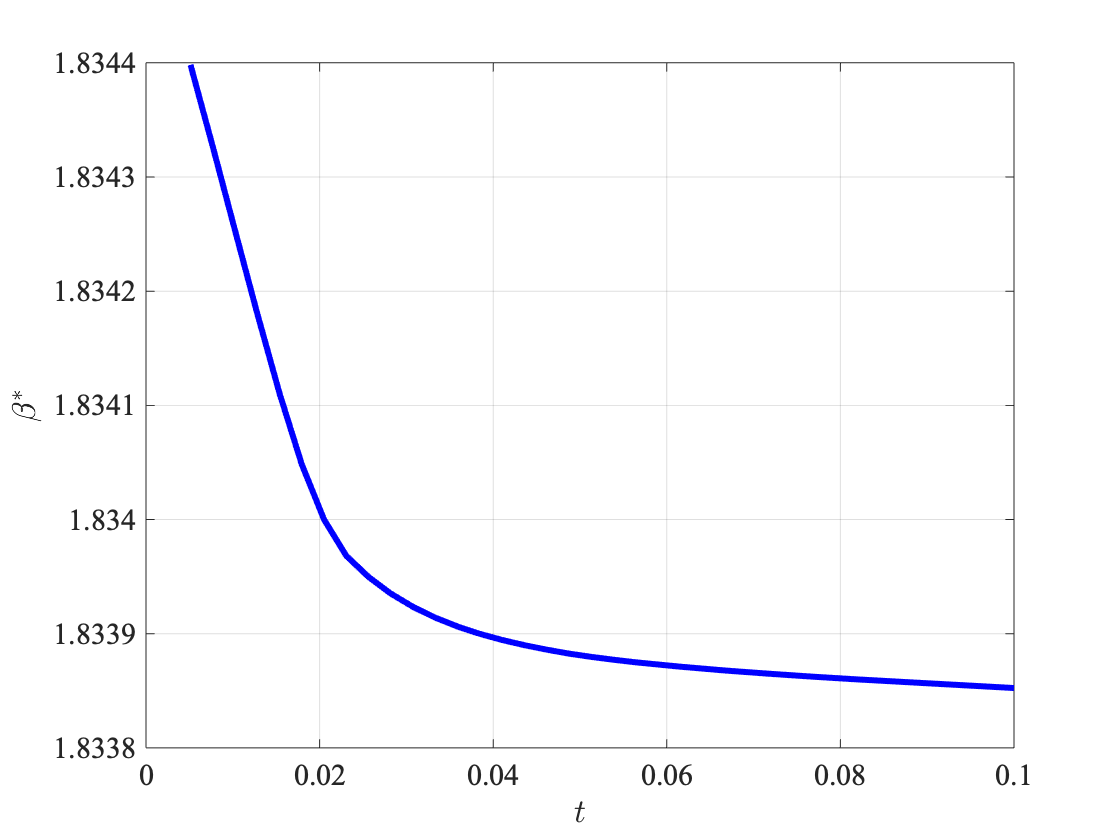}}
    \subfloat{\includegraphics[width=0.25\linewidth]{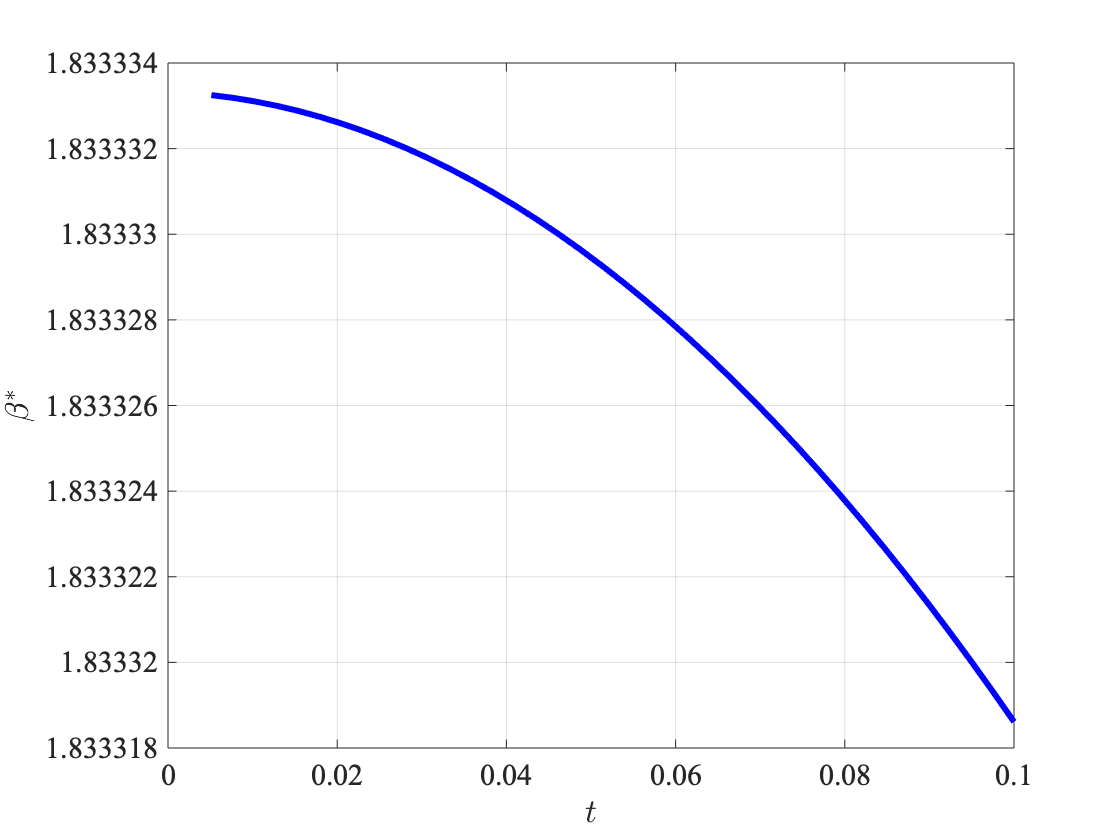}}
    \hspace{0.1in}
    \subfloat{\includegraphics[width=0.25\linewidth]{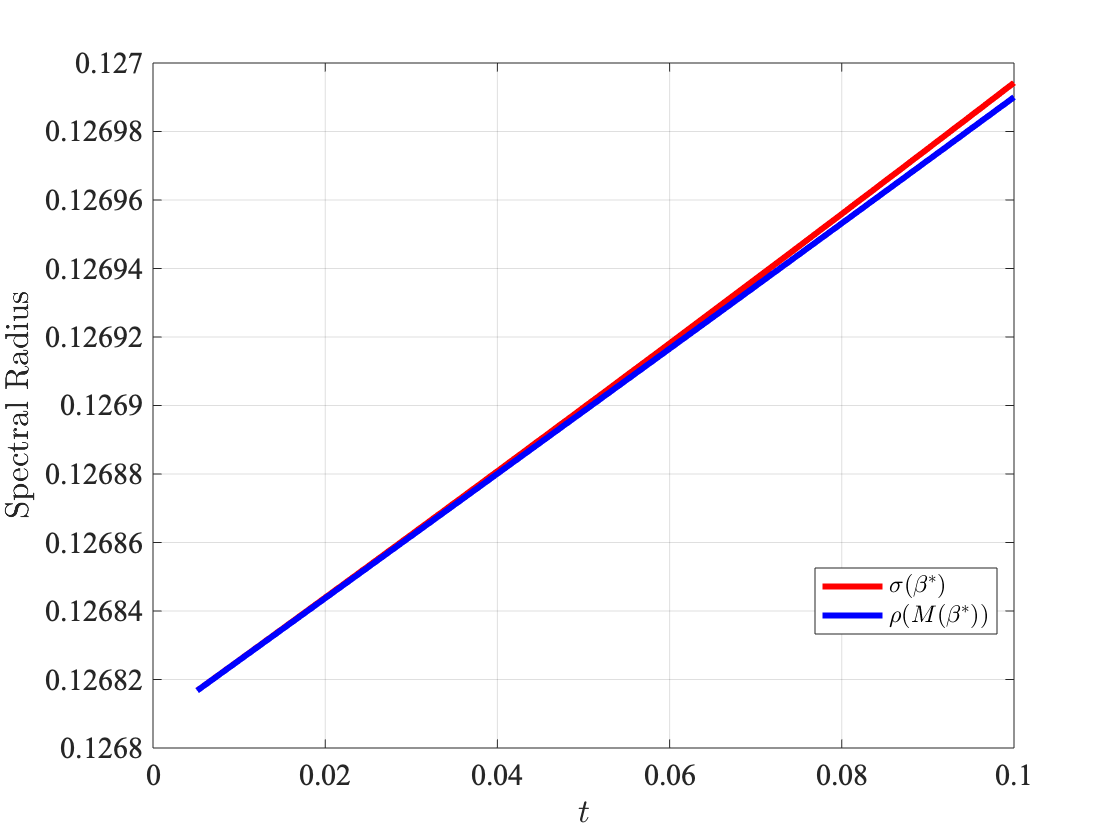}}
    \subfloat{\includegraphics[width=0.25\linewidth]{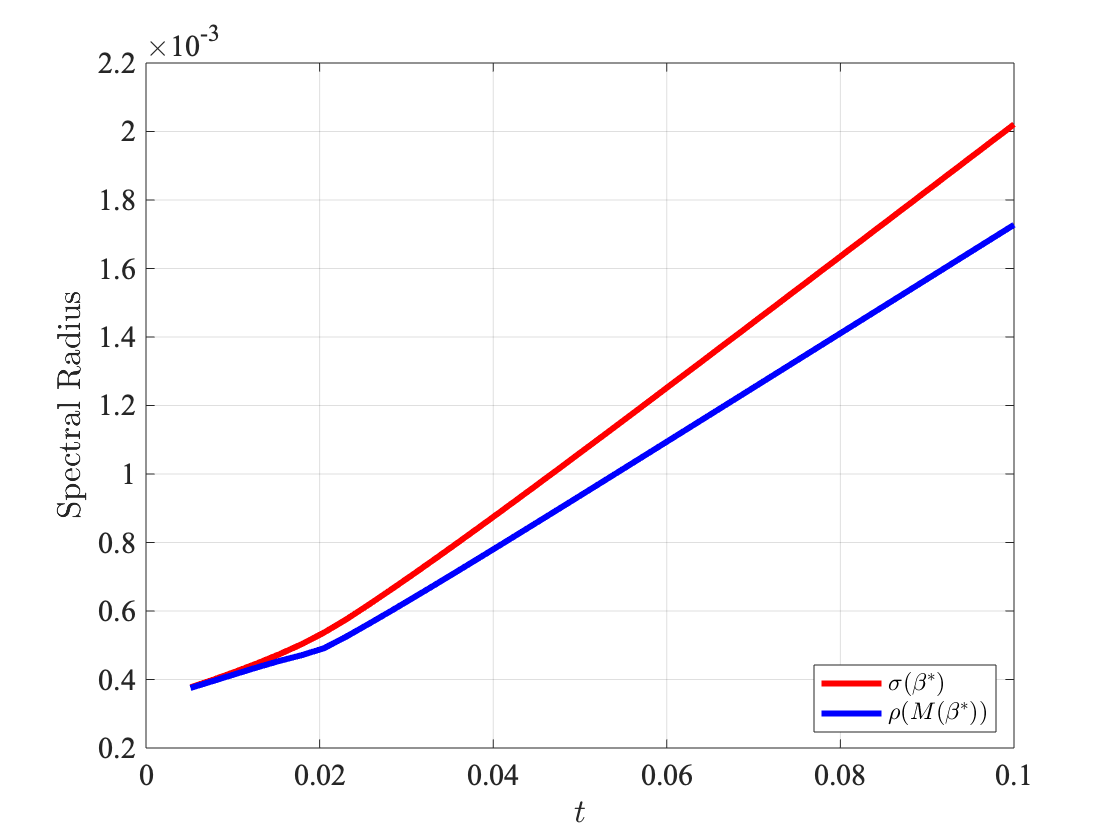}}
    \subfloat{\includegraphics[width=0.25\linewidth]{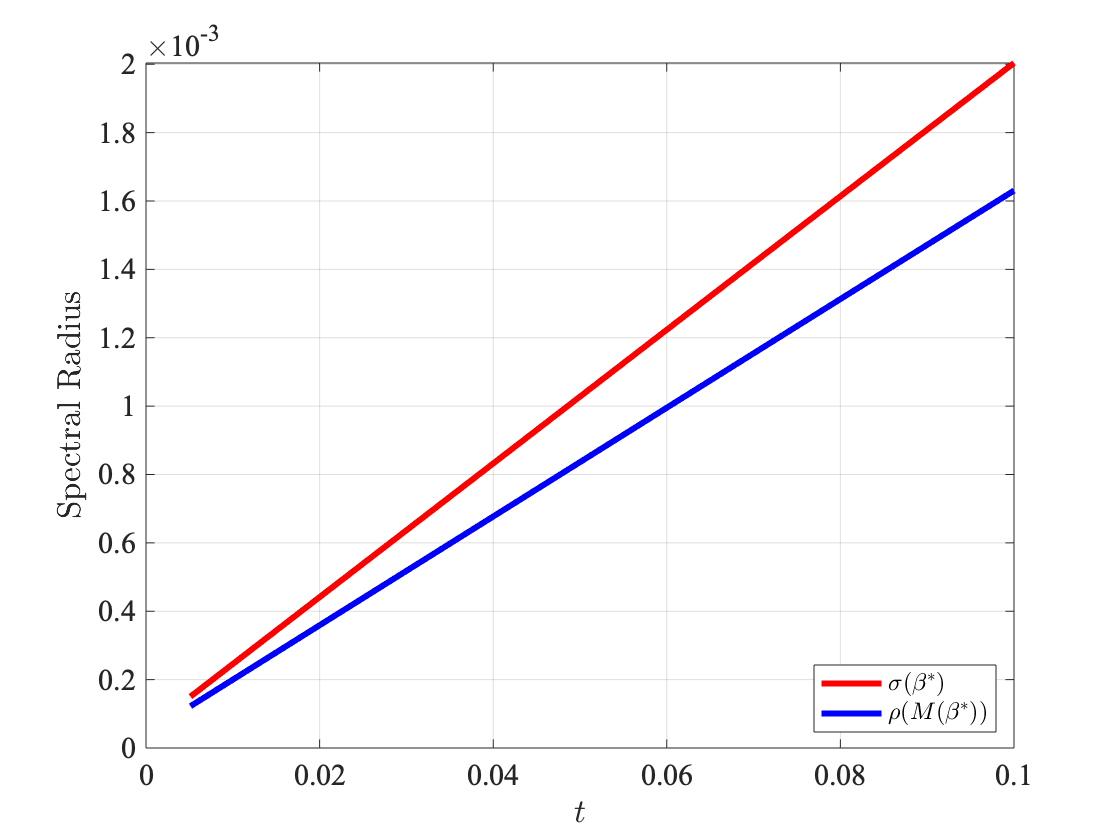}}
\caption{The estimated $\beta^{*}$, the spectral radius of $\rho(M(\beta^{*}))$ and its upper bound $\sigma(\beta^*)$ over time for scheme B in 3D. Top two rows: $N_x=N_y=N_z=4$, $N_t=40$; Bottom two rows: $N_x=N_y=N_z=8$, $N_t=40$. From left to right: $\alpha=0.5,0.001,0$.}
    \label{fig:4}
\end{figure}

Figures~\ref{fig:3} and~\ref{fig:4} show that, for both schemes A and B in 3D, the spectral radius $\rho(M(\beta^))$ remains below the corresponding upper bound $\sigma(\beta^)$ for all tested grid resolutions and damping parameters. Moreover, the spectral radius stays well below one throughout the simulations, numerically supporting the theoretical spectral-radius estimate in the three-dimensional setting.

\begin{figure}[htbp]
	\centering
	\subfloat{\label{N_v1}\includegraphics[width=2in]{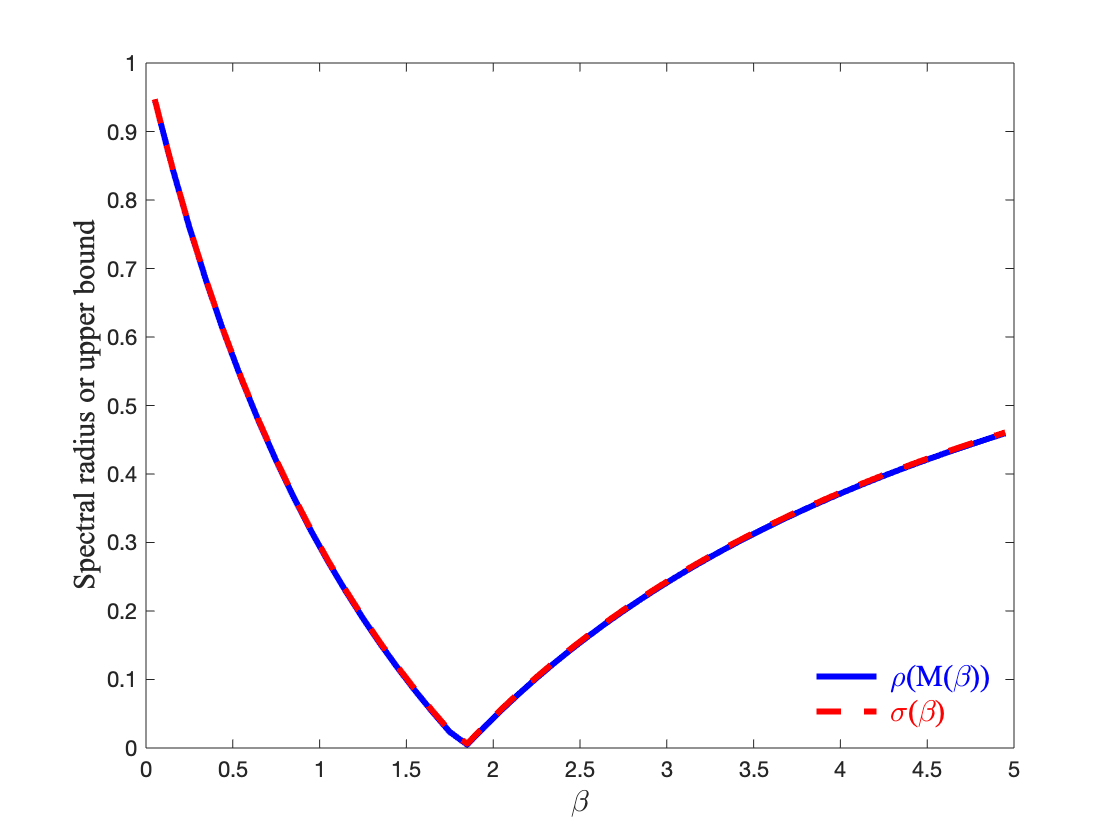}}
	\subfloat{\label{N_v2}\includegraphics[width=2in]{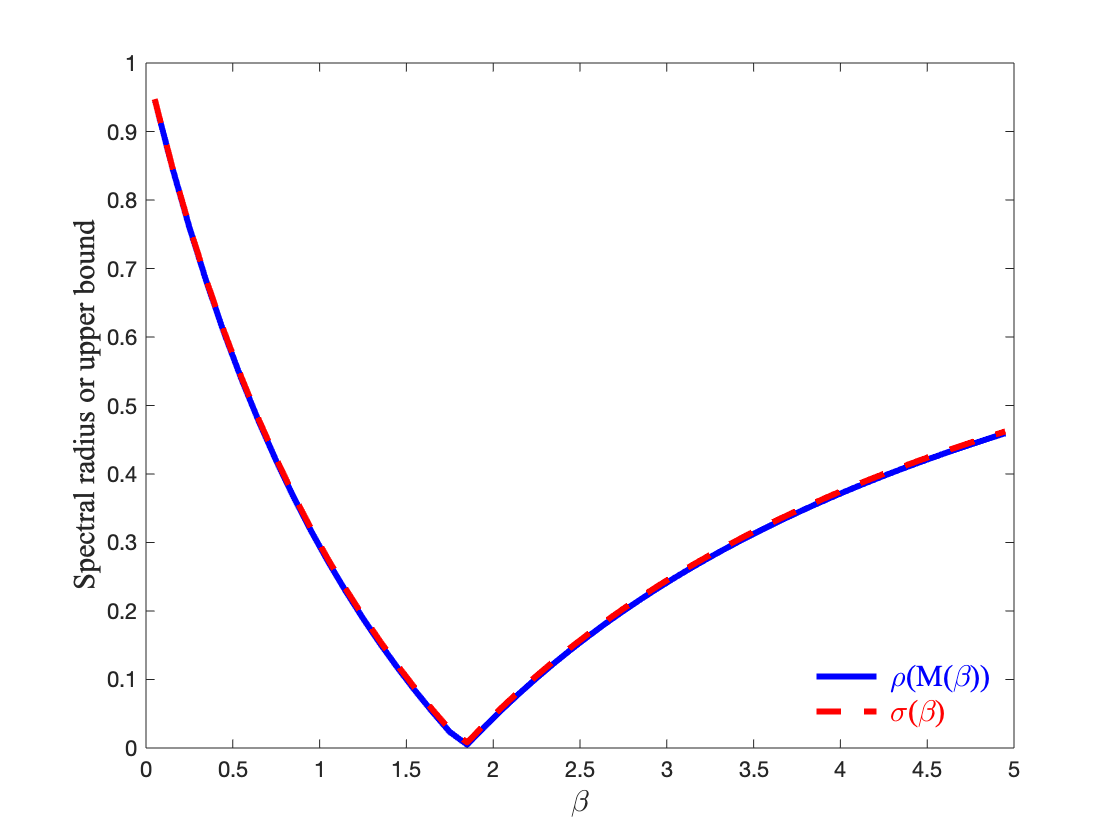}}
	\subfloat{\label{N_v3}\includegraphics[width=2in]{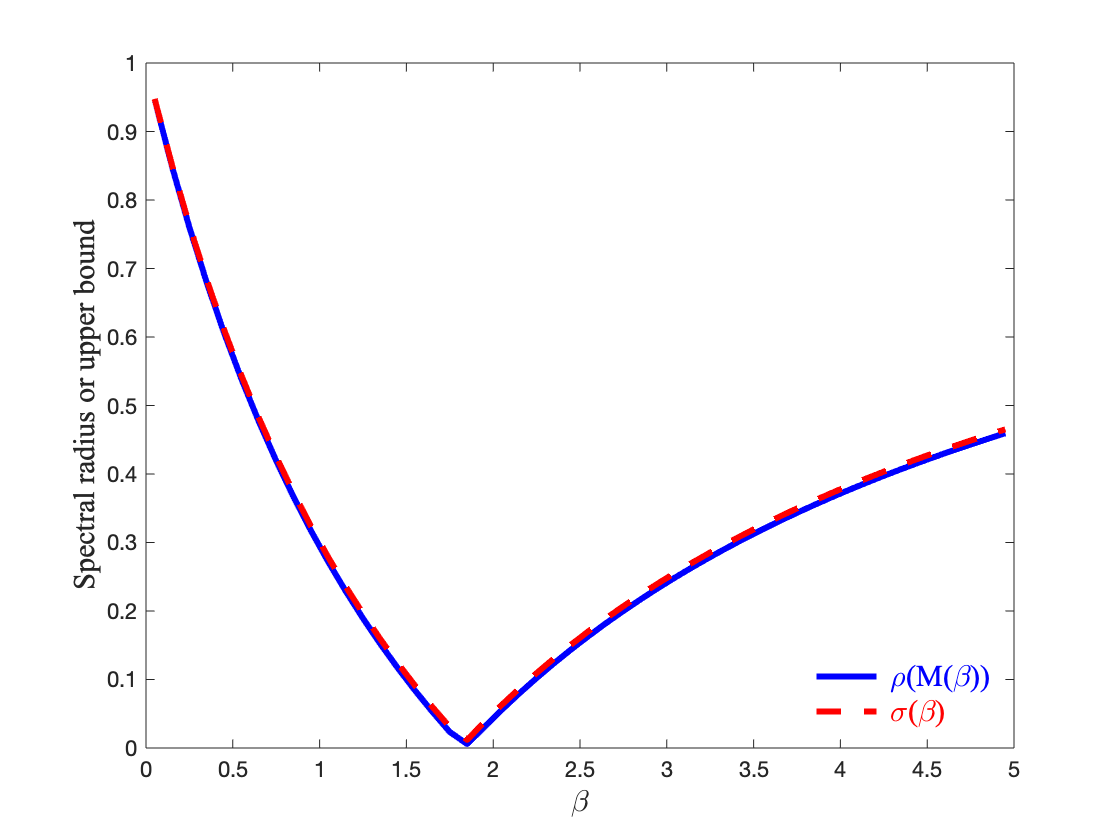}}
    \hspace{0.1in}
    \subfloat{\label{N_v1}\includegraphics[width=2in]{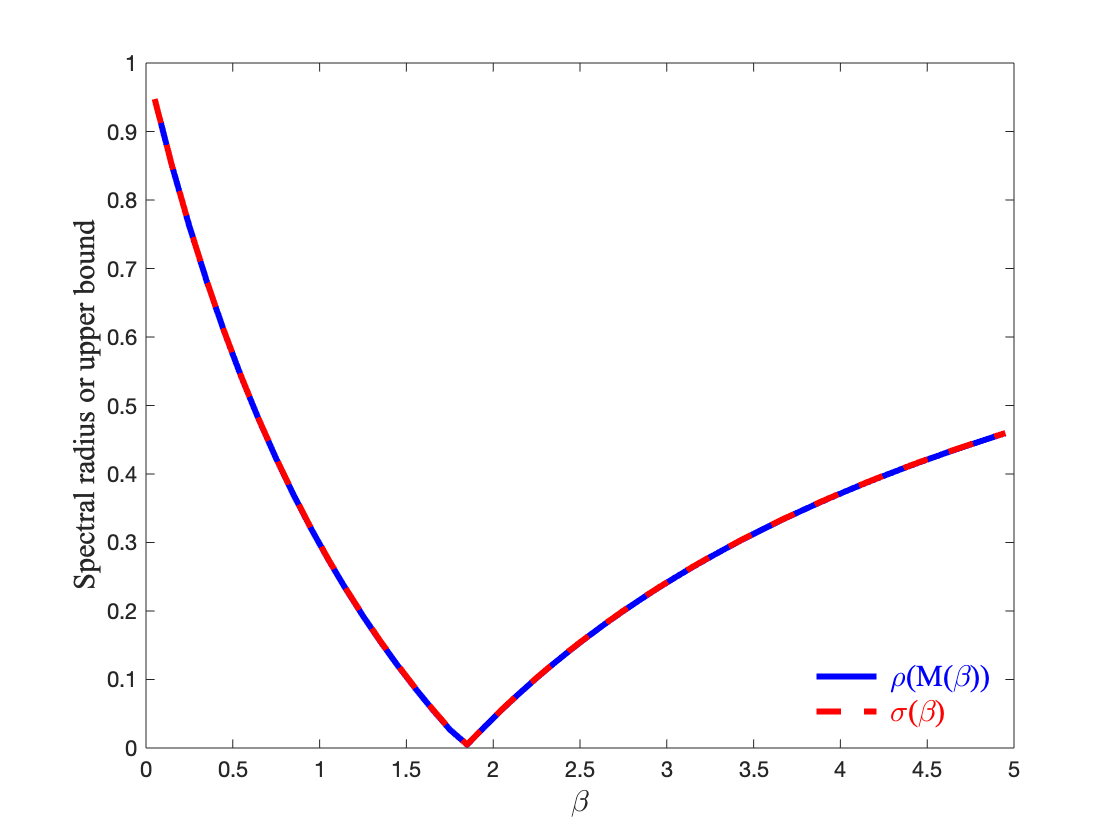}}
	\subfloat{\label{N_v2}\includegraphics[width=2in]{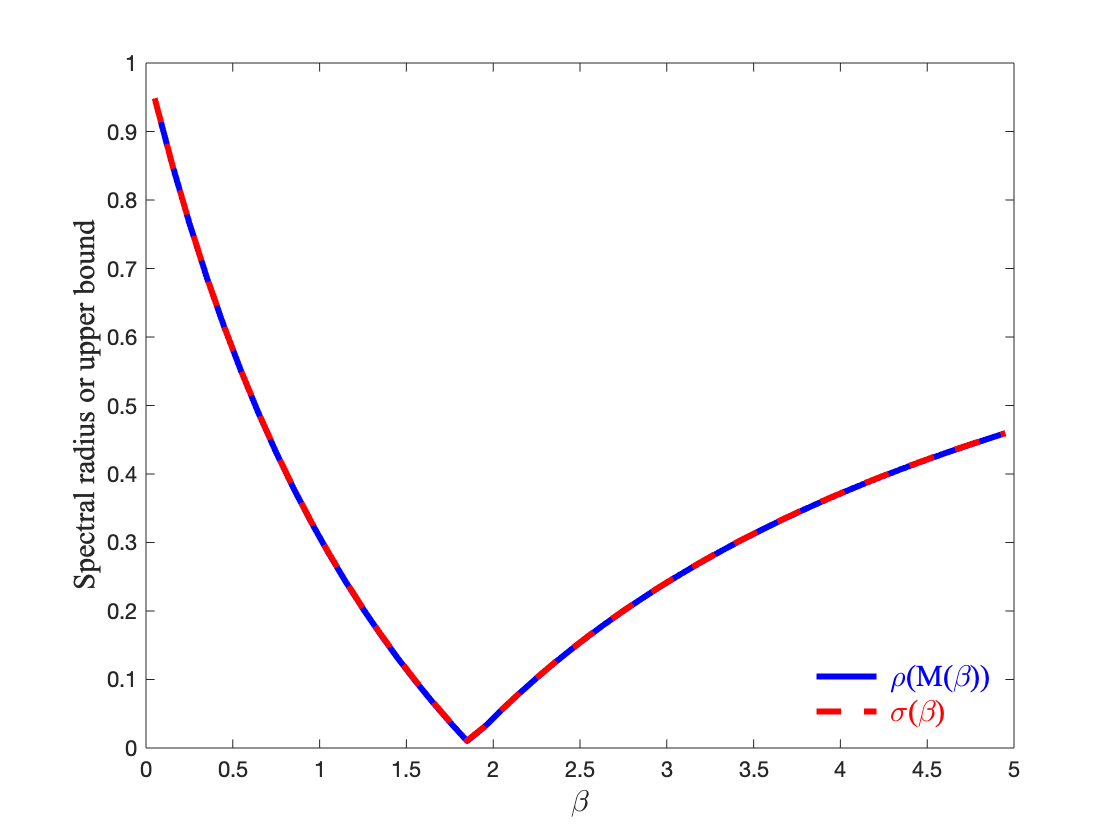}}
	\subfloat{\label{N_v3}\includegraphics[width=2in]{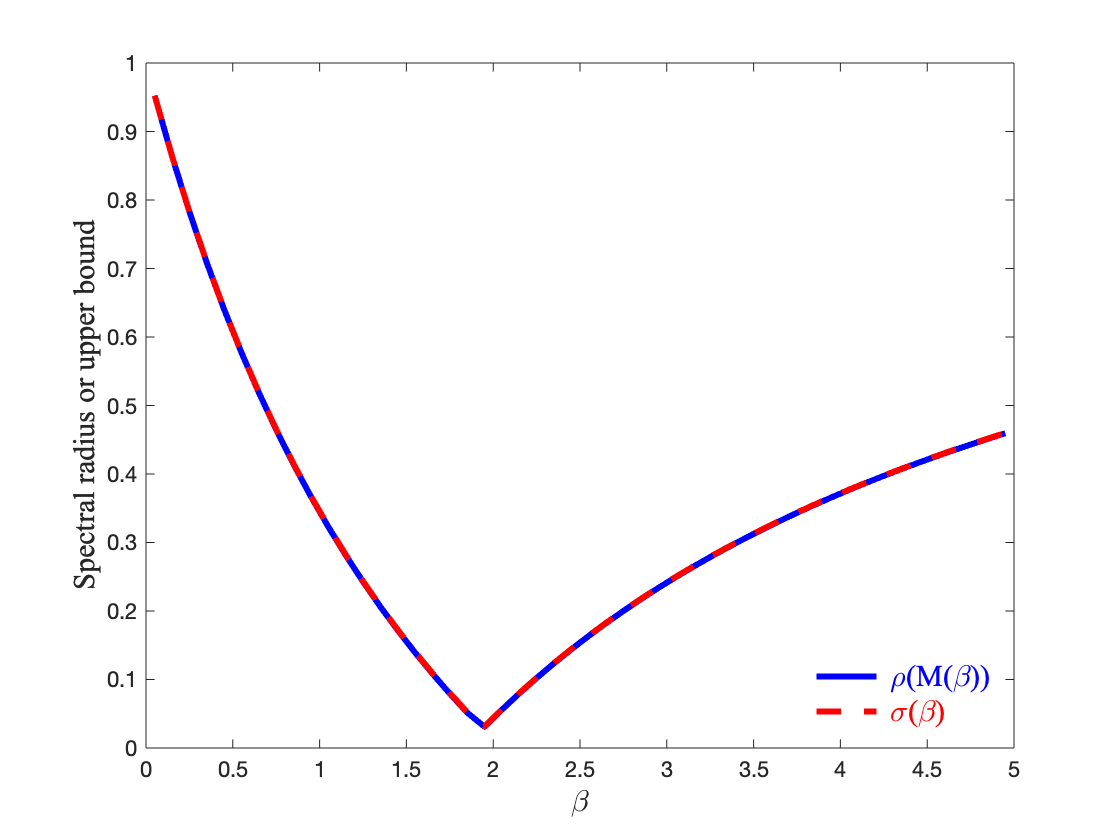}}
    \hspace{0.1in}
    \subfloat{\label{N_v1}\includegraphics[width=2in]{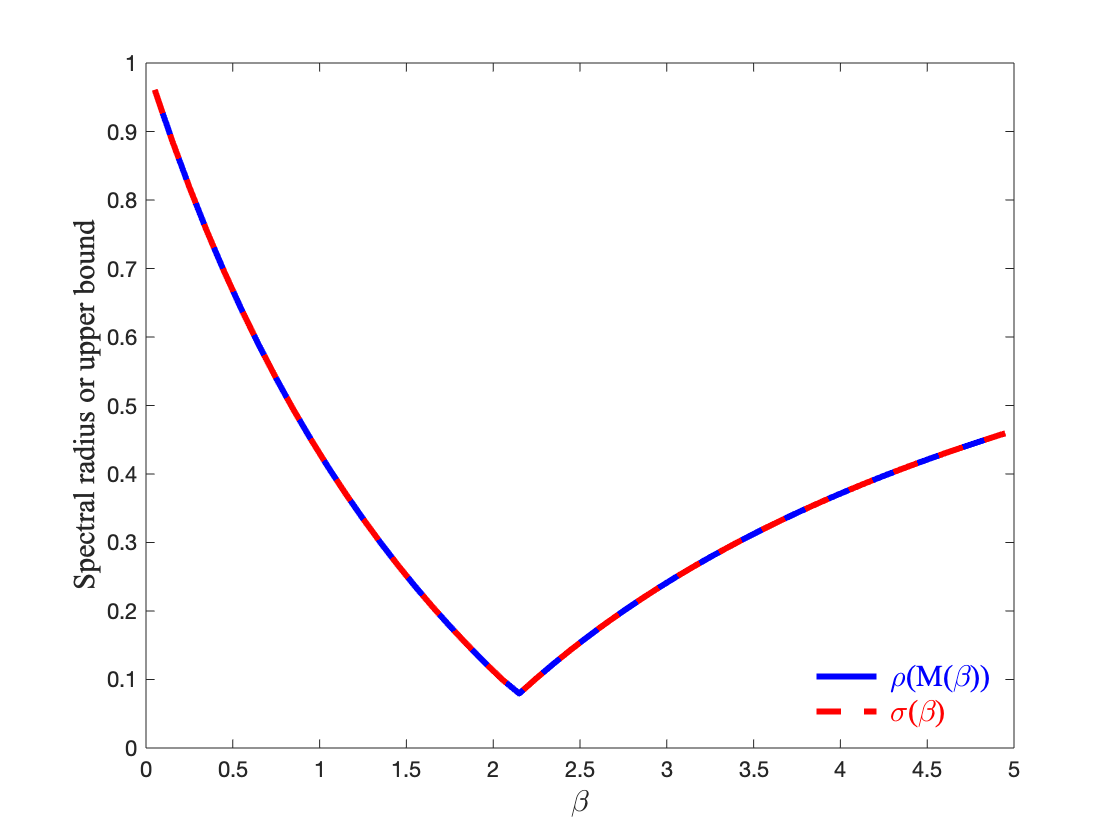}}
	\subfloat{\label{N_v2}\includegraphics[width=2in]{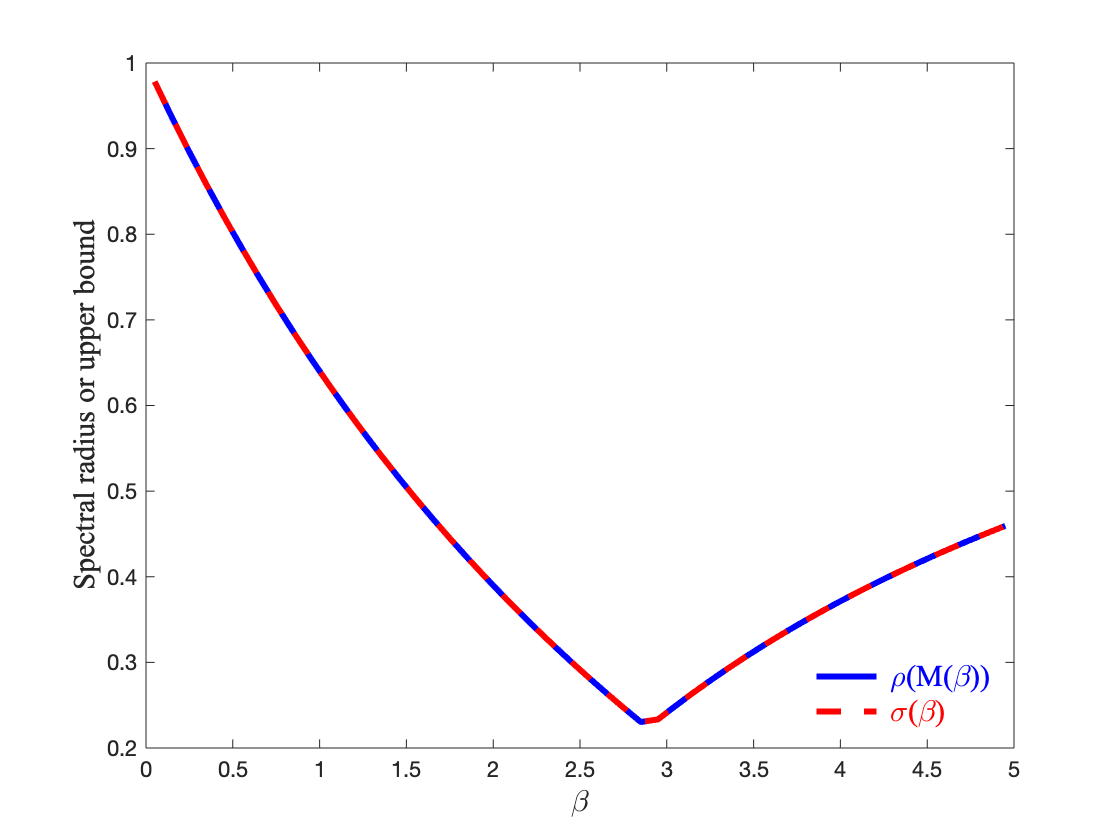}}
	\subfloat{\label{N_v3}\includegraphics[width=2in]{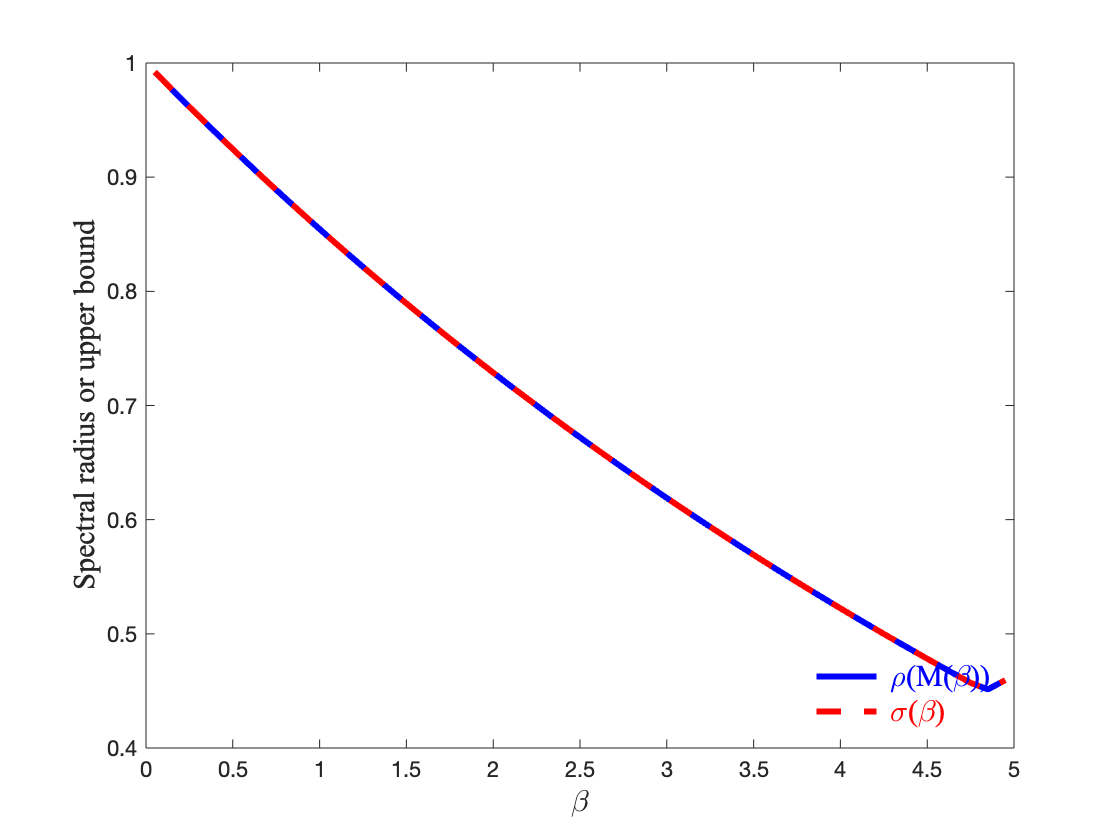}}
\caption{The spectral radius $\rho(M(\beta))$ of the iteration matrices for different $\beta$: ``$-$'' and the upper bound $\sigma(\beta)$ for different $\beta$: ``$--$'' with $T=0.1$, for scheme A in 1D. From left to right: $N_x=16,32,64$. From top to bottom: $\alpha=0,0.01,0.5$.}\label{fig:5}
\end{figure}

\begin{figure}[htbp]
	\centering
	\subfloat{\label{N_v1}\includegraphics[width=2in]{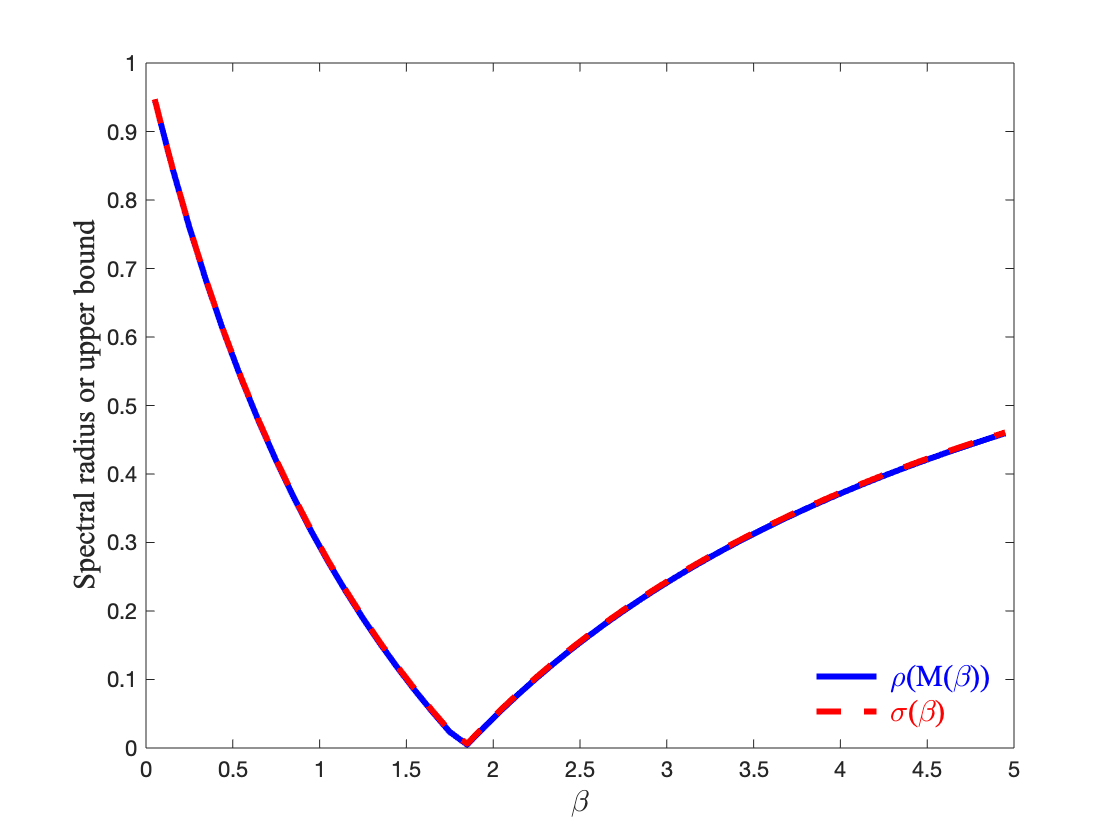}}
	\subfloat{\label{N_v2}\includegraphics[width=2in]{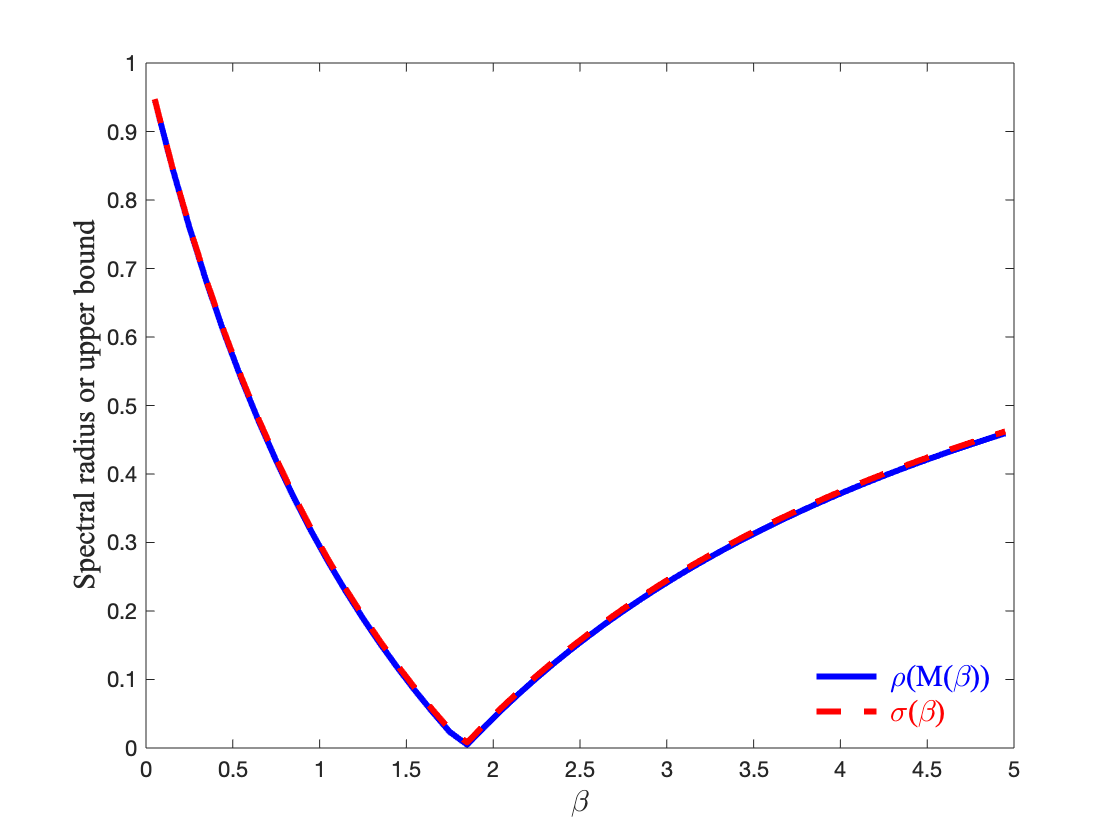}}
	\subfloat{\label{N_v3}\includegraphics[width=2in]{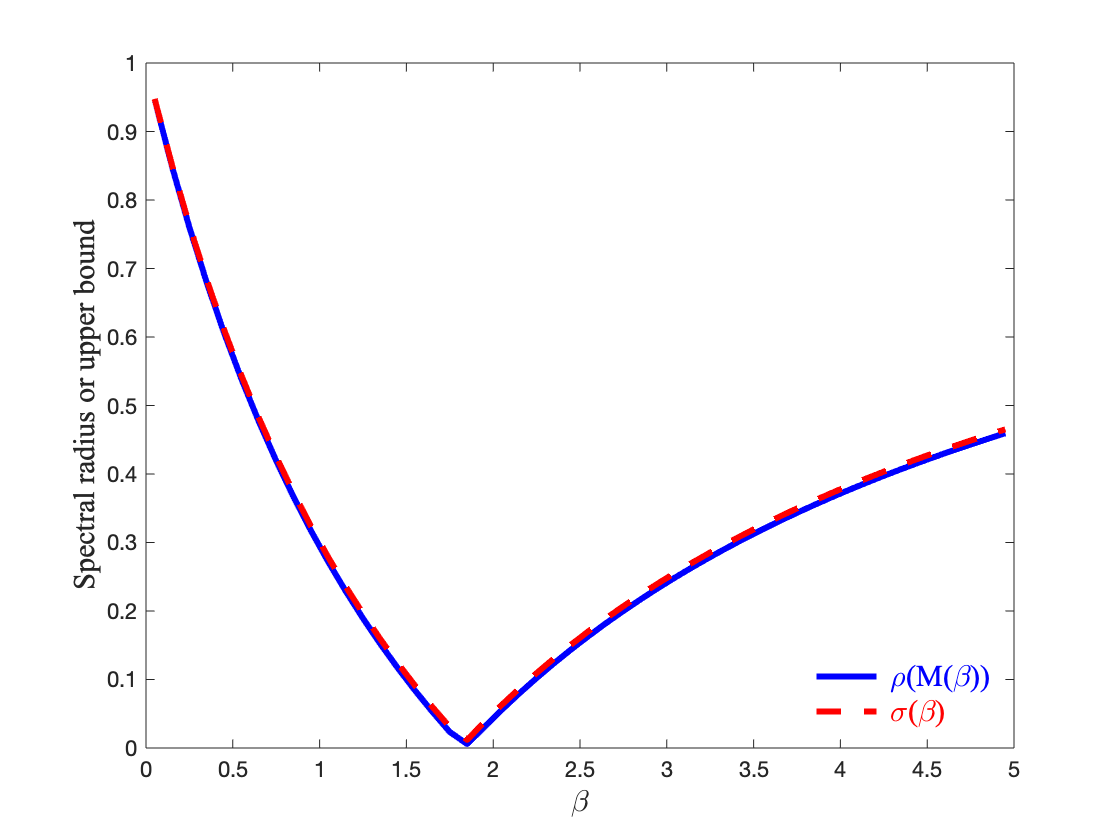}}
    \hspace{0.1in}
    \subfloat{\label{N_v1}\includegraphics[width=2in]{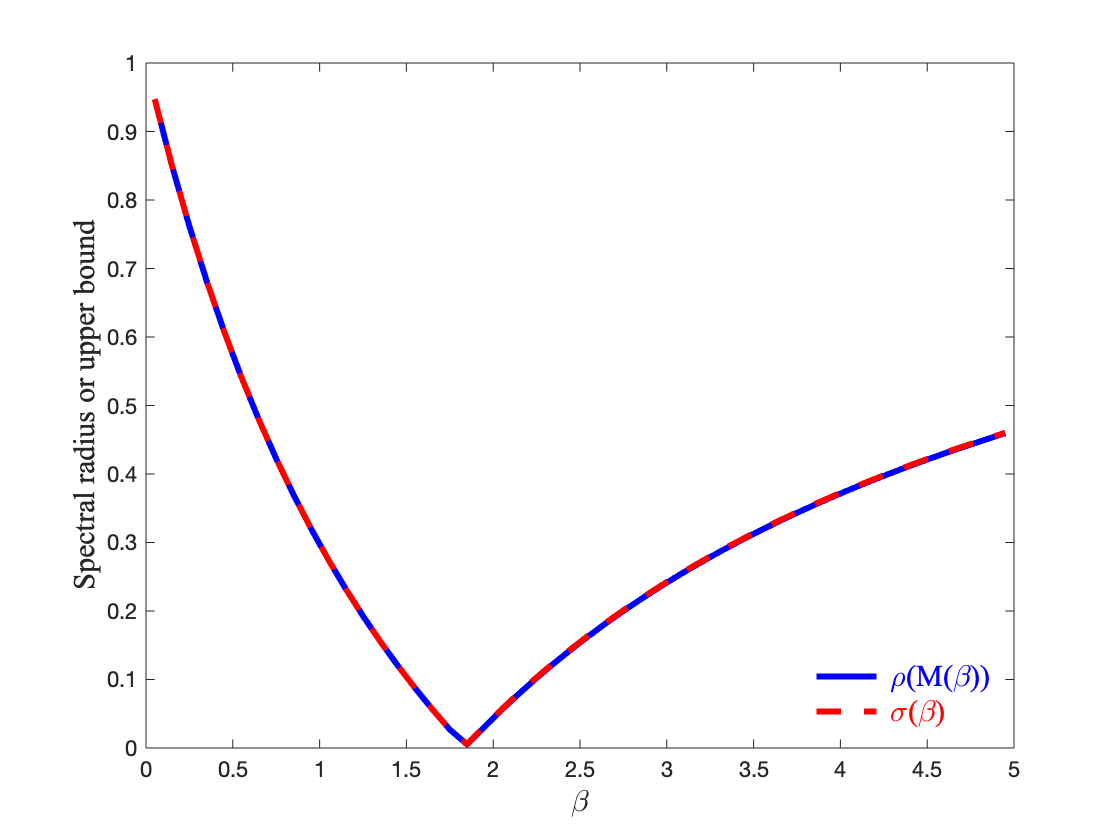}}
	\subfloat{\label{N_v2}\includegraphics[width=2in]{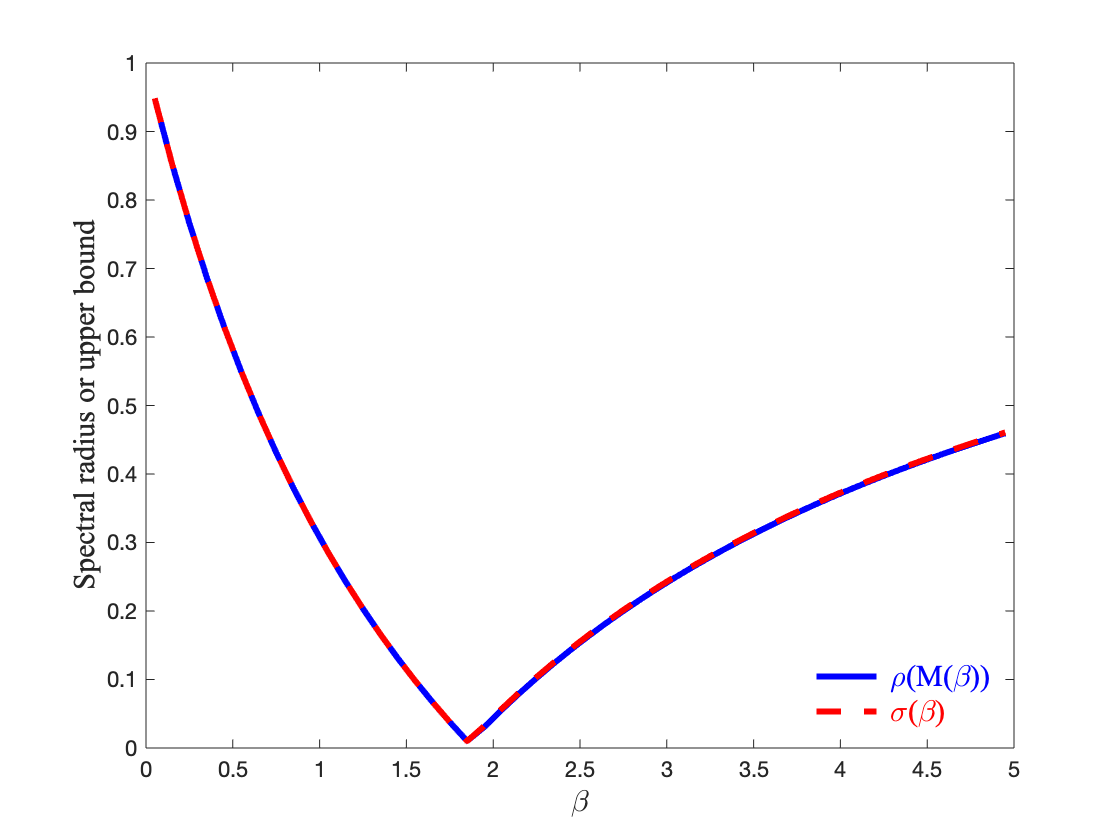}}
	\subfloat{\label{N_v3}\includegraphics[width=2in]{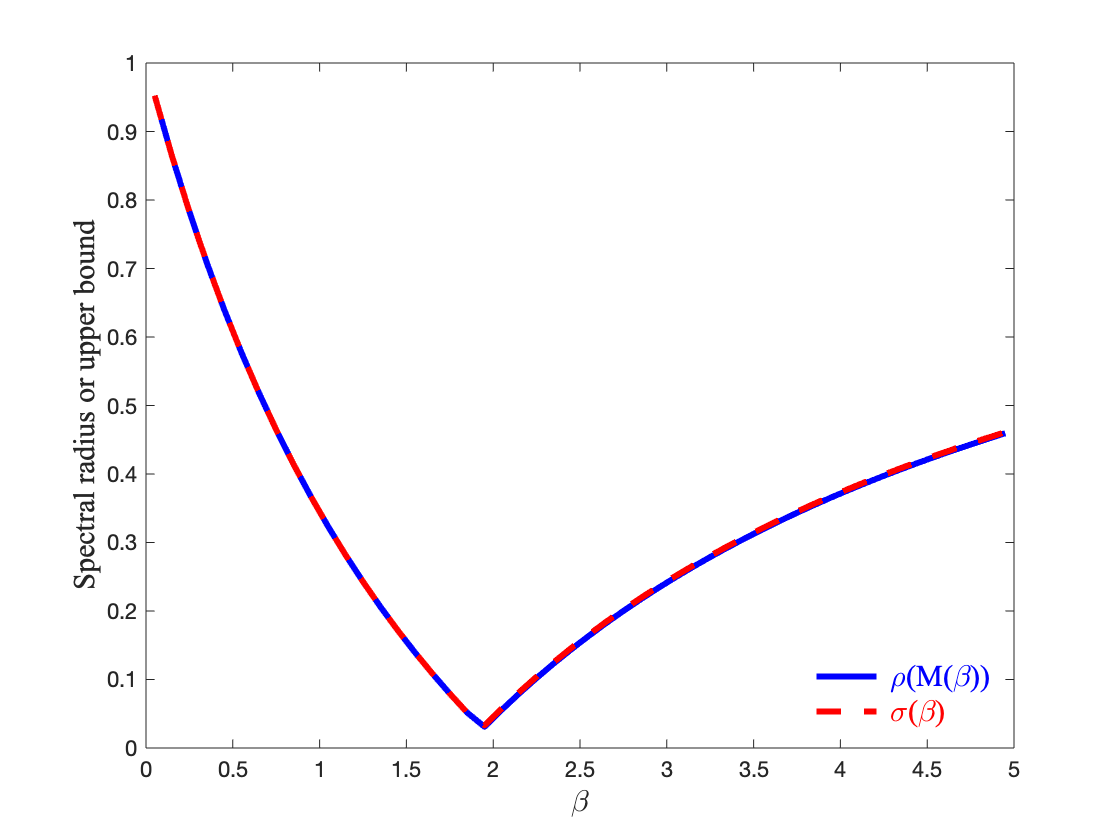}}
    \hspace{0.1in}
    \subfloat{\label{N_v1}\includegraphics[width=2in]{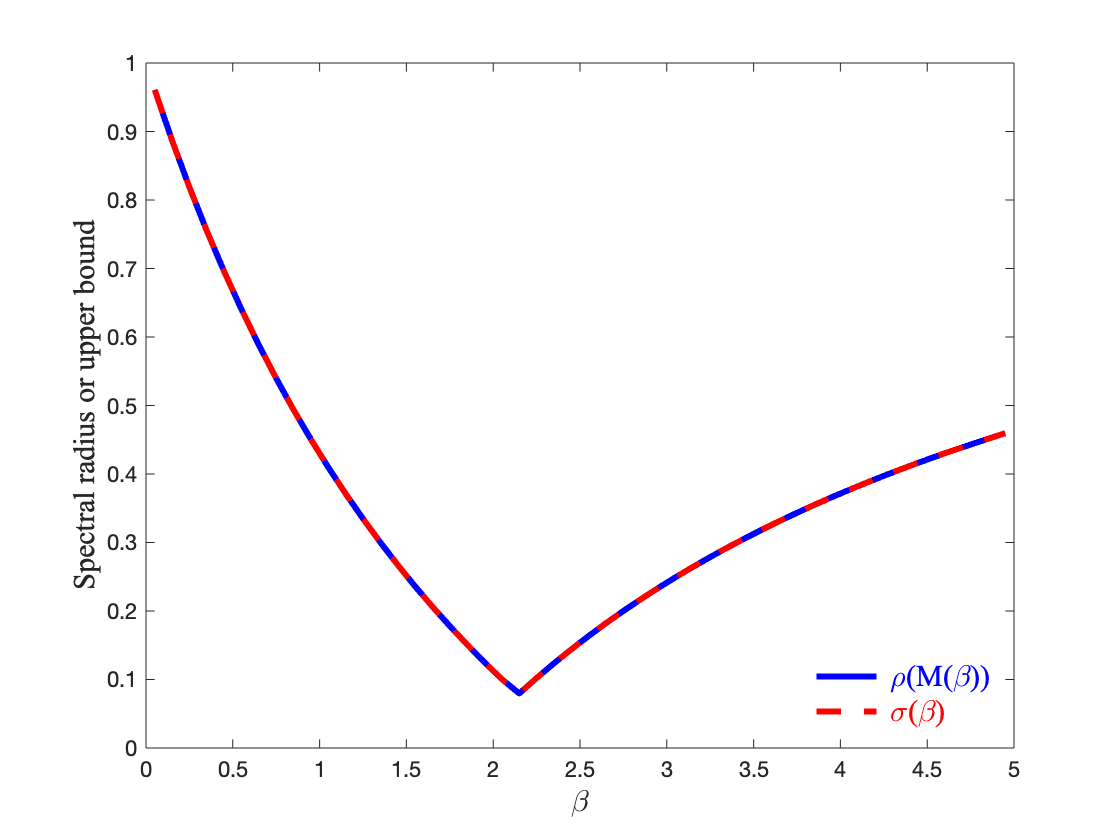}}
	\subfloat{\label{N_v2}\includegraphics[width=2in]{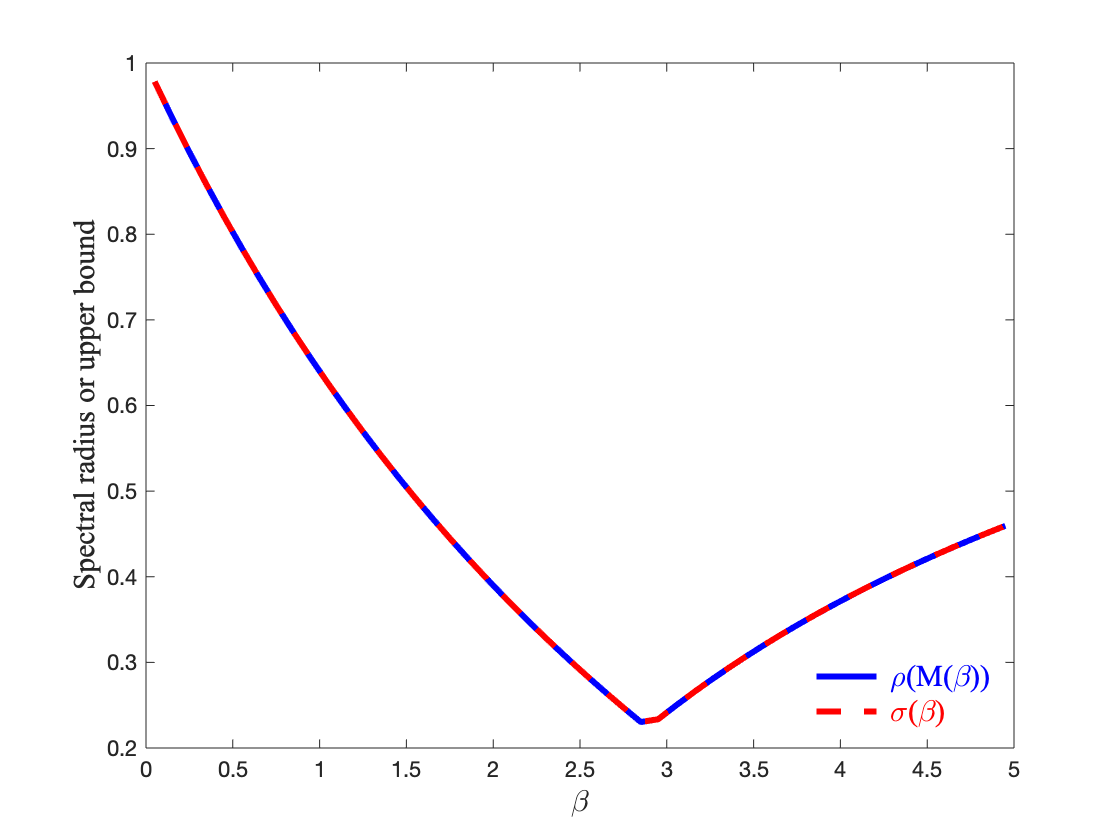}}
	\subfloat{\label{N_v3}\includegraphics[width=2in]{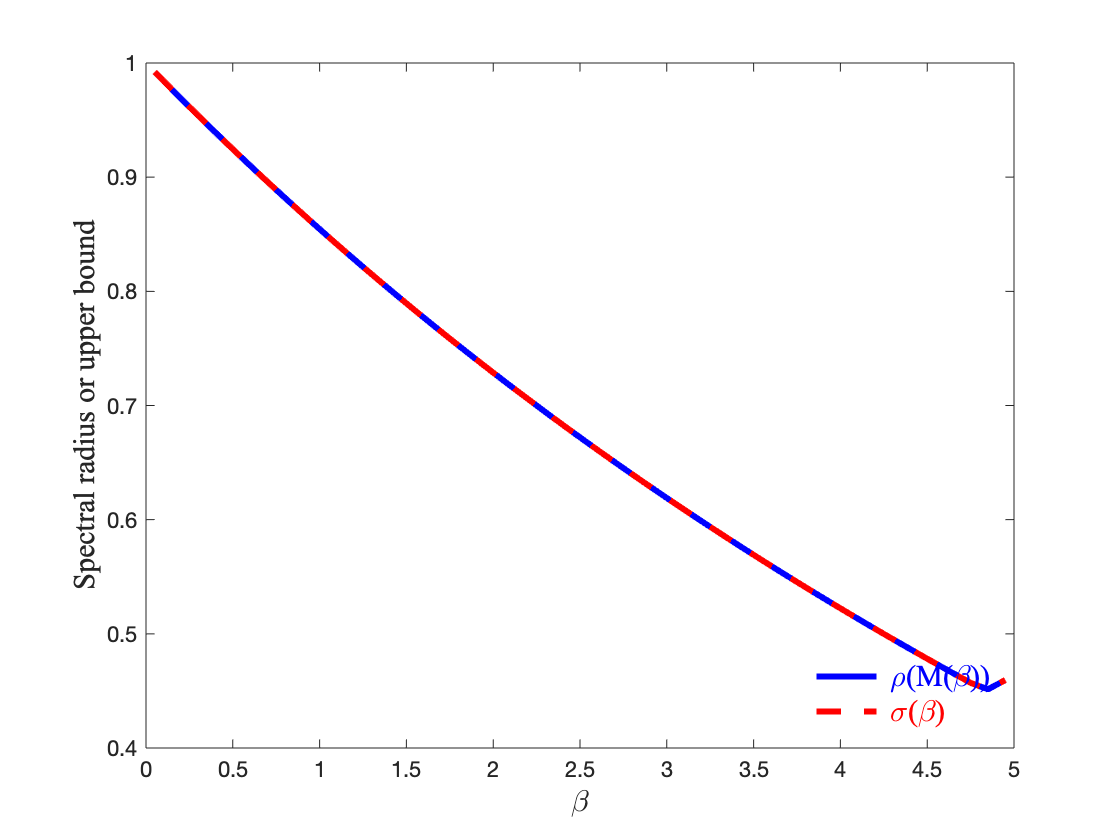}}
\caption{The spectral radius $\rho(M(\beta))$ of the iteration matrices for different $\beta$: ``$-$'' and the upper bound $\sigma(\beta)$ for different $\beta$: ``$--$'' with $T=0.1$, for scheme B in 1D. From left to right: $N_x=16,32,64$. From top to bottom: $\alpha=0,0.01,0.5$.}\label{fig:6}
\end{figure}

\begin{figure}[htbp]
	\centering
	\subfloat{\label{N_v1}\includegraphics[width=2in]{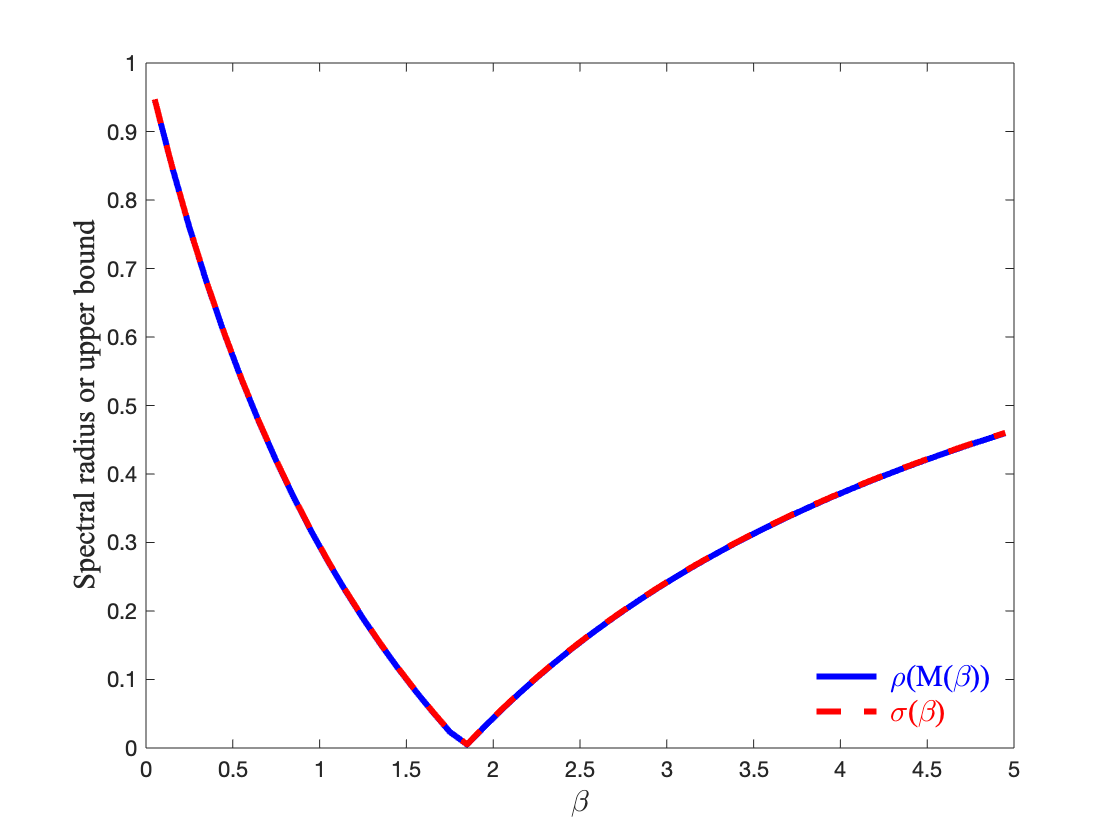}}
	\subfloat{\label{N_v2}\includegraphics[width=2in]{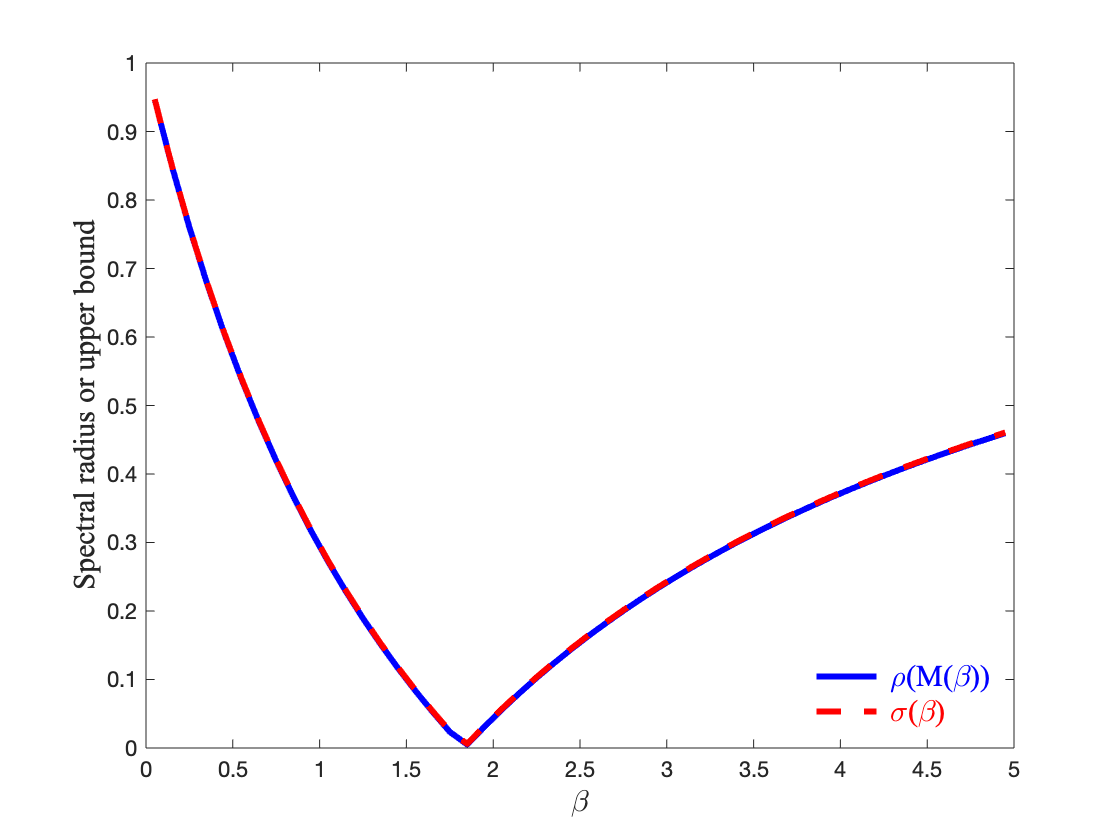}}
	\subfloat{\label{N_v3}\includegraphics[width=2in]{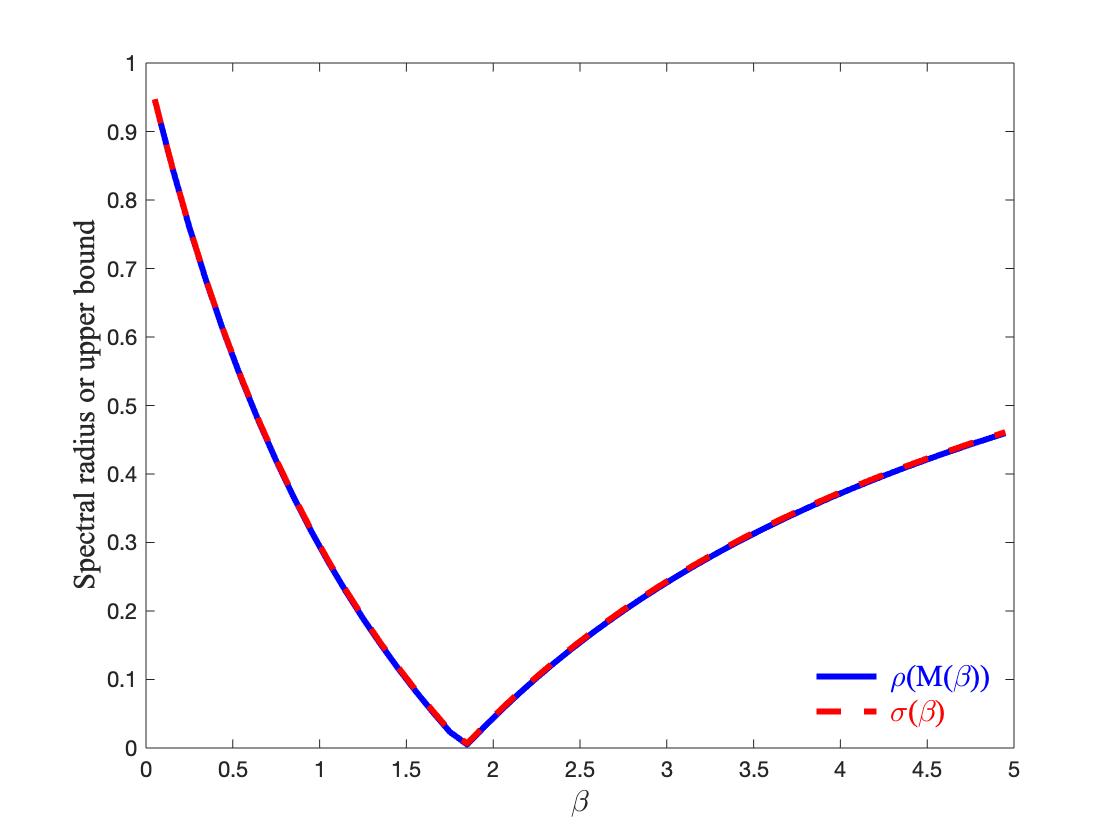}}
    \hspace{0.1in}
    \subfloat{\label{N_v1}\includegraphics[width=2in]{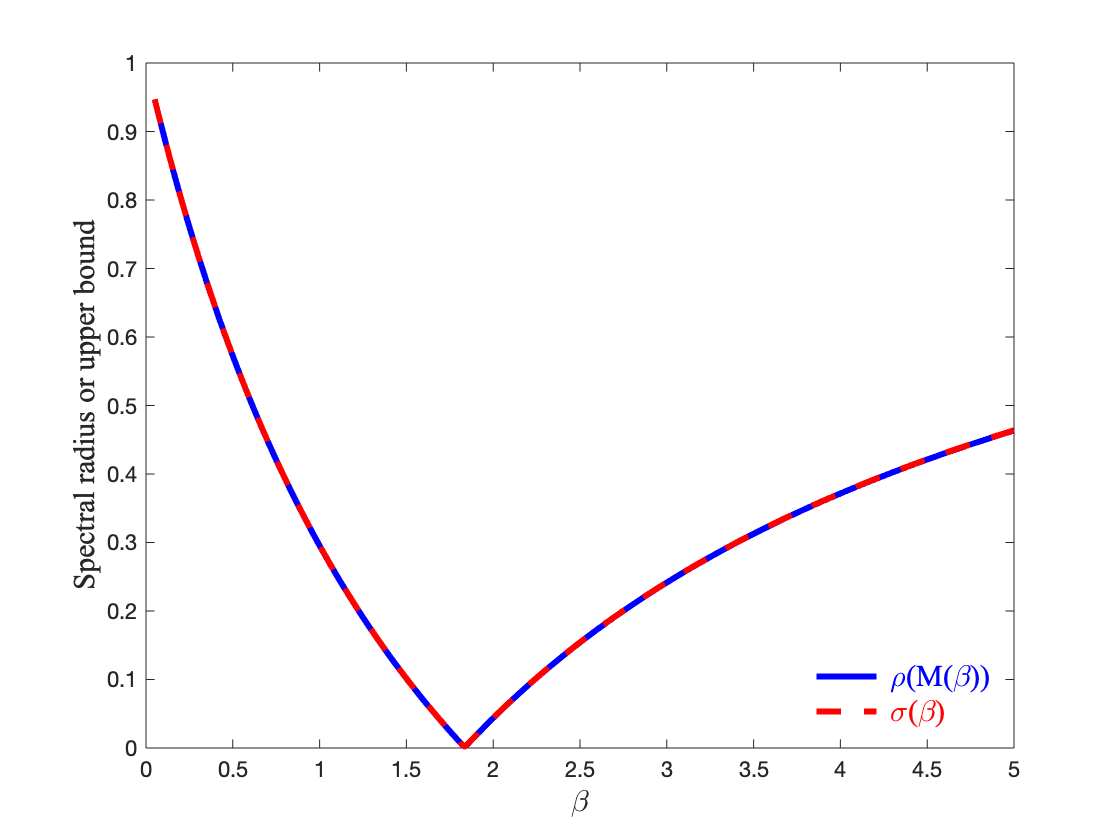}}
	\subfloat{\label{N_v2}\includegraphics[width=2in]{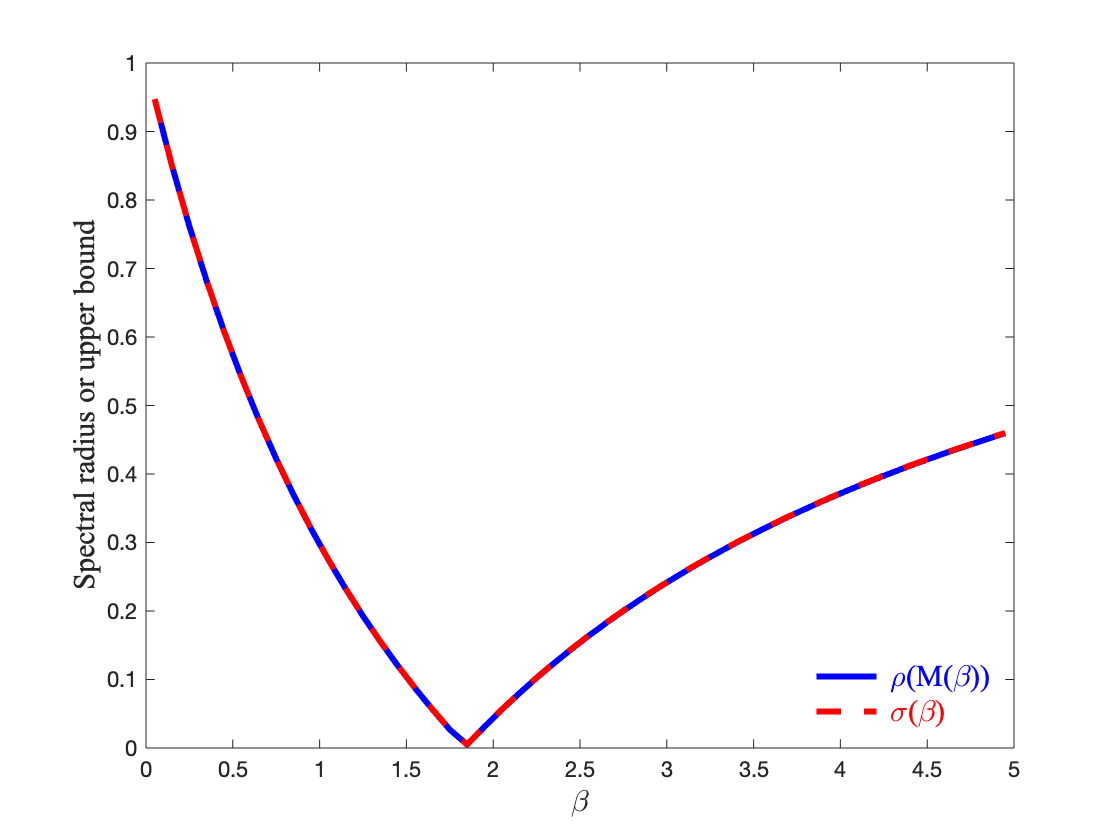}}
	\subfloat{\label{N_v3}\includegraphics[width=2in]{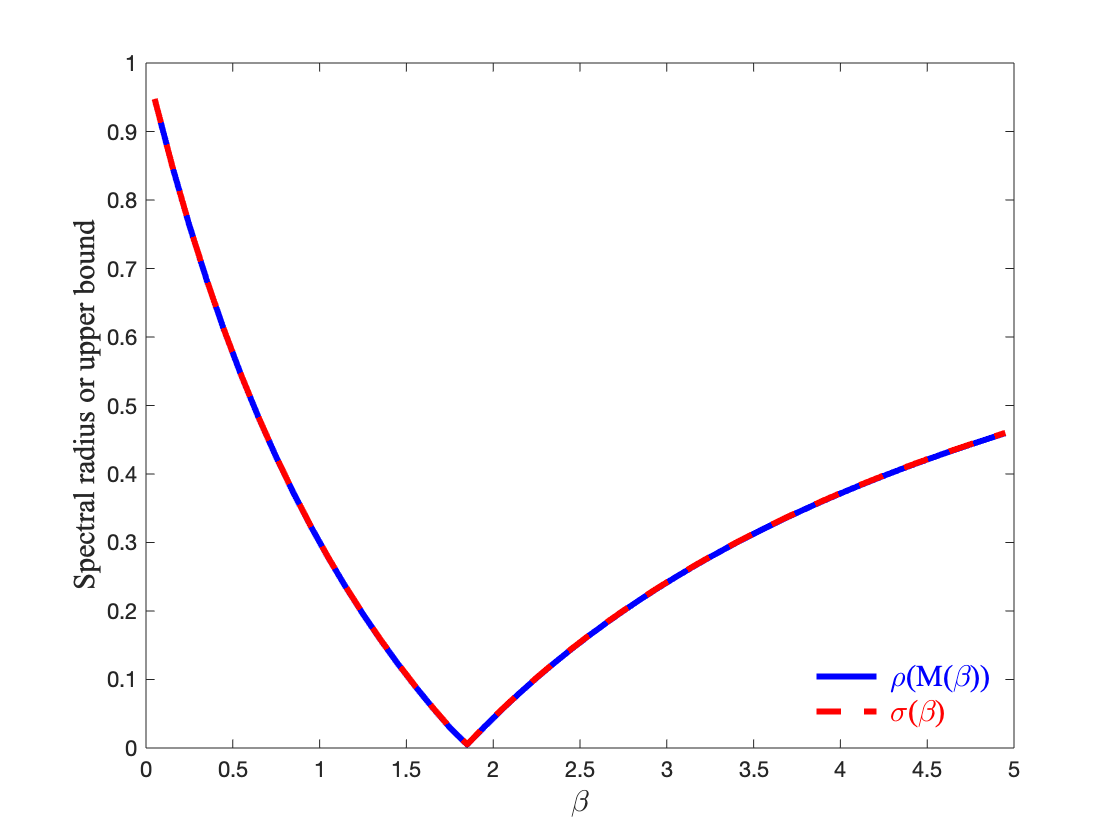}}
    \hspace{0.1in}
    \subfloat{\label{N_v1}\includegraphics[width=2in]{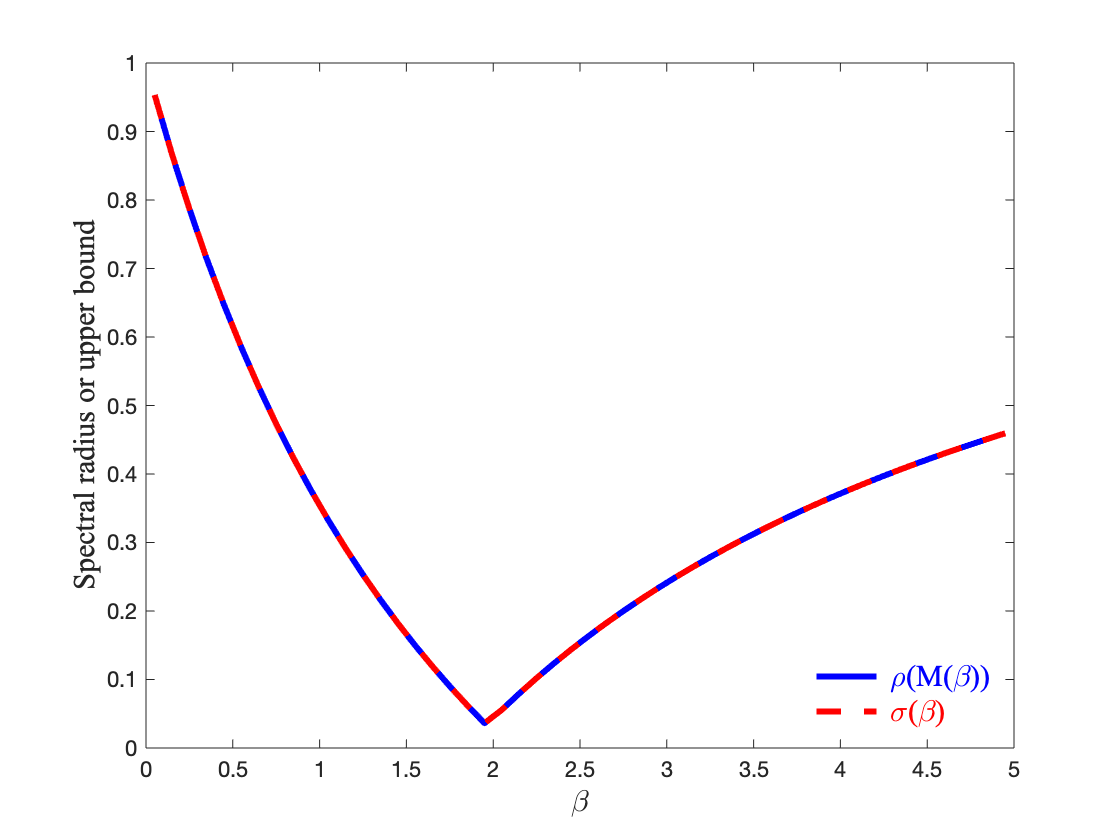}}
	\subfloat{\label{N_v2}\includegraphics[width=2in]{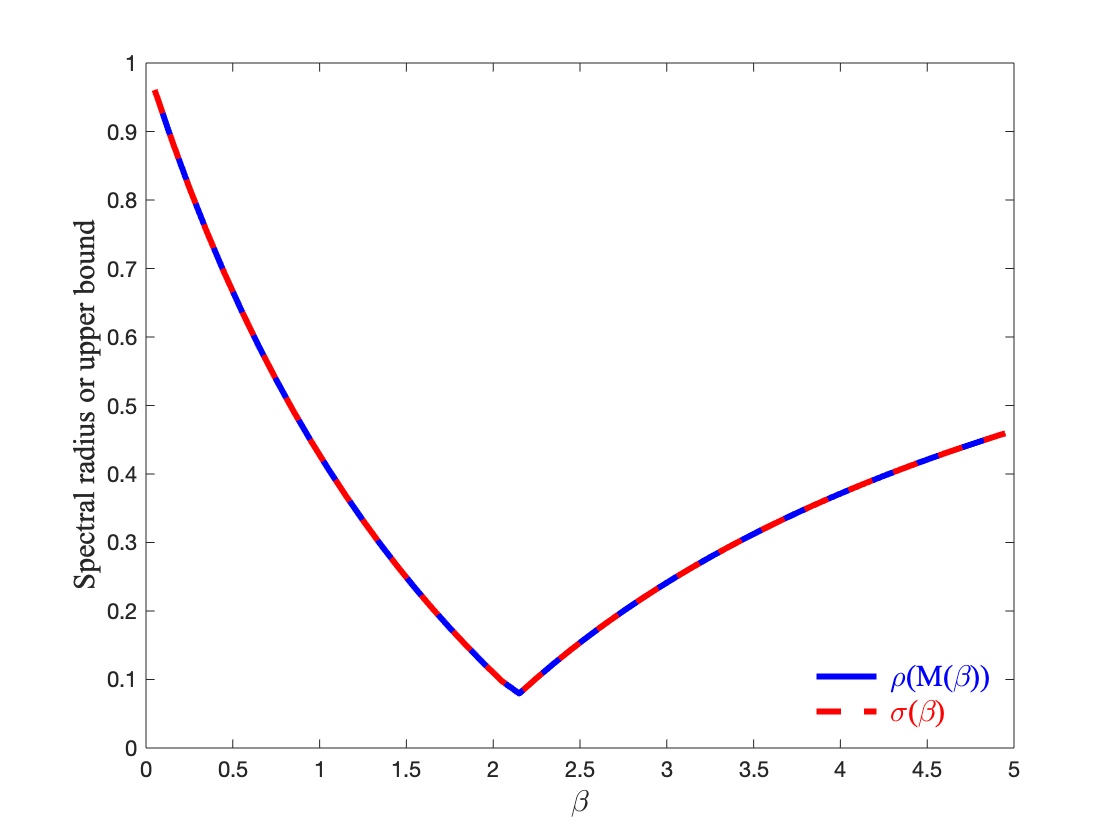}}
	\subfloat{\label{N_v3}\includegraphics[width=2in]{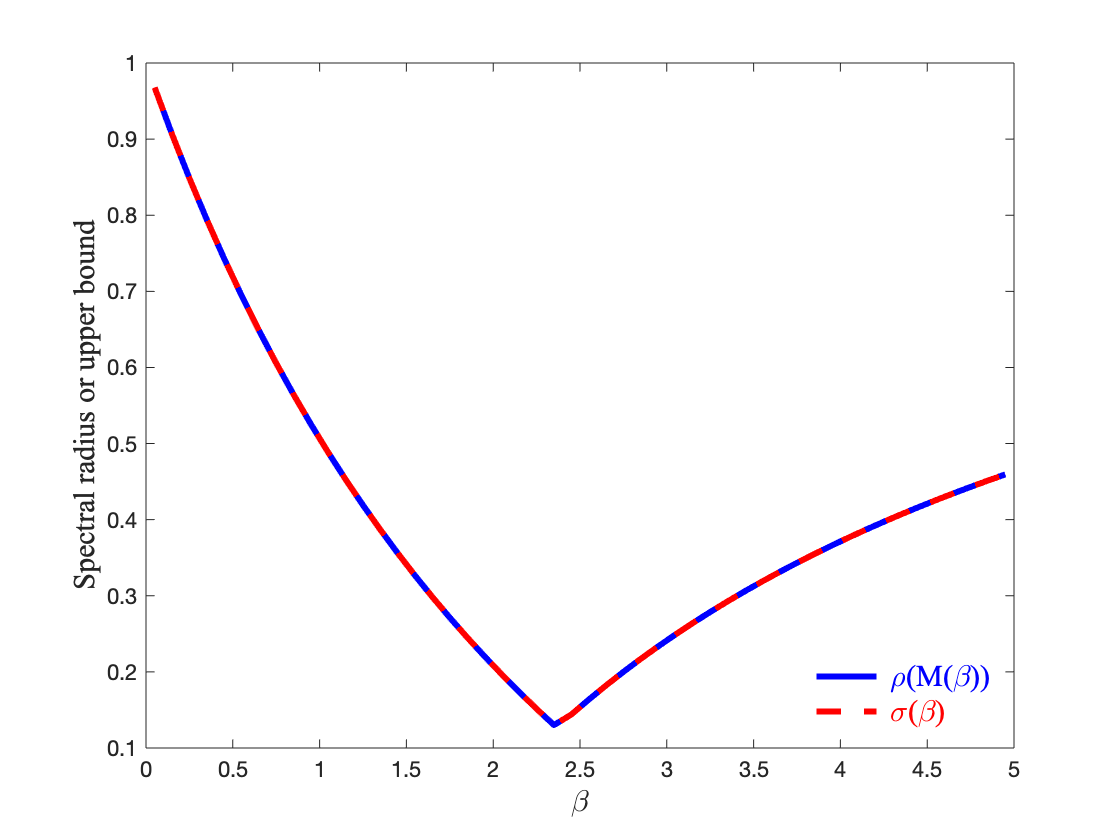}}
\caption{The spectral radius $\rho(M(\beta))$ of the iteration matrices for different $\beta$: ``$-$'' and the upper bound $\sigma(\beta)$ for different $\beta$: ``$--$'' with $T=0.1$ for scheme A in 3D. From left to right: $N_x=N_y=N_z=4,6,8$. From top to bottom: $\alpha=0,0.01,0.5$.}\label{fig:7}
\end{figure}

\begin{figure}[htbp]
	\centering
	\subfloat{\label{N_v1}\includegraphics[width=2in]{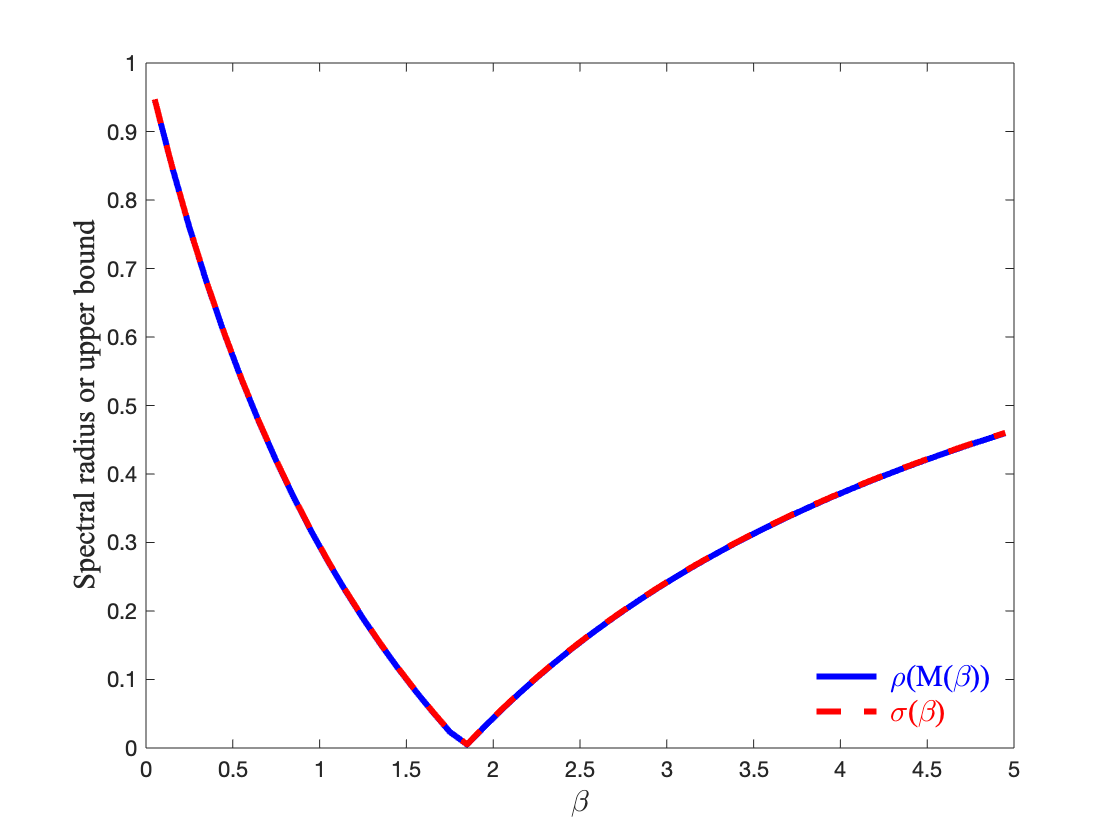}}
	\subfloat{\label{N_v2}\includegraphics[width=2in]{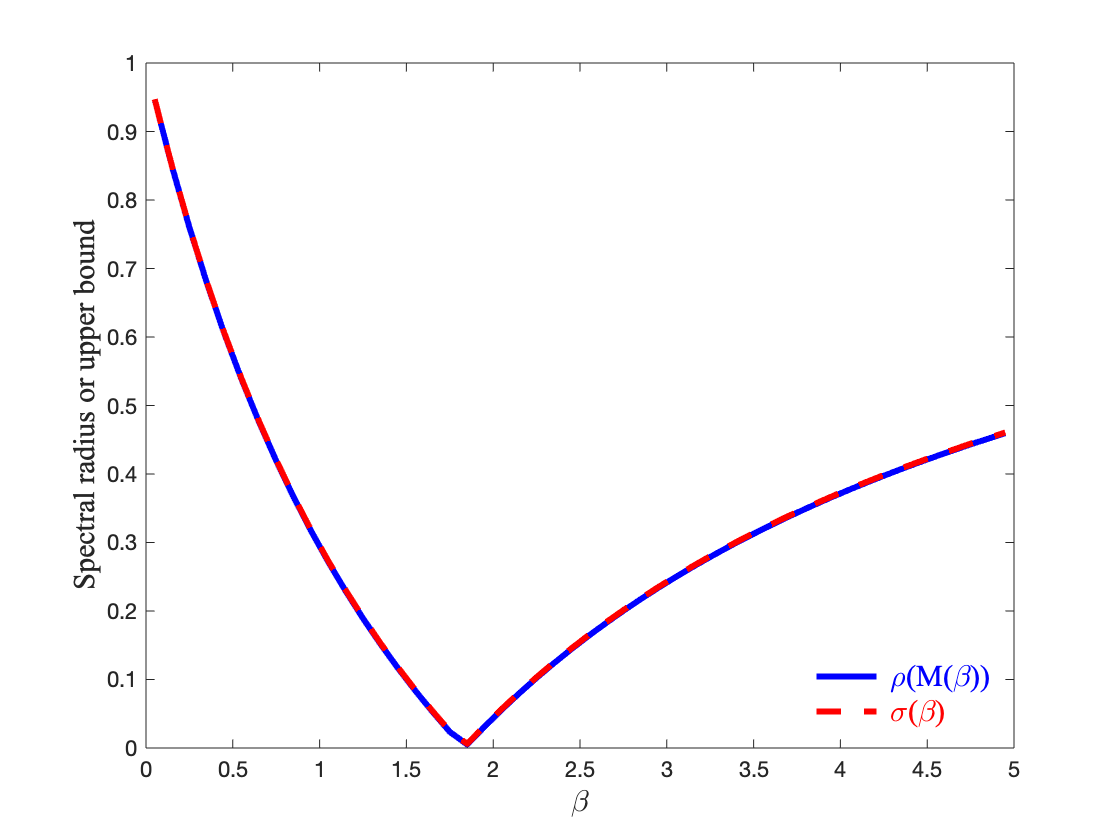}}
	\subfloat{\label{N_v3}\includegraphics[width=2in]{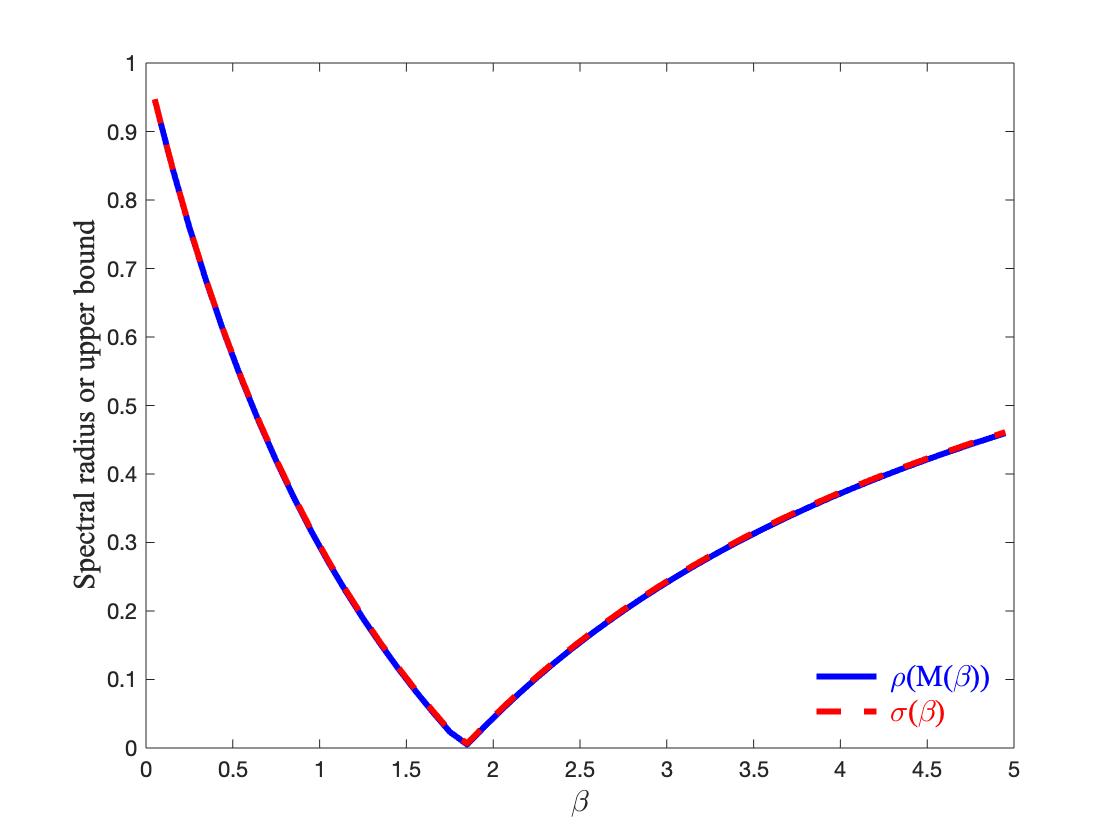}}
    \hspace{0.1in}
    \subfloat{\label{N_v1}\includegraphics[width=2in]{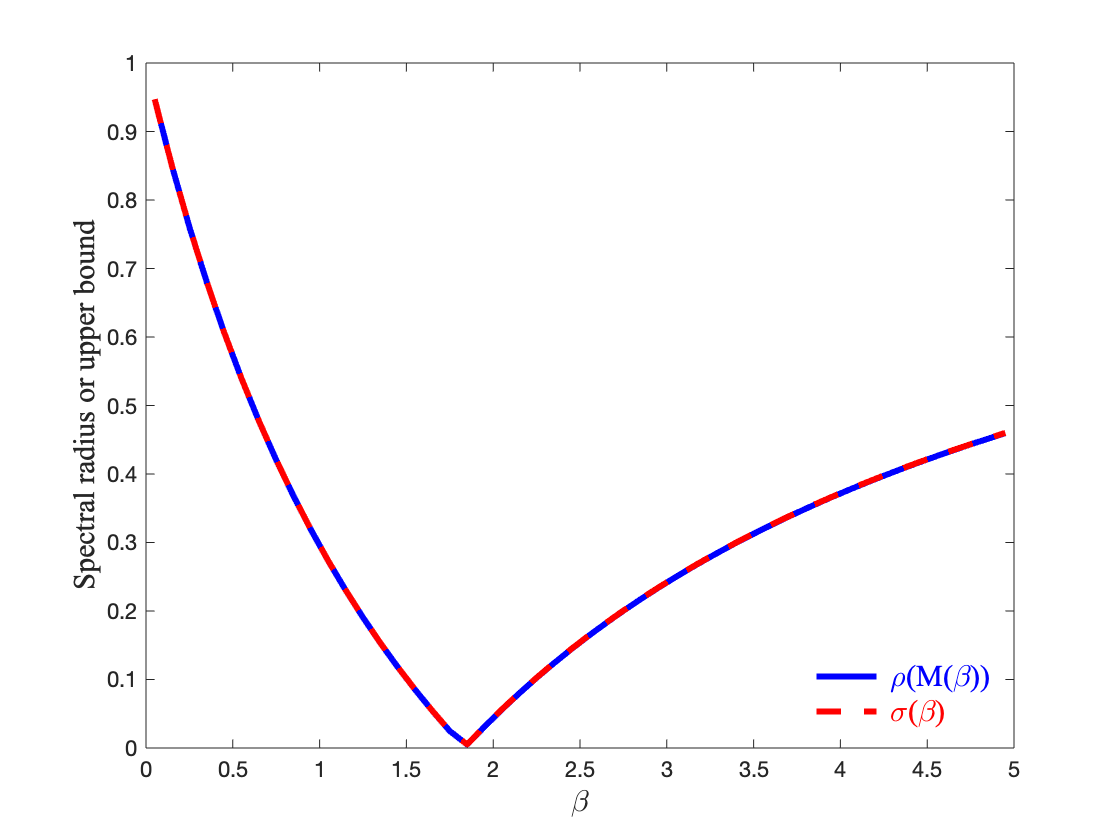}}
	\subfloat{\label{N_v2}\includegraphics[width=2in]{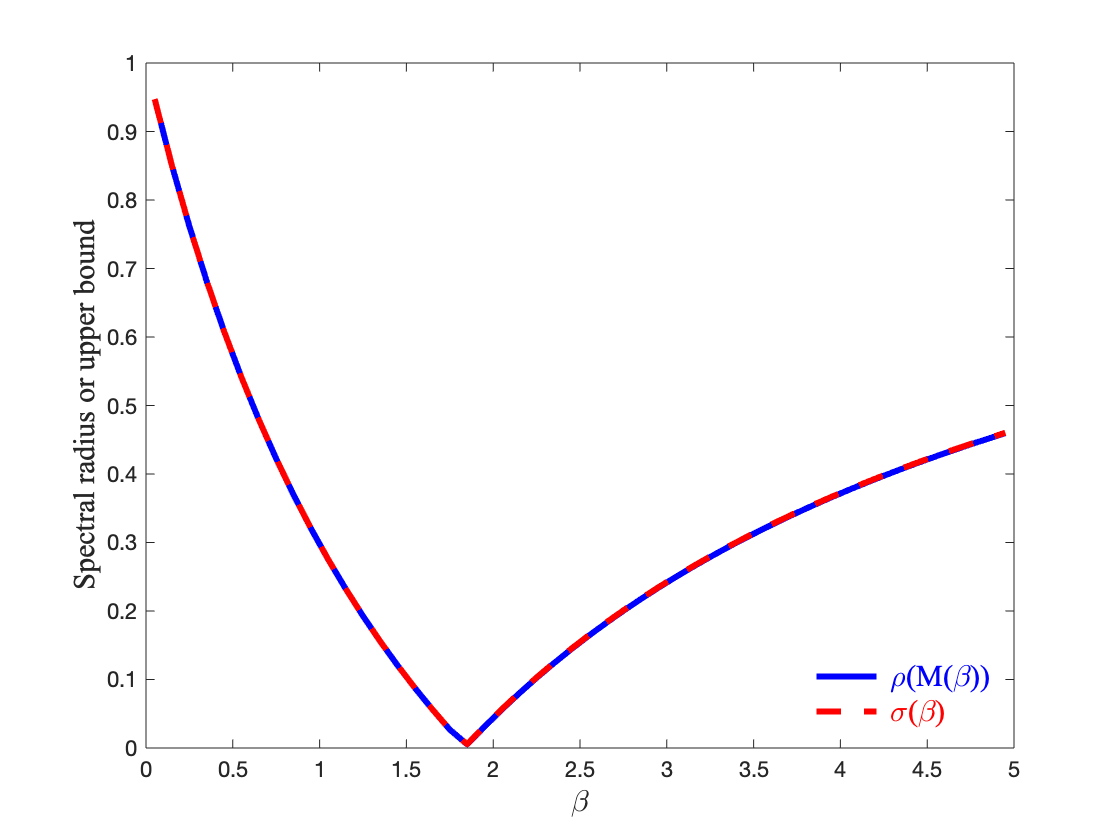}}
	\subfloat{\label{N_v3}\includegraphics[width=2in]{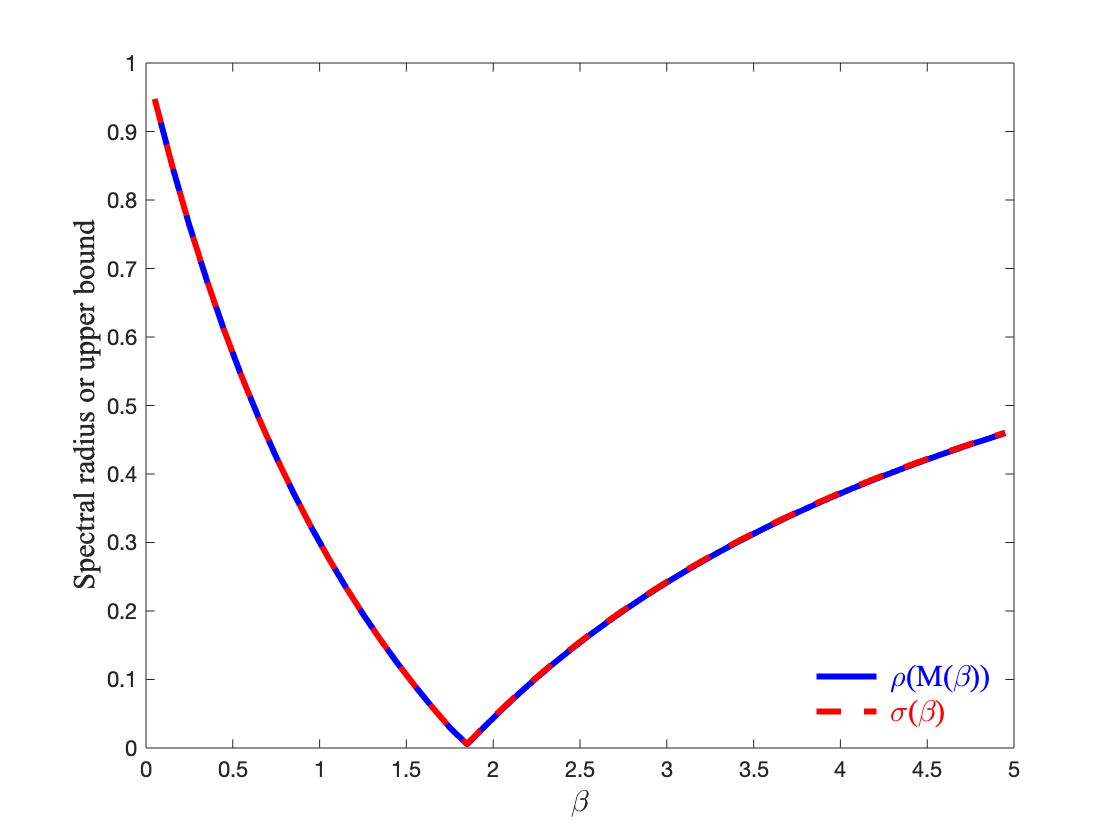}}
    \hspace{0.1in}
    \subfloat{\label{N_v1}\includegraphics[width=2in]{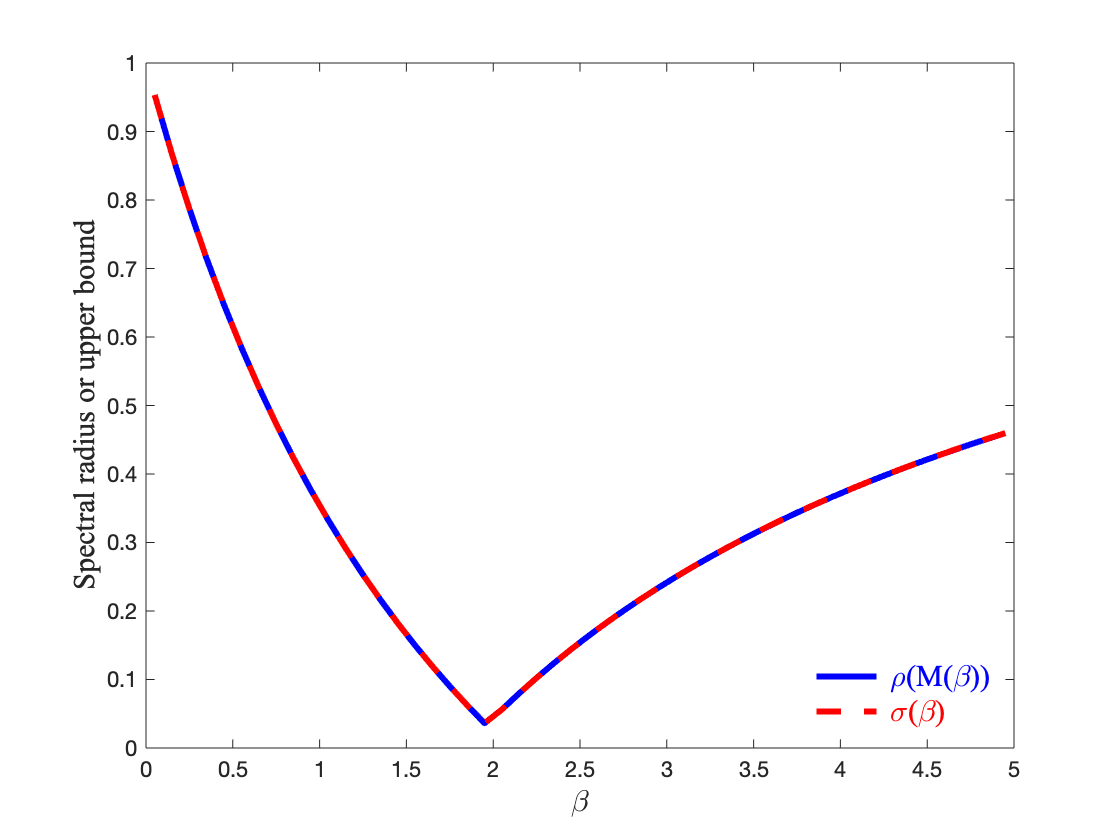}}
	\subfloat{\label{N_v2}\includegraphics[width=2in]{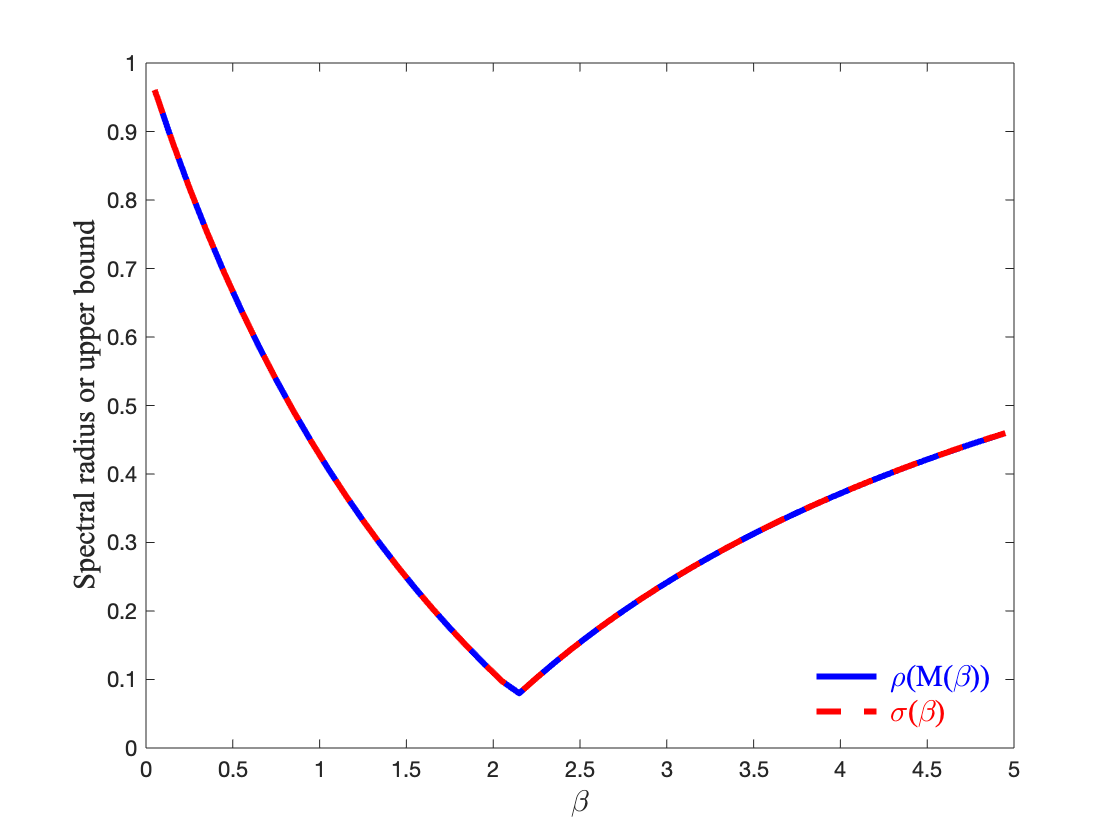}}
	\subfloat{\label{N_v3}\includegraphics[width=2in]{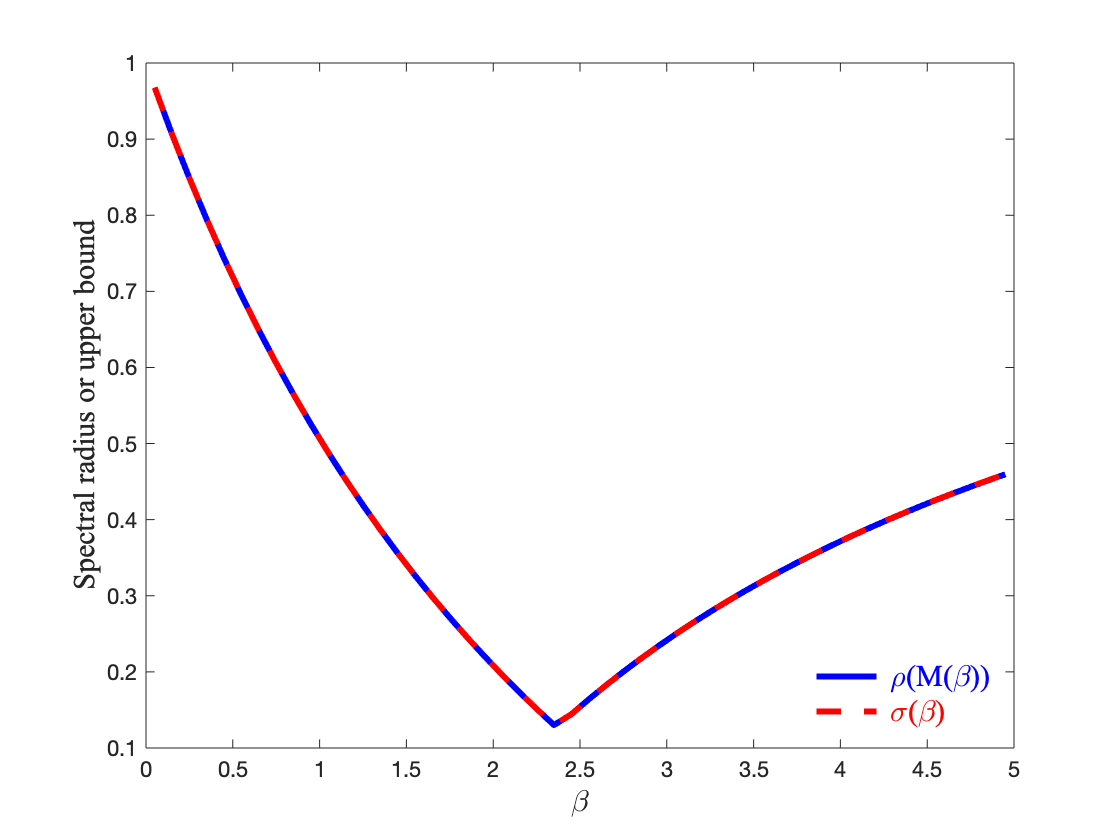}}
\caption{The spectral radius $\rho(M(\beta))$ of the iteration matrices for different $\beta$: ``$-$'' and the upper bound $\sigma(\beta)$ for different $\beta$: ``$--$'' with $T=0.1$ for scheme B in 3D. From left to right: $N_x=N_y=N_z=4,6,8$. From top to bottom: $\alpha=0,0.01,0.5$.}\label{fig:8}
\end{figure}
Figures~\ref{fig:5}--\ref{fig:8} further compare $\rho(M(\beta))$ with its upper bound $\sigma(\beta)$ for different values of $\beta$. In all tested cases, the upper bound closely follows the actual spectral radius, with very similar behavior observed for both schemes and in both one- and three-dimensional settings. These results indicate that the theoretical bound provides an accurate and relatively sharp estimate of $\rho(M(\beta))$ across the considered spatial resolutions and damping parameters.

\subsection{Results for IHSS iterations}\label{sec:ihss}

In our calculations, the inner Preconditioned Conjugate Gradient (PCG) and Aggregation‑based Algebraic MultiGrid (AGMG) iterates are terminated if the current residuals (tolerances) of the inner iterations satisfy
\begin{equation}
\varepsilon_\ell \leq \max\{0.1\delta^\ell, 1 \times 10^{-7}\}, \quad \eta_\ell \leq \max\{0.1\delta^\ell, 1 \times 10^{-6}\},
\end{equation}
with the controlled constant $\delta \in (0,1)$ need to be determined.

\subsubsection{The IHSS iteration method for scheme A and B}

We present the HSS iteration, for solving the system of linear equations $A\x=b$.

We use the HSS method for such a system, we then have
\begin{align*}
	H=\frac12 (A+A^{*})=D=\begin{pmatrix}
	\frac{11}{6}-k \alpha \epsilon \Delta_h&0&0\\
	0&\frac{11}{6}-k \alpha \epsilon \Delta_h&0\\
	0 & 0 & \frac{11}{6}-k \alpha \epsilon \Delta_h
	\end{pmatrix}
\end{align*}
and 
\begin{align*}
	S=\frac12 (A-A^{*})=\begin{pmatrix}
	0 & -k\epsilon \hat{w} \Delta_h &k\epsilon \hat{v} \Delta_h\\
	k\epsilon \hat{w} \Delta_h & 0 &-k \epsilon \hat{u} \Delta_h\\
	-k \epsilon \hat{v} \Delta_h & k\epsilon \hat{u} \Delta_h &0
	\end{pmatrix}
\end{align*}

\begin{table}[htbp]
    \centering
    \begin{tabular}{c|c|c|c|c}
    \hline
       Scheme & $N_x$&$N_t$& IHSS Average Iter & IHSS Average Relres \\
         \hline
       Scheme A &100 &100 &12 &4.190372514249884e-07 \\
       &500&100&53.142857142857146&9.376726839567723e-07\\
       &500&10&9.046250000000000e+02
&1.040222806415436e-06\\
&100&10&17.375000000000000&8.172751288649430e-07\\
         \hline
         Scheme B &100 &100 & 12 &4.190365609222554e-07 \\
         &500&100&40.387755102040813&9.355751951082065e-07\\
         &500&10&4.567500000000000e+02&9.947215325307641e-07\\
         &100&10&15&7.401248485752497e-07\\
         \hline
    \end{tabular}
    \begin{tabular}{c|c|c|c|c}
    \hline
       Scheme & $N_x$&$N_t$& IHSS Average Iter & IHSS Average Relres \\
         \hline
       Scheme A &100 &100 &1.208877551020408e+02&9.647518720100208e-07 \\
        &50&10&3.622500000000000e+02&9.899276620140013e-07\\
        &80&10&1001&9.592343749108075e-05\\
&100&10&1001&NaN\\
         \hline
         Scheme B &100 &100 & 77.438775510204081&9.577557668650449e-07 \\
         &50&10&2.236250000000000e+02&9.878883236020164e-07\\
         &80&10&4.721250000000000e+02&9.933686265851382e-07\\
         &100&10&8.066250000000000e+02&9.965410113137439e-07\\
         \hline
    \end{tabular}
    \caption{The estimated average of iteration and relative residual for scheme A and B in 1D. From top to bottom: $\alpha=0.01$ (top two), $\alpha=0.5$ (bottom two), $\beta=1$ and $\delta=0.8$.}
    \label{tab:p1}
\end{table}
Table~\ref{tab:p1} reports the average number of IHSS iterations
and the corresponding relative residuals for schemes A and B in 1D.
In most tested cases, the relative residual is reduced to approximately
$10^{-7}$, indicating that the prescribed stopping criterion is successfully
satisfied. However, the number of IHSS iterations increases significantly
for some parameter combinations, showing that the convergence of the inner
iteration becomes slower in these cases.

\begin{table}[htbp]
    \centering
    \begin{tabular}{c|c|c|c|c}
    \hline
       Scheme & $\beta$&$\delta$& IHSS Average Iter & IHSS Average Relres \\
         \hline
       Scheme A &1 &0.8 &12 &4.190372514249884e-07 \\
       &10&0.8&38&7.583533555835701e-07\\
       &1.833&0.8&2.051020408163265&7.025401232937075e-08\\
       &0.5&0.8&25&8.395218695759273e-07\\
         \hline
         Scheme B &1 &0.8 & 12 &4.190365609222554e-07 \\
         &10&0.8&38&7.583547136407340e-07\\
         &1.833&0.8&2.051020408163265&2.900336508329003e-08\\
         &0.5&0.8&25&8.395211632058914e-07\\
         \hline
    \end{tabular}
    \begin{tabular}{c|c|c|c|c}
    \hline
       Scheme & $\beta$&$\delta$& IHSS Average Iter & IHSS Average Relres \\
         \hline
       Scheme A &1 &0.8 &1.208877551020408e+02&9.647518720100208e-07 \\
       &10&0.8&38&7.583060228684486e-07\\
       &1.833&0.8&18.795918367346939&8.985023279541960e-07\\
       &0.5&0.8&3.222142857142857e+02
&9.814206421844759e-07\\
         \hline
         Scheme B &1 &0.8 &77.438775510204081 &9.577557668650449e-07 \\
         &10&0.8&38&7.583496321009146e-07\\
         &1.833&0.8&11.581632653061224&8.985479928942138e-07\\
         &0.5&0.8&2.473673469387755e+02 &9.812418035351897e-07\\
         \hline
    \end{tabular}
 \caption{The estimated average of iteration and relative residual for scheme A and B in 1D. From top to bottom: $\alpha=0.01$ (top two), $\alpha=0.5$ (bottom two), $N_x=100,N_t=100$.}
    \label{tab:p2}
\end{table}
Table~\ref{tab:p2} further shows that the average number of IHSS
iterations varies noticeably with the choice of $\beta$. This indicates that
the convergence behavior of the IHSS iteration is sensitive to the parameter
$\beta$.

\begin{table}[htbp]
    \centering
    \begin{tabular}{c|c|c|c|c}
    \hline
       Scheme & $N_x=N_y=N_z=N$&$N_t$& IHSS Average Iter & IHSS Average Relres \\
         \hline
       Scheme A &4 &40 &2&1.426085545565550e-07 \\
       &8&40&3.184210526315789&4.152191532052564e-07\\
       &16&40&2.736842105263158&1.686275419182831e-07\\
         \hline
         Scheme B &4 &40 &2&1.425241630814856e-07 \\
       &8&40&3.184210526315789&4.142937901376481e-07\\
       &16&40&2.736842105263158&1.663960123673145e-07\\
         \hline
    \end{tabular}
    \begin{tabular}{c|c|c|c|c}
    \hline
       Scheme & $N_x=N_y=N_z=N$&$N_t$& IHSS Average Iter & IHSS Average Relres \\
         \hline
       Scheme A &4 &40 &2.947368421052631&1.551262001884746e-07 \\
       &8&40&4.842105263157895&4.629733973736334e-07\\
       &16&40&4.289473684210527&3.929050307973152e-07\\
         \hline
          Scheme B &4 &40 &2.894736842105263&1.610855633609793e-07 \\
       &8&40&4.842105263157895&4.582489059923120e-07\\
       &16&40&4.289473684210527&3.480260626391152e-07\\
         \hline
    \end{tabular}
    \caption{The estimated average of iteration and relative residual for scheme A and B in 3D. From top to bottom: $\alpha=0.01$ (top two), $\alpha=0.5$ (bottom two), $\beta=1.833$ and $\delta=0.8$.}
    \label{tab:p3}
\end{table}
\noindent
Table~\ref{tab:p3} reports the average number of IHSS iterations
and the corresponding relative residuals for schemes A and B in 3D.
For all tested cases, the relative residual remains on the order of
$10^{-7}$, while the average number of iterations stays relatively small.
This indicates that the IHSS iteration maintains stable convergence behavior
for the considered three-dimensional test cases.

\begin{table}[htbp]
    \centering
    \begin{tabular}{c|c|c|c|c}
    \hline
       Scheme & $\beta$&$\delta$& IHSS Average Iter & IHSS Average Relres \\
         \hline
       Scheme A &1 &0.8 &12&4.191083671222556e-07 \\
       &10&0.8&38&7.579925106214455e-07\\
       &1.833&0.8&3.184210526315789&4.152191532052564e-07\\
       &0.5&0.8&25&8.386107912804774e-07\\
         \hline
         Scheme B &1 &0.8 &12 &4.191062544411093e-07 \\
         &10&0.8& 38&7.579943747665589e-07\\
         &1.833&0.8&3.184210526315789&4.142937901376481e-07 \\
         &0.5&0.8&25&8.386103294718192e-07\\
         \hline
    \end{tabular}
    \begin{tabular}{c|c|c|c|c}
    \hline
       Scheme & $\beta$&$\delta$& IHSS Average Iter & IHSS Average Relres \\
         \hline
       Scheme A &1 &0.8 &12&4.197122611007259e-07 \\
       &10&0.8&38&7.579073244529126e-07\\
       &1.833&0.8&4.842105263157895&4.629733973736334e-07\\
       &0.5&0.8&25&8.430547262289169e-07\\
         \hline
         Scheme B &1 &0.8 &12 &4.195493519525654e-07 \\
         &10&0.8&38&7.579885175392655e-07\\
         &1.833&0.8&4.842105263157895&4.582489059923120e-07\\
         &0.5&0.8&25&8.400858479968356e-07\\
         \hline
    \end{tabular}
 \caption{The estimated average of iteration and relative residual for scheme A and B in 3D. From top to bottom: $\alpha=0.01$ (top two), $\alpha=0.5$ (bottom two), $N_x=N_y=N_z=8,N_t=40$.}
    \label{tab:p4}
\end{table}
Table~\ref{tab:p4} reports the average number of IHSS iterations
and the corresponding relative residuals for schemes A and B in 3D.
The two schemes exhibit very similar iteration counts for the tested values
of $\beta$, while the relative residuals remain on the order of $10^{-7}$.
This indicates that both schemes show comparable and stable convergence
behavior for the considered three-dimensional test cases.

\begin{figure}[htbp]
    \centering
    \subfloat[arrow, scheme A]{\includegraphics[width=0.3\linewidth]{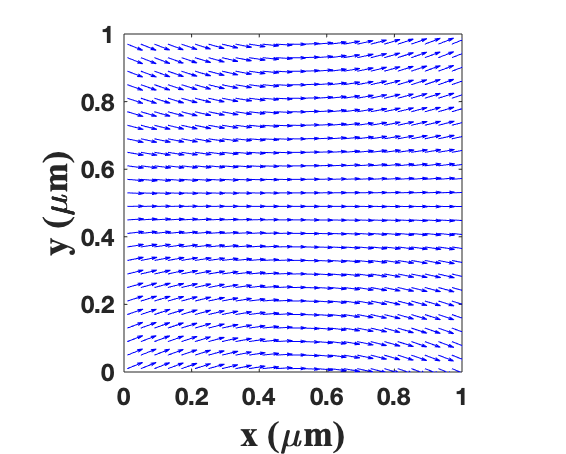}}
    \subfloat[color, scheme A]{\includegraphics[width=0.3\linewidth]{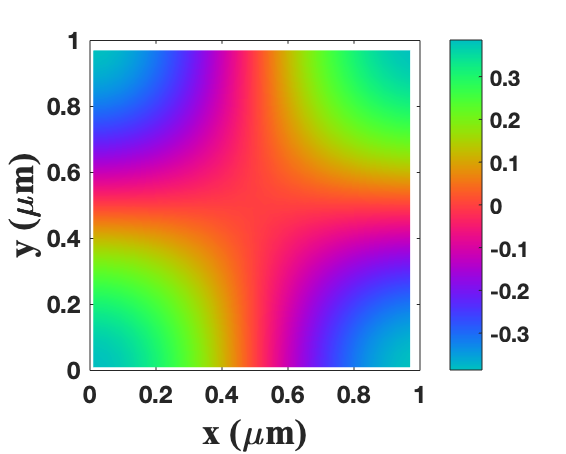}}
    \hspace{0.1in}
    \subfloat[arrow, scheme B]{\includegraphics[width=0.3\linewidth]{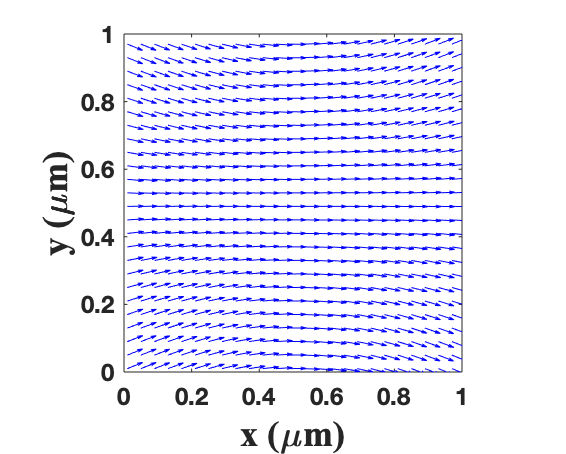}}
    \subfloat[color, scheme B]{\includegraphics[width=0.3\linewidth]{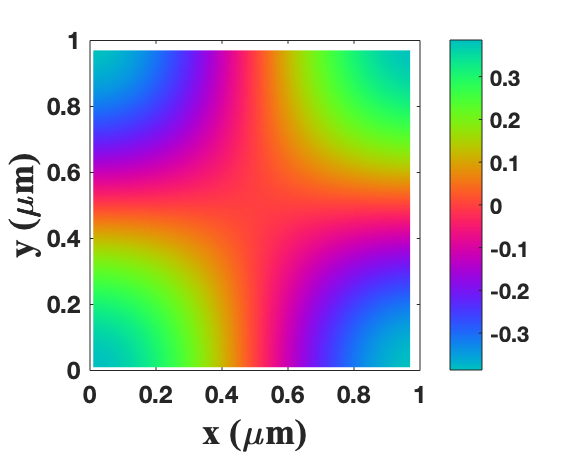}}
    \caption{The solution profiles for scheme A and B in 3D given the initial condition $\m_0=[\cos(\cos(\pi x)\cos(\pi y)\cos(\pi z))\sin(0),\sin(\cos(\pi x)\cos(\pi y)\cos(\pi z))\sin(0),\cos(0)]^T$ with $N_x=N_y=50$, $N_z=4$ at the final time $T=0.1$. $\beta=1.833$ and $\delta=0.8$, $\alpha=0.5$.}
    \label{fig:p5}
\end{figure}
As shown in Figure~\ref{fig:p5}, schemes A and B produce
very similar solution profiles at the final time. Both the arrow plots and
the color maps exhibit nearly identical spatial distributions, indicating
that the two schemes provide consistent numerical solutions for the
considered three-dimensional test problem.

\FloatBarrier

\section{Conclusions}
\label{sec:conclusions}

In this work, we consider two third-order semi-implicit schemes for solving the Landau-Lifshitz-Gilbert equation. The resulting linear systems inherit a large-sparse structure and are non-Hermitian with all real eigenvalues being positive. Such linear systems are well-suited for iterative solution via the HSS/IHSS methods, which enjoy well-established convergence theory and favourable convergence rates. Numerical experiments are presented to validate the effectiveness of these approaches for micromagnetic simulations.

\section*{Acknowledgments}
This work is supported in part by the Jiangsu Science and Technology Programme-Fundamental Research Plan Fund (BK20250468), Research and
 Development Fund of XJTLU (RDF-24-01-015).

\bibliographystyle{amsplain}
\bibliography{references}

\end{document}